\documentclass[10pt,a4paper,reqno,english]{amsart} 
\usepackage{pstricks,amsmath,here,latexsym,amssymb}
\usepackage{epsf,psfrag,epsfig,color,graphicx}
\usepackage[latin1]{inputenc}

\newtheorem{prop}{Proposition} 
\newtheorem{lemma}[prop]{Lemma}
\newtheorem{cor}[prop]{Corollary} 

\newtheorem{thm}[prop]{Theorem}
\newcommand{\dem}{\noindent \textbf{Proof. }}

\theoremstyle{definition} 
\newtheorem{Def}[prop]{Definition}

\newcommand{\Ref}[1]{(\ref{#1})}
\newcommand{\findem}{\vspace{-.4cm} \begin{flushright} $\square~$ \end{flushright} \vspace{.4cm} }

\newcommand{\ite}{\noindent $\bullet~$}
\newcommand{\iten}{\noindent -~}
\newcommand{\titre}[1]{\noindent \textbf{#1}}

\newcommand{\B}[1]{\overline{#1}}

\newcommand{\mG}{\mathcal{G}}
\newcommand{\mI}{\mathcal{I}}
\newcommand{\mE}{\mathcal{E}}

\newcommand{\mO}{\mathcal{O}}
\newcommand{\mS}{\mathcal{S}}

\newcommand{\fir}{\textrm{first}}
\newcommand{\pos}{\textrm{post}}

\newcommand{\connect}[1]{{#1 \atop \rightarrow}}
\newcommand{\topple}[1]{{#1 \atop  \dashrightarrow}}

\catcode`\@=11
\def\section{\@startsection{section}{1}%
 \z@{1.5\linespacing\@plus\linespacing}{.5\linespacing}%
 {\normalfont\bfseries\centering}}

\def\subsection{\@startsection{subsection}{2}%
  \z@{.6\linespacing\@plus\linespacing}{.5\linespacing}%
  {\normalfont\bfseries}}

\def\subsubsection{\@startsection{subsubsection}{3}%
  \z@{.5\linespacing\@plus.7\linespacing}{-.5em}%
  {\normalfont\itshape}}
\catcode`\@=12
\title[Tutte polynomial, orientations and Sandpile model]
{Tutte polynomial, subgraphs, orientations and sandpile model: new connections via embeddings}
\author{Olivier Bernardi}
\address{LaBRI, Université Bordeaux 1, 351 cours de la Libération, 33405 Talence, France}
\email{bernardi@labri.fr}

\date{\today}

\begin{document}
\maketitle

\begin{abstract}
For any graph $G$ with $n$ edges, the spanning subgraphs and the orientations of $G$ are both counted by the evaluation $T_G(2,2)=2^n$ of its Tutte polynomial. We define a bijection $\Phi$ between spanning subgraphs and orientations and explore its enumerative consequences regarding the Tutte polynomial. The bijection $\Phi$ is closely related to a recent characterization of the Tutte polynomial relying on a \emph{combinatorial embedding} of the graph $G$, that is, on a choice of cyclic order of the edges around each vertex. Among other results, we obtain a combinatorial interpretation for each of the evaluations $T_G(i,j),0\leq i,j \leq 2$ of the Tutte polynomial in terms of orientations. The strength of our approach is to derive all these interpretations by specializing the bijection $\Phi$ in various ways. For instance, we obtain a bijection between  the connected subgraphs of~$G$ (counted by $T_G(1,2)$) and the root-connected orientations. We also obtain  a bijection between  the forests (counted by $T_G(2,1)$) and outdegree sequences which specializes into a bijection between spanning trees (counted by $T_G(1,1)$) and root-connected outdegree sequences. We also define a bijection between spanning trees and recurrent configurations of the sandpile model. Combining our results we obtain a bijection between recurrent configurations and root-connected outdegree sequences which leaves the configurations at level 0 unchanged.
\end{abstract}

\section{INTRODUCTION}

In 1947, Tutte defined a graph invariant that he named the \emph{dichromate} because he thought of it as bivariate generalization of the chromatic polynomial \cite{Tutte:dichromate-ring}. Since then, the dichromate, now known as the \emph{Tutte polynomial}, has been widely studied (see \cite{Bollobas:Tutte-poly-chapter10,Brylawski:Tutte-poly}).\\


There are several alternative definitions of the Tutte polynomial \cite{OB:Tutte-plongement-def,Gessel:Tutte-poly+DFS,Vergnas:Morphism-matroids-2,Tutte:dichromate}. The most straightforward definition for a connected graph $G=(V,E)$ is 
\begin{eqnarray}
T_G(x,y)=\sum_{S\textrm{ spanning subgraph}} (x-1)^{c(S)-1}(y-1)^{c(S)+|S|-|V|}, \label{eq:intro-Tutte-subgraph}
\end{eqnarray}
where the sum is over all spanning subgraphs $S$ (equivalently, subsets of edges), $c(S)$ denotes the number of connected components of $S$ and $|.|$ denotes cardinality. From this definition, it is easy to see that $T_G(1,1)$ (resp. $T_G(2,1)$,  $T_G(1,2)$) counts the spanning trees (resp. forests, connected subgraphs) of $G$. A somewhat less interesting specialization is $T_G(2,2)=2^{|E|}$ counting the spanning subgraphs of $G$. Note that this is also the number of orientations of $G$. As a matter of fact, all the specializations $T_G(i,j),0\leq i,j\leq 2$ as well as  some of their refinements have nice interpretations in terms of orientations \cite{Brylawski:Tutte-poly,Gessel:Tutte-poly+DFS, Gioan:enumerating-degree-sequences, Greene:interpretation-Tutte-poly,  Vergnas:Morphism-matroids-2, Lass:interpretation-Tutte-poly, Stanley:acyclic-orientations}. \\

As one can see, there is a lot of interesting specializations of the Tutte polynomial and a number of articles are devoted to combinatorial proofs of these specializations \cite{Gebhard-sagan:acyclic-orientations, Gessel:Tutte-poly+DFS, Gessel:enumerative-csq-DFS,Gioan:enumerating-degree-sequences, Gioan-bij-tree-orientation, Lass:interpretation-Tutte-poly}. In this paper, we give bijective proofs for the interpretation of each of the evaluations $T_G(i,j),0\leq i,j\leq 2$ in terms of orientations. The strength of our approach is to derive all these interpretations from a single bijection between subgraphs and orientations that we specialize in various ways. For instance, we derive a bijection between connected subgraphs (counted by $T_G(1,2)$) and root-connected orientations. We also derive a bijection between forests (counted by $T_G(2,1)$) and outdegree sequences.   In particular, we derive a bijection between spanning trees (counted by $T_G(1,1)$) and root-connected outdegree sequences.  \\

We shall also deal with the \emph{sandpile model} \cite{Bak:sandpile,Dhar:sandpile} (equivalently \emph{chip firing game} \cite{Bjorner:chip-firing}). It is known that the recurrent configurations of the sandpile model on $G$ (equivalently  $G$-parking functions \cite{Shapiro:G-parking-function}) are counted by $T_G(1,1)$ \cite{Dhar:sandpile}. Observe that this is the number of spanning trees. The following refinement is also true: the coefficient of $y^k$ in $T_G(1,y)$ is the number of recurrent configurations at \emph{level} $k$ \cite{Merino:external-activity=sandpile-level}. A bijective proof of this result was given in \cite{Borgne-cori-activite-externe-sable}. We give an alternative bijective proof. We also answer a question of Gioan  \cite{Gioan:enumerating-degree-sequences} by establishing a  bijection between recurrent configurations of the sandpile model and root-connected outdegree sequences that leaves the configurations at level 0 unchanged.\\

Our bijections require a choice of a combinatorial embedding of the graph $G$, that is, a choice of a cyclic ordering of the edges around each vertex. In \cite{OB:Tutte-plongement-def}  the \emph{internal} and \emph{external} \emph{embedding-activities} of spanning trees  were defined for embedded graphs. It was proved that for any embedding of the graph $G$, the Tutte polynomial of $G$ is given by 
\begin{eqnarray}
 T_G(x,y)=\sum_{T\textrm{ spanning tree}} x^{\mI(T)}y^{\mE(T)}, \label{eq:intro-Tutte-OB}
\end{eqnarray}
where the sum is over all spanning trees $T$ and $\mI(T)$ (resp. $\mE(T)$) denotes the  internal (resp. external) embedding-activity. This characterization of the Tutte polynomial is reminiscent but inequivalent to the one given by Tutte in \cite{Tutte:dichromate}. The characterization \Ref{eq:intro-Tutte-OB} is our main tool in order to obtain enumerative corollary from our bijections. In this respect, our approach is close to the one used by Gessel and Sagan in \cite{Gessel:enumerative-csq-DFS,Gessel:Tutte-poly+DFS} in order to obtain enumerative consequences from a new notion of \emph{external activity}. \\

The outline of this paper is as follows. \\
\ite In Section \ref{section:definition}, we recall some definitions and preliminary results obout graphs, orientations and the sandpile model.  \\
\ite In Section \ref{section:glimpse}, we take a glimpse at the results to be developed in the following sections. We first establish some elementary results about the \emph{tour} of spanning trees and their \emph{embedding-activities}. 
Then we define a mapping $\Phi$ from spanning trees to orientations. We highlight a connection between the embedding-activities of a spanning tree $T$ and the acyclicity or strong-connectivity of the orientation $\Phi(T)$. Building on the mapping $\Phi$ we also define a bijection $\Gamma$ between spanning trees to root-connected outdegree sequences and a closely related bijection $\Lambda$ between spanning trees and recurrent configurations of the sandpile model.  \\
\ite In Section  \ref{section:partition}, we define a partition $\Pi$ of the set of subgraphs. Each part of this partition is  an interval in the boolean lattice of the  set of subgraphs and is associated to a spanning tree. The interval associated with a spanning tree $T$ is closely related to the embedding-activities of $T$. 
We show how the partition  $\Pi$ explains the link between the subgraph expansion \Ref{eq:intro-Tutte-subgraph} and the spanning tree expansion \Ref{eq:intro-Tutte-OB} of the Tutte polynomial. We also consider several criteria for subgraphs: connected, forest, \emph{internal}, \emph{external} and prove that the families of subgraphs that can be defined by combining these criteria are counted by one of the evaluations  $T_G(i,j),0\leq i,j\leq 2$ of the Tutte polynomial.\\
\ite In Section \ref{section:subgraphs-orientations}, we extend the mapping $\Phi$ to the set of all subgraphs. This definition makes use of the partition  $\Pi$ of the set of subgraphs. We prove that $\Phi$ is a bijection between subgraphs and orientations.\\
\ite In Section \ref{section:specializations}, we study the specializations of the bijection $\Phi$ to the families of subgraphs  defined by the criteria \emph{connected, forest, internal, external}. We prove that $\Phi$ induces bijections between these families of subgraphs and the families of orientations defined by the criteria \emph{root-connected, minimal, acyclic, strongly connected}. As a consequence, we obtain an interpretation for each of the evaluations $T_G(i,j)$, $0\leq i,j \leq 2$ of the Tutte polynomial in terms of orientations or outdegree sequences.\\
\ite In Section \ref{section:bij-sandpile}, we study the bijection  $\Lambda$ between spanning trees and recurrent configurations of the sandpile model. \\
\ite Lastly, in Section \ref{section:conclusion} we comment on the case of planar graphs.\\

\section{DEFINITIONS} \label{section:definition}
We denote by $\mathbb{N}$ the set of non-negative integers. For any set $S$, we denote by $|S|$ its cardinality. For any sets $S_1,S_2$, we denote by $S_1\vartriangle S_2$ the symmetric difference of  $S_1$ and $S_2$. If $S\subseteq S'$ and $S'$ is clear from the context, we denote by $\B{S}$ the complement of $S$, that is, $S'\setminus S$.
If $S\subseteq S'$ and $s\in S'$, we write $S+s$ and $S-s$ for  $S\cup\{s\}$ and $S\setminus\{s\}$ respectively (whether $s$ belongs to $S$ or not).

\subsection{Graphs}
\begin{figure}[htb!]
\begin{center}
\input{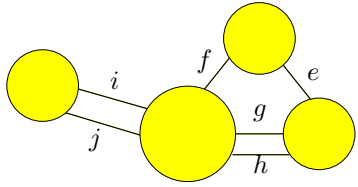}
\caption{The cut $\{e,f,g,h,i,j\}$ and the connected components after deletion of this cut (shaded regions).}\label{fig:cut-cocycle}
\end{center}
\end{figure}

In this paper we consider finite, undirected graphs. Loops and multiple edges are allowed but, for simplicity, we shall only consider \emph{connected} graphs. Formally, a \emph{graph} $G=(V,E)$ is a finite set of \emph{vertices} $V$, a finite set of \emph{edges} $E$ and a relation of \emph{incidence} in $V\times E$ such that each edge $e$ is incident to either one  or two vertices. The \emph{endpoints} of an edge $e$ are the vertices incident to $e$. A \emph{cycle} is a set of edges that form a simple closed path. A \emph{cut} is a set of edges $C$ whose deletion increases the number of connected components and such that the endpoints of every edge in $C$ are in distinct components of the resulting graph. A cut is shown in Figure \ref{fig:cut-cocycle}. Given a subset of vertices $U$, the  \emph{cut defined  by $U$} is the set of edges with one endpoint in $U$ and one endpoint in $\B{U}$.  A \emph{cocycle} is a cut which is minimal for inclusion (equivalently it is a cut whose deletion increases the number of connected components by \emph{one}). For instance, the set of edges $\{f,g,h\}$ in Figure \ref{fig:cut-cocycle} is a cocycle.\\

Let $G=(V,E)$ be a graph. A \emph{spanning subgraph} of $G$ is a graph  $G'=(V,E')$ where $E'\subseteq E$. \emph{All the subgraphs considered in this paper are spanning} and we shall not further precise it. A subgraph is entirely determined by its edge set and, by convenience, we shall identify the subgraph with its edge set.  A forest is an acyclic graph. A \emph{tree} is a connected forest. A \emph{spanning tree} is a (spanning) subgraph which is a tree. Given a tree $T$ and a vertex distinguished as the \emph{root-vertex} we shall use the usual \emph{family vocabulary} and talk about the \emph{father}, \emph{son}, \emph{ancestors} and \emph{descendants} of vertices in $T$. By convention, a vertex is considered to be an ancestor and a descendant of itself.
If a vertex of the graph $G$ is distinguished as the \emph{root-vertex} we implicitly consider it to be the root-vertex of every spanning tree.\\

Let  $T$  be a spanning tree of the graph $G$. An edge of $G$ is said to be \emph{internal} if it is in $T$ and \emph{external} otherwise. The \emph{fundamental cycle} (resp. \emph{cocycle}) of an external (resp. internal) edge $e$ is the set of edges $e'$ such that the
 subgraph $T-e'+e$ (resp.  $T-e+e'$) is a spanning tree. Observe that the fundamental cycle $C$ of an external edge $e$ is a cycle contained in $T+e$ ($C$ is made of $e$ and the path of $T$ between the endpoints of $e$). Similarly,  the fundamental cocycle $D$ of an internal edge $e$ is a cocycle contained in $\B{T}+e$ ($D$ is made of the edges linking the two subtrees obtained from $T$ by removing $e$). Observe also that, if $e$ is internal and $e'$ is external, then $e$ is in the fundamental cycle of $e'$ if and only if $e'$ is in the fundamental cocycle of $e$.

\subsection{Embeddings} \label{section:embedding+tour}
We recall the notion of \emph{combinatorial map} \cite{Cori:These-asterisque,Cori-Machi:survey}.  A \emph{combinatorial map} (or \emph{map} for short) $\mG=(H,\sigma,\alpha)$ is a set of \emph{half-edges} $H$, a permutation $\sigma$  and an involution without fixed point $\alpha$ on $H$ such that the group generated by $\sigma$ and $\alpha$ acts transitively on $H$. A map is \emph{rooted} if one of the half-edges is distinguished as the \emph{root}. For $h_0\in H$, we denote by $\mG=(H,\sigma,\alpha,h_0)$ the map $(H,\sigma,\alpha)$ rooted on $h_0$. From now on \emph{all our maps are rooted}. \\

\begin{figure}[htb!]
\begin{center}
\input{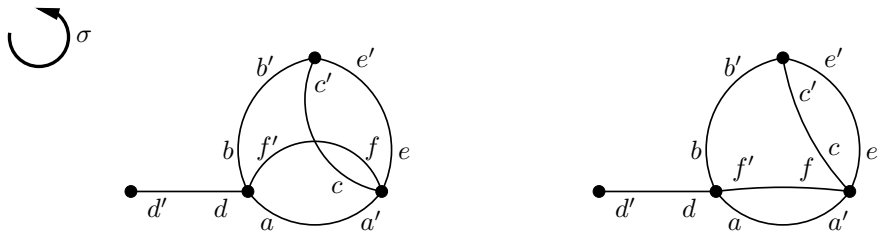}
\caption{Two embeddings of the same graph.}\label{fig:exp-embedding-connection}
\end{center}
\end{figure}

Given a  map $\mG=(H,\sigma,\alpha,h_0)$, we consider the \emph{underlying} graph $G=(V,E)$, where $V$ is the set of cycles of $\sigma$, $E$ is the set of cycles of $\alpha$ and the incidence relation is to have at least one common half-edge.  We represented the underlying graph of the map $\mG=(H,\sigma,\alpha)$  on the left of Figure \ref{fig:exp-embedding-connection}, where the set of half-edges is $H=\{a,a',b,b',c,c',d,d',e,e',f,f'\}$, the involution $\alpha$  is $(a,a')(b,b')(c,c')(d,d')(e,e')(f,f')$ in cyclic notation and the permutation $\sigma$ is $(a,f',b,d)(d')(a',e,f,c)(e',b',c')$. 
Graphically, we keep track of the cycles of $\sigma$ by drawing the half-edges of each cycle in counterclockwise order around the corresponding vertex. Hence, our drawing characterizes the map $\mG$ since the order around vertices give the cycles of the permutation $\sigma$ and the edges give the cycles of the involution $\alpha$. On the right of Figure \ref{fig:exp-embedding-connection},  we represented the map $\mG'=(H,\sigma',\alpha)$, where $\sigma'=(a,f',b,d)(d')(a',e,c,f)(e',b',c')$. The maps  $\mG$ and $\mG'$ have isomorphic underlying graphs. \\

Note that the underlying graph of a map $\mG=(H,\sigma,\alpha)$ is always connected since  $\sigma$ and $\alpha$ act transitively on $H$. A  \emph{combinatorial embedding} (or \emph{embedding} for short) of a connected graph $G$ is a   map $\mG=(H,\sigma,\alpha)$ whose underlying graph is isomorphic to $G$ (together with an explicit bijection between the set $H$ and the set of half-edges of $G$).  When an embedding $\mG$  of $G$ is given we shall write the edges of $G$  as pairs of half-edges (writing for instance $e=\{h,h'\}$). Moreover, we call \emph{root-vertex} the vertex incident to the root and \emph{root-edge} the edge containing the root.  In the following, we use the terms \emph{combinatorial map} and \emph{embedded graph} interchangeably. \emph{We do not require our graphs to be planar}.\\

Intuitively, a combinatorial embedding corresponds to the choice of a cyclic order on the edges around each vertex. This order  can also be seen as a local planar embedding. In fact there is a one-to-one correspondence between  \emph{combinatorial embeddings of graphs} and  the \emph{cellular embeddings of graphs in orientable surfaces} (defined up to homeomorphism); see \cite[Thm. 3.2.4]{Mohar:graphs-on-surfaces}.

\subsection{Orientations and outdegree sequences} 
\begin{figure}[ht!]
\begin{center}
\input{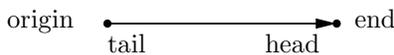}
\caption{Half-edges and endpoints of arcs.} \label{fig:head-and-tail-Tutte}
\end{center}
\end{figure}

Let $G$ be a graph and let $\mG$ be an embedding of $G$. An orientation is a choice of a direction for each edge of $G$, that is to say,  a function $\mO$ which associates to any edge $e=\{h_1,h_2\}$ one of the ordered pairs $(h_1,h_2)$ or $(h_1,h_2)$. Note that loops have two possible directions. We call $\mO(e)$ an \emph{arc}, or \emph{oriented edge}. If $\mO(e)=(h_1,h_2)$ we call $h_1$ the \emph{tail} and $h_2$ the \emph{head}. We call \emph{origin} and \emph{end} of  $\mO(e)$ the endpoint of the tail and head respectively. Graphically, we represent an  arc by an arrow going from the origin to the end (see Figure~\ref{fig:head-and-tail-Tutte}).

A \emph{directed path} is a sequence of arcs $(a_1,a_2,\ldots,a_k)$ such that the end of $a_i$ is the origin of $a_{i+1}$ for $1\leq i \leq k-1$. A \emph{directed cycle} is a simple directed closed path.
A \emph{directed cocycle} is a set of arcs $a_1,\ldots,a_k$  whose deletion disconnects the graph into two components and such that all arcs are directed toward the same component. If the orientation $\mO$ is not clear from the context, we shall say that a path, cycle, or cocycle is \emph{$\mO$-directed}. An orientation is said to be \emph{acyclic} (resp. \emph{totally cyclic} or \emph{strongly connected}) if there is no directed cycle (resp. cocycle).\\

We say that a vertex $v$ is  \emph{reachable} from a vertex $u$ if there is a directed path  $(a_1,a_2,\ldots,a_k)$ such that the origin of $a_1$ is $u$ and the end of $a_k$ is $v$. If  $v$ is  \emph{reachable} from  $u$ in the orientation $\mO$ denote it by  $u\connect{\mO} v$. An orientation is said to be \emph{$u$-connected} if every vertex is reachable from $u$. It is known that  every edge in an oriented graph is either in a directed cycle but not both \cite{Minty:cycle-cocycle}. Hence, an orientation $\mO$ is strongly connected if and only if the origin of every arc is reachable from its end. Equivalently, $\mO$ is strongly connected if every pair of vertices are reachable from one another.\\ 


The \emph{outdegree sequence} (or \emph{score vector}) of an orientation $\mO$ of the graph $G=(V,E)$ is the function $\delta: V\mapsto \mathbb{N}$ which associates to every vertex the number of incident tails. We say that $\mO$ is a \emph{$\delta$-orientation}. The outdegree sequences are strongly related to the \emph{cycle flips}, that is, the reversing of every edge in a directed cycle. Indeed, it is known that the outdegree sequences are in one-to-one correspondence with the equivalence classes of orientations up to \emph{cycle flips}\cite{Felsner:lattice}.\\

\begin{figure}[ht!]
\begin{center}
\input{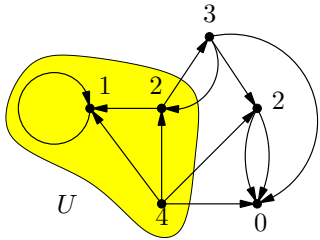}
\caption{The excess of the subset $U$ is $exc_\delta(U)=(4+2+1)-4=3$.} \label{fig:exp-excess}
\end{center}
\end{figure}

There are nice characterizations of the functions  $\delta: V\mapsto \mathbb{N}$ which are the outdegree sequence of an orientation. Given a function $\delta: V\mapsto \mathbb{N}$, we define the \emph{excess} of a subset of vertices $U\subseteq V$ by
$$exc_\delta(U)=\left(\sum_{u \in U}\delta(u)\right) - |G_U|,$$
where $|G_U|$ is the number of edges of $G$ having both endpoints in $U$. By definition, if $\delta$ is the outdegree sequence of an orientation $\mO$, the sum  $\sum_{u \in U}\delta(u)$ is the number of tails incident with vertices in $U$. From this number, exactly $|G_U|$ are part of edges with both endpoints in $U$. Hence, the excess $exc_\delta(U)$ corresponds to the number of tails incident with vertices in $U$ in the cut defined by $U$. This property is illustrated in Figure \ref{fig:exp-excess}. It is clear that if $\delta: V\mapsto \mathbb{N}$ is an outdegree sequence, then the excess of $V$ is 0 and the excess of any subset $U\subseteq V$ is non-negative. In fact, the converse is also true: every function $\delta: V\mapsto \mathbb{N}$ satisfying these two conditions is an outdegree sequence \cite{Felsner:lattice}.\\

We now prove that the reachability between two vertices in a directed graph only depends on the outdegree sequence of the orientation.

\begin{lemma}\label{thm:accessibility-outsequence}
Let $G=(V,E)$ be a graph and let $u,~v$ be two vertices. Let $\mO$ be an orientation of $G$ and let $\delta$ be its outdegree sequence. Then $v$ is reachable from  $u$ if and only if there is no subset of vertices $U\subseteq V$ containing $u$ and not $v$ and such that $exc_\delta(U)=0$.
\end{lemma}

\begin{figure}[ht!]
\begin{center}
\input{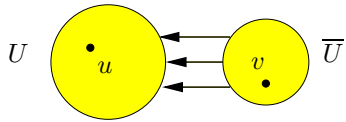}
\caption{Reachability is a property of the outdegree sequence.} \label{fig:accessibility-outsequence}
\end{center}
\end{figure}

\dem Lemma \ref{thm:accessibility-outsequence} is illustrated in Figure \ref{fig:accessibility-outsequence}. Observe that the excess of a subset $U\subseteq V$ is 0 if and only if the cut defined by $U$ is directed toward $U$. \\
\ite Suppose  there is a subset of vertices $U\subseteq V$ containing $u$ and not $v$ such that $exc_\delta(U)=0$. Then, the cut defined by $U$ is directed toward $U$. Thus, there is no directed path from $U$ to $\B{U}$. Hence $v$ is not reachable  from $v$.\\
\ite  Conversely, suppose   $v$ is not reachable from $u$. Consider the set of vertices $U$ reachable from $u$. The subset $U$ contains $u$ but not $v$. Moreover, the cut defined by $U$ is directed toward $U$, hence $exc_\delta(U)=0$.
\findem

Since reachability only depends on the outdegree sequence of the orientation, we can define an outdegree sequence $\delta$ to be \emph{$u$-connected} or \emph{strongly connected} if the $\delta$-orientations are. The $u$-connected outdegree sequences were considered in  \cite{Gioan:enumerating-degree-sequences} in connection with the \emph{cycle/cocycle reversing system} (see Subsection \ref{section:cycle-reversing-system}).\\

\noindent \textbf{Remark:} From the characterization of outdegree sequences given above and Lemma~\ref{thm:accessibility-outsequence} it is possible to  characterize  \emph{$u$-connected} and \emph{strongly connected} outdegree sequences. Let $G=(V,E)$ be a graph and $\delta: V\mapsto \mathbb{N}$ be a mapping such that $\sum_{v \in V}\delta(v) = |E|$.  The mapping $\delta$ is a strongly connected outdegree sequence if and only if the excess of any subset $U\subsetneq V$ is positive. The mapping $\delta$ is a $u$-connected outdegree sequence if and only if the excess of any subset $U\subsetneq V$ is non-negative and is positive whenever $u\in U$.

\subsection{The sandpile model} \label{section:def-sandpile}
The sandpile model is a dynamical system introduced in statistical physics in order to study self-organized criticality \cite{Bak:sandpile,Dhar:sandpile-self-organized}. It appeared independently in combinatorics as the \emph{chip firing game} \cite{Bjorner:chip-firing}. Roughly speaking, the model consists of grains of sand toppling through edges when there are too many on the same vertex. \emph{Recurrent configurations} play an important role in the model: they correspond to configurations that can be observed after a long period of time. The recurrent configuration are also equivalent to the $G$-parking functions introduced by Shapiro and Postnikov in the study of certain quotient of the polynomial ring \cite{Shapiro:G-parking-function}. Despite its simplicity, the sandpile model displays interesting enumerative \cite{Borgne-cori-activite-externe-sable,Dhar:sandpile,Merino:external-activity=sandpile-level} and algebraic properties \cite{Cori:dual-sandpile,Dhar:sandpile-algebraic}.\\

\begin{figure}[ht!]
\begin{center}
\input{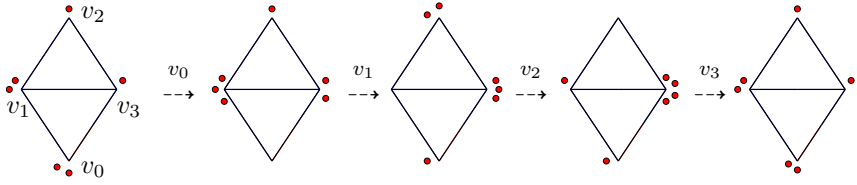}
\caption{A recurrent configuration and the evolution rule.} \label{fig:sable-diamant}
\end{center}
\end{figure}

Let $G=(V,E)$ be a graph with a vertex $v_0$ distinguished as the \emph{root-vertex}. A \emph{configuration of the sandpile model} (or \emph{sandpile configuration} for short) is a function $\mS:V\mapsto \mathbb{N}$, where $\mS(v)$ represents the number of grains of sand on $v$. A vertex $v$ is \emph{unstable} if $\mS(v)$ is  greater than or equal to its degree $\deg(v)$. An unstable vertex $v$ can \emph{topple} by sending a grain of sand through each of the incident edges. This leads to the new sandpile configuration $\mS'$ defined by $\mS'(u)=\mS(u)+\deg(u,v)$ for all $u\neq v$ and $\mS'(v)=\mS(v)-\deg(v,*)$, where $\deg(u,v)$ is the number of edges with endpoints $u,v$ and $\deg(v,*)$ is the number of non-loop edges incident to $v$. We denote this transition by $\mS\topple{v} \mS'$. An evolution of the system is represented in Figure \ref{fig:sable-diamant}.\\

A sandpile configuration is \emph{stable} if every vertex $v\neq v_0$ is stable. A stable configuration $\mS$  is \emph{recurrent} if  $\mS(v_0)=\deg(v_0)$ and if there is a labeling of the vertices in $V$ by $v_0,v_1,\ldots,v_{|V|-1}$ such that  $\mS\topple{v_0}\mS_1 \topple{v_1}\ldots  \topple{v_{|V|-1}} \mS_{|V|}=\mS$. This means that after toppling the root-vertex $v_0$, there is a valid sequence of toppling involving each vertex once that gets back to the initial configuration. For instance, the configuration at the left of Figure \ref{fig:sable-diamant} is recurrent. Lastly, the \emph{level} of a recurrent configuration $\mS$ is given by: 
$\textrm{level}(\mS)=\left(\sum_{v\in V}\mS(v) \right)-|E|.$\\

\section{A GLIMPSE AT THE RESULTS}\label{section:glimpse}
\subsection{Tour of spanning trees and embedding-activities}
\begin{figure}[htb!]
\begin{center}
\input{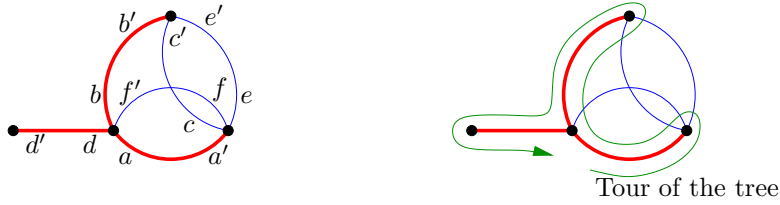}
\caption{Intuitive representation of the tour of a spanning tree (indicated by thick lines).}\label{fig:tour-of-tree}
\end{center}
\end{figure}

We first define the tour of spanning  trees. Informally, the tour of a tree is a walk around the tree that follows internal edges and crosses external edges.  A graphical representation of the tour is given in Figure \ref{fig:tour-of-tree}.\\

Let $\mathcal{G}=(H,\sigma,\alpha)$ be an embedding of the graph $G=(V,E)$. Given a spanning tree $T$, we define the \emph{motion function} $t$ on the set $H$ of half-edges by:
\begin{eqnarray} \label{def:t-def-Tutte}
\begin{array}{l|ll}
t(h)=& \sigma(h) & \textrm{ if } h  \textrm{ is external,}\\
& \sigma  \alpha(h) & \textrm{ if } h \textrm{ is internal.}
\end{array}
\end{eqnarray}
It was proved in \cite{OB:Tutte-plongement-def} that $t$ is a cyclic permutation on $H$. For instance, for the embedded graph of Figure
\ref{fig:tour-of-tree}, the motion function is the cyclic permutation  $(a,e,f,c,a',f',$ $b,c',e',b',d,d')$. The cyclic order defined by the
motion function $t$ on the set of half-edges is what we call the \emph{tour of the tree} $T$. \\


We will now define the \emph{embedding-activities} of spanning trees introduced in \cite{OB:Tutte-plongement-def} in order to characterize the Tutte polynomial (see Theorem \ref{thm:activity-OB} below).


\begin{Def}\label{def:GT-order}
Let $\mathcal{G}=(H,\sigma,\alpha,h)$ be an embedded graph and let $T$ be a spanning tree. We define the \emph{$(\mG,T)$-order} on the set $H$ of half-edges by $h<t(h)<t^2(h)<\ldots<t^{|H|-1}(h)$, where $t$ is the motion function.  (Note that the $(\mG,T)$-order is a linear order on $H$ since $t$ is a cyclic permutation.) We define the  $(\mG,T)$-order on the edge set by setting $e=\{h_1,h_2\}<e'=\{h_1',h_2'\}$ if $\min(h_1,h_2)<\min(h_1',h_2')$. (Note that this is also a linear order.)
\end{Def}

\noindent \textbf{Example:} Consider the embedded graph $\mG$ rooted on $a$ and the spanning tree $T$ represented in Figure \ref{fig:tour-of-tree}.  The $(\mG,T)$-order on the half-edges is $a<e<f<c<a'<f'<b<c'<e'<b'<d<d'$. Therefore, the $(\mG,T)$-order on the edges is $\{a,a'\}<\{e,e'\}<\{f,f'\}<\{c,c'\}<\{b,b'\}<\{d,d'\}$.\\


\begin{Def}\label{def:GT-activity}
Let  $\mG$ be a rooted embedded graph and $T$ be a spanning tree.   We say that an external (resp. internal) edge is
\emph{$(\mG,T)$-active} (or \emph{embedding-active} if $\mG$ and $T$
are clear from the context) if it is minimal for the $(\mG,T)$-order in its fundamental cycle (resp. cocycle).  
\end{Def}

\noindent \textbf{Example:} In Figure \ref{fig:tour-of-tree}, the $(\mG,T)$-order on the edges is $\{a,a'\}<\{e,e'\}<\{f,f'\}<\{c,c'\}<\{b,b'\}<\{d,d'\}$. Hence, the internal active edges are $\{a,a'\}$ and $\{d,d'\}$ and there is no external active edge. For instance, $\{e,e'\}$ is not active since $\{a,a'\}$ is in its fundamental cycle. \\

The following characterization of the Tutte polynomial was proved in \cite{OB:Tutte-plongement-def}.

\begin{thm} \label{thm:activity-OB}
Let $\mG$ be any rooted embedding of the connected graph $G$ (with at least one edge).  The Tutte polynomial of $G$ is equal to 
\begin{eqnarray}
T_G(x,y)=\sum_{T \textrm{ spanning tree}} x^{\mI(T)}y^{\mE(T)}, \label{eq:Tutte-embedded}
\end{eqnarray}
where the sum is over all spanning trees and $\mI(T)$ (resp. $\mE(T)$) is the number of $(\mG,T)$-active internal (resp. external) edges.
\end{thm}

\begin{figure}[htb!]
\begin{center}
\input{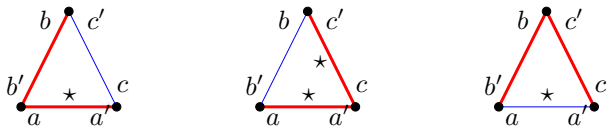}
\caption{The embedding-activities of the spanning trees of $K_3$.}\label{fig:exp-Tutte-triangle-OB}
\end{center}
\end{figure}

\titre{Example: } We represented the spanning trees of $K_3$ in Figure \ref{fig:exp-Tutte-triangle-OB}. If the embedding is rooted on the half-edge $a$, then the embedding-active edges are the one indicated by a $\star$. Hence, the spanning trees (taken from left to right) have respective contributions $x$, $x^2$ and $y$ and the Tutte polynomial is $T_{K_3}(x,y)=x^2+x+y$.\\

Note that the characterization \Ref{eq:Tutte-embedded} of the Tutte polynomial \ref{thm:activity-OB} implies that the sum in \Ref{eq:Tutte-embedded} does not depend on the embedding, whereas the summands clearly depends on it. This characterization is reminiscent but inequivalent to the one given by Tutte in \cite{Tutte:dichromate}. \\

From now on we adopt the following conventions. If an embedding $\mG$ and a spanning tree $T$ are clear from the context, the $(\mG,T)$ order is denoted by $<$. If $F$ is a set of edges and $h$ is a half-edge, we say that $h$ is \emph{in} $F$ if the edge $e$ containing $h$ is in $F$. A half-edge $h$ is said to be \emph{internal}, \emph{external} or \emph{$(\mG,T)$-active}  if the edge $e$ is. \\

We now make some elementary remarks about embedding-activities that will be useful throughout the paper.\\

\begin{lemma} \label{lemma:between-h1-and-h2}
Let $\mG$ be an embedded graph. Let $T$ be a spanning tree and let $e=\{h_1,h_2\}$ be an internal edge. Assume that $h_1<h_2$ (for the $(\mG,T)$-order) and denote by $v_1$ and $v_2$ the endpoints of $h_1$ and $h_2$ respectively. Then, $v_1$ is the father of $v_2$ in $T$. Moreover, the half-edges $h$ such that $h_1<h\leq h_2$ are the half-edges incident to a descendant of $v_2$. 
\end{lemma}

\dem
Let $t$ be the motion function associated to the tree $T$ ($t$ is defined by \Ref{def:t-def-Tutte}). We consider the subtrees $T_1$ and $T_2$ obtained from $T$ by deleting $e$ with the convention that $h_1$ is incident to $T_1$ and $h_2$ is incident to $T_2$. 
Let $h$ be any half-edge distinct from $h_1$ and $h_2$. By definition of $t$, the half-edges $h$ and $t(h)$ are incident to the same subtree $T_i$. Therefore, the $(\mG,T)$-order is such that $h_0<l_1<\cdots <l_i<h_1<l_1'<\cdots <l_j' <h_2<l_1''< \cdots< l_k''$ where $l_1',\ldots ,l_j',h_2$ are the half-edges incident with the subtree $T_2$ not containing the root-vertex $v_0$. Since the subtree $T_2$ does not contain $v_0$ its vertices are the descendants of $v_2$ in~$T$.
\findem

\begin{lemma}\label{lemma:cycle-h1-h1-h2-h2}
With the same assumption as in Lemma \ref{lemma:between-h1-and-h2}, let $e=\{h_1,h_2\}$ with $h_1<h_2$ be an internal edge and let $e'=\{h_1',h_2'\}$ with $h_1'<h_2'$  be an external edge.\\ 
\ite Then, $e$ is in the fundamental cycle of $e'$ (equivalently, $e'$ is in the fundamental cocycle of $e$) if and only if $h_1<h_1'<h_2<h_2'$ or  $h_1'<h_1<h_2'<h_2$. \\
\ite  Suppose that  $e$ is in the fundamental cycle of $e'$ and denote by  $v_1,v_2,v_1',v_2'$ the endpoints of $h_1,h_2,h_1',h_2'$ respectively. Recall that $v_1$ is the father of $v_2$ in $T$ (Lemma \ref{lemma:between-h1-and-h2}) and that exactly one of the vertices $v_1',~v_2'$ is a descendant of $v_2$. If $e<e'$, then $v_1'$ is the descendant of $v_2$,  else it is $v_2'$.
\end{lemma}

\dem\\
\ite Let $V_2$ be the set of descendants of $v_2$. Recall that the edge $e'$ is in the fundamental cocycle of $e$ if and only if it has one endpoint in $V_2$ and the other in $\B{V_2}$. By Lemma~\ref{lemma:between-h1-and-h2}, this is equivalent to the fact that exactly one of the half-edges $h_1',h_2'$ is in $\{h': h_1<h'\leq h_2\}$. Thus, $e'$ is in the fundamental cocycle of $e$ if and only if $h_1<h_1'<h_2<h_2'$ or  $h_1'<h_1<h_2'<h_2$.\\
\ite Suppose that  $e$ is in the fundamental cycle of $e'$. By the preceding point, $e<e'$ implies $h_1<h_1'<h_2<h_2'$. In this case, $h_1'$ is incident to a descendant of $v_2$ by Lemma \ref{lemma:between-h1-and-h2}. Similarly, $e'<e$ implies $h_1'<h_1<h_2'<h_2$, hence $h_2'$ is incident to a descendant of $v_2$.
\findem

\begin{lemma}\label{lem:external-active}
An external edge $e'=\{h_1',h_2'\}$ with $h_1'<h_2'$ is $(\mG,T)$-active if and only if the endpoint of $h_1'$ is an ancestor of the endpoint of $h_2'$.
\end{lemma}

\dem Denote by $v_1'$ and $v_2'$  the endpoints of $h_1'$ and $h_2'$ respectively.  \\
\ite Suppose $v_1'$ is an ancestor of $v_2'$. We want to prove that  $e'$ is active. Let $e=\{h_1,h_2\}$ with $h_1<h_2$ be an internal edge in the fundamental cycle of $e'$. The edge $e$ is in the path of $T$ between $v_1'$ and $v_2'$.  Denote by $v_1$ and $v_2$  the endpoints of $h_1$ and $h_2$ respectively. Recall that $v_1$ is the father of $v_2$ (Lemma \ref{lemma:between-h1-and-h2}). Since $v_2'$ is a descendant of $v_2$, we have $e'<e$  by  Lemma  \ref{lemma:cycle-h1-h1-h2-h2}.  The edge $e'$ is less than any edge in its fundamental cycle hence it is $(\mG,T)$-active.\\
\ite Suppose that $v_1'$ is not an ancestor of $v_2'$. Then the edge $e=\{h_1,h_2\}$ with $h_1<h_2$ linking $v_1'$ to its father in $T$ is in the fundamental cycle of $e'$. If we denote by  $v_1$ and $v_2$ the endpoints of $h_1$ and $h_2$ respectively, we get $v_2=v_1'$ by  Lemma \ref{lemma:between-h1-and-h2}. Since the endpoint  $v_1'$ of $h_1'$ is  a descendant of the endpoint $v_2$ of $h_2$, we get $e<e'$ by Lemma \ref{lemma:cycle-h1-h1-h2-h2}. Thus, $e'$ is not $(\mG,T)$-active. 
\findem

\subsection{A mapping from spanning trees to orientations and some related bijections} 
We now take a glimpse at the results to be developed in the following sections. In order to present these results, we define a mapping $\Phi$ from spanning trees to orientations. The mapping $\Phi$ will  be extended into a bijection between subgraphs and orientations  in Section \ref{section:subgraphs-orientations}. Related to the mapping $\Phi$, we define two other mappings $\Gamma$ and $\Lambda$ on the set of spanning trees. The mapping $\Gamma$ is a bijection between spanning trees and root-connected outdegree sequences  while $\Lambda$ is a bijection between spanning trees and recurrent sandpile configurations.\\

\begin{figure}[ht!]
\begin{center}
\input{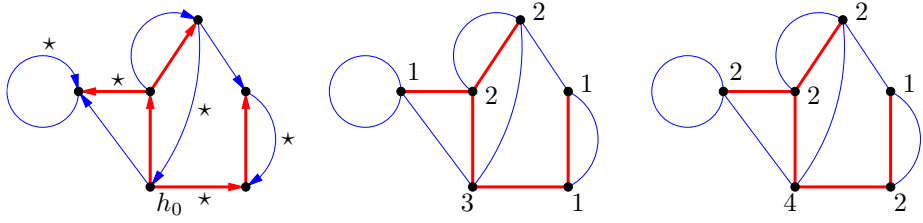}
\caption{Left: Orientation $\mO_T$ associated the spanning tree $T$ (indicated by thick lines) and active edges (indicated by a star). Middle: outdegree sequence $\Gamma(T)$. Right: recurrent sandpile configuration $\Lambda(T)$.} \label{fig:canonical-orientation}
\end{center}
\end{figure}

Consider an embedded graph $\mG=(H,\sigma,\alpha,h_0)$ and a spanning tree $T$. Recall that the tour of $T$ defines a linear order, the  \emph{$(\mG,T)$-order}, on $H$ for which the root $h_0$ is the least element. We associates with the spanning tree $T$ the orientation $\mO_T$ of $\mG$ defined by:
\begin{eqnarray} \label{def:O_T}
\hspace{-.2cm} \begin{array}{l|ll} \textrm{For any edge } e=\{h_1,h_2\} \textrm{ with } h_1<h_2, ~
\mO_T(e)=& (h_1,h_2) & \textrm{ if } e  \textrm{ is internal,}\\
& (h_2,h_1) & \textrm{ if } e \textrm{ is external.}
\end{array}
\end{eqnarray}
This definition is illustrated in Figure \ref{fig:canonical-orientation} (left). 

Observe that the spanning tree $T$ is oriented from its root-vertex $v_0$ to its leaves in $\mO_T$. Indeed, it is clear from the definitions and Lemma \ref{lemma:between-h1-and-h2} that every internal edge is oriented from father to son. This property implies that for every spanning tree $T$ the orientation $\mO_T$ is $v_0$-connected.\\

The  mapping $\Phi:T\mapsto \mO_T$ from spanning trees to $v_0$-connected orientations is not bijective. However, it is injective and  in Section \ref{section:subgraphs-orientations}  we will extend it into a bijection between subgraphs and orientations. For the time being, let us observe (the proof will be given in Section \ref{section:subgraphs-orientations}) that the tree $T$ can be recovered from the orientation $\mO_T$ by the following procedure:\\

\noindent \textbf{Procedure Construct-tree:}\\
\noindent \textbf{Initialization:} Initialize the \emph{current half-edge} $h$ to be the root $h_0$. Initialize the tree $T$ and the set of visited arcs $F$ to be empty.\vspace{0cm}\\
\textbf{Core:} Do: \vspace{.0cm}\\
\textbf{C1:} \hspace{-1pt} If the edge $e$ containing $h$ is not in $F$ and  $h$ is a tail then add $e$ to $T$.\\
\indent \indent Add $e$ to $F$.\vspace{.0cm}\\
\textbf{C2:}  \hspace{-1pt} Move to the next half-edge around $T$: \\
\indent \indent If $e$ is in $T$, then set the current half-edge $h$ to be $\sigma\alpha(h)$, else set it to be $\sigma(h)$.\vspace{.0cm}\\
Repeat until the current half-edge $h$ is $h_0$.\vspace{.0cm}\\
\textbf{End:} Return the tree $T$.\\

In the procedure \textbf{Construct-tree} we keep track of the set $F$ of edges already visited. The decision of adding an edge $e$ to the tree $T$ or not is taken when $e$ is visited for the \emph{first} time.  The principle of procedure \textbf{Construct-tree}, which consists in constructing a tree $T$ while making its tour, will appear again in the next sections.\\

Building on the mapping $\Phi: T\mapsto \mO_T$, we define two mappings $\Gamma$ and $\Lambda$.
\begin{Def} \label{def:Gamma}
Let $\mG$ be an embedded graph. The mapping $\Gamma$ associates with any spanning tree $T$ the outdegree sequence of the orientation $\mO_T$.
\end{Def}

\begin{Def} \label{def:Lambda}
Let $\mG$ be an embedded graph and let $V$ be the vertex set. The mapping $\Lambda$ associates with any spanning tree $T$ the sandpile configuration $\mS_T:V\mapsto \mathbb{N}$, where $\mS_T(v)$ is the number of tails plus the number of external $(\mG,T)$-active heads incident to $v$ in the orientation $\mO_T$.
\end{Def}

The mappings $\Gamma$ and  $\Lambda$  are illustrated in  Figure \ref{fig:canonical-orientation}.  \\

As observed above, the orientation $\mO_T$ is always $v_0$-connected hence the image of any spanning tree by the mapping $\Gamma$ is a $v_0$-connected outdegree sequence. We shall prove in Section \ref{section:specializations}  that $\Gamma$ is a bijection between spanning tree and $v_0$-connected outdegree sequences. We will also show how to extend it into a bijection between forests and outdegree sequences. 
Regarding the mapping  $\Lambda$, we shall prove in  Section~\ref{section:bij-sandpile} that it is a bijection between spanning trees and recurrent sandpile configurations. Moreover, the number of external $(\mG,T)$-active edges is easily seen to be the level of the configuration $\Lambda(T)$. This gives a new bijective proof of a result by Merino linking external activities to the level of recurrent sandpile configurations \cite{Borgne-cori-activite-externe-sable,Merino:external-activity=sandpile-level}.\\

\begin{figure}[ht!]
\begin{center}
\input{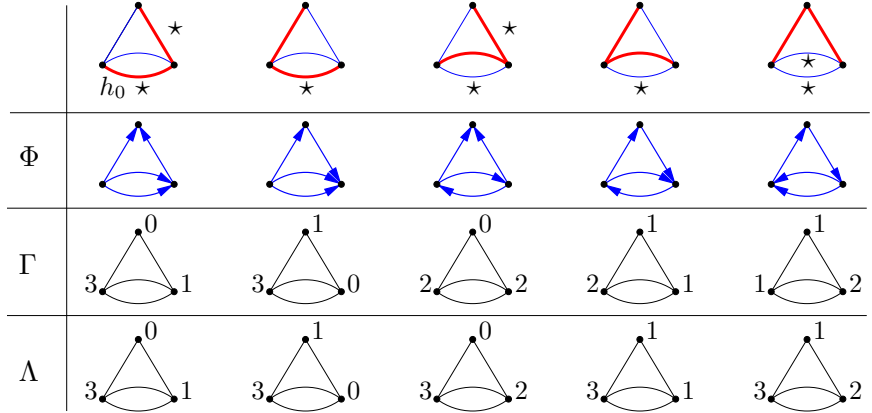}
\caption{Spanning trees (embedding-active edges are indicated by a star) and their image by the mappings $\Phi$, $\Gamma$ and $\Lambda$.} \label{fig:exp-Gamma-lambda}
\end{center}
\end{figure}

The two mappings $\Gamma$ and  $\Lambda$  coincide on \emph{internal} trees, that is, trees that have external activity 0. Thus, the mapping $\Gamma\circ \Lambda^{-1}$ is a bijection between recurrent sandpile configurations and $v_0$-connected outdegree sequences that leaves the configurations at level 0 unchanged. This answers a problem raised by Gioan~\cite{Gioan:enumerating-degree-sequences}. As an illustration we represented the 5 spanning trees of a graph in Figure \ref{fig:exp-Gamma-lambda} and their image by the mappings $\Phi$, $\Gamma$ and $\Lambda$ (the first two spanning trees are internal).\\

We now highlight a relation (to be exploited in  Section \ref{section:specializations}) between the embedding-activities of the spanning tree $T$ and the acyclicity or strong connectivity of the associated orientation $\mO_T$.
\begin{lemma}\label{lem:fundamental-is-directed}
Let $\mG$ be an embedded graph ant let $T$ be a spanning tree. The fundamental cycle (resp. cocycle) of an external (resp. internal) edge $e$ is $\mO_T$-directed if and only if $e$ is $(\mG,T)$-active.
\end{lemma}

\begin{figure}[ht!]
\begin{center}
\input{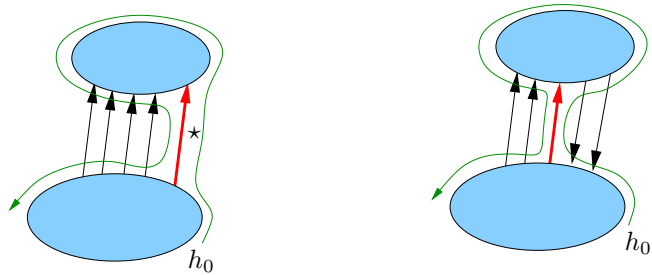}
\caption{Fundamental cocycles of an active internal edge (left) and of a non-active internal edge (right).} \label{fig:cocycle-oriented}
\end{center}
\end{figure}

\begin{figure}[ht!]
\begin{center}
\input{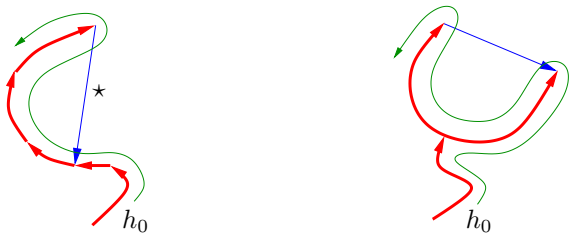}
\caption{Fundamental cycles of an active external edge (left) and of a non-active external edge (right).} \label{fig:cycle-oriented}
\end{center}
\end{figure}

Lemma \ref{lem:fundamental-is-directed} is illustrated by  Figures \ref{fig:cocycle-oriented} and \ref{fig:cycle-oriented}. From this lemma   we deduce that if $\mO_T$ is acyclic (resp. strongly connected) then $T$ is \emph{internal} (resp. \emph{external}), that is, has no external (resp. internal) active edge. In fact, we shall prove in Section \ref{section:specializations} that the converse is true: if the tree $T$ is internal (resp. external), then the orientation $\mO_T$ is acyclic (resp. strongly connected). For instance, in Figure \ref{fig:exp-Gamma-lambda} the two first (last) spanning trees are internal (resp. external) and the corresponding orientations are acyclic (resp. strongly connected).\\ 

\dem Consider an edge  $e=\{h_1,h_2\}$ with $h_1<h_2$ and denote by $v_1$ and $v_2$ the endpoints of $h_1$ and $h_2$ respectively.  \\
\ite Suppose that $e$ is internal. We want to prove that the fundamental cocycle $D$ of $e$ is directed if and only if  $e$ is  $(\mG,T)$-active. Recall that $v_1$ is the father of $v_2$ by Lemma \ref{lemma:between-h1-and-h2}.  Let $V_2$ be the set of descendants of~$v_2$. Recall that $D$ is the cocycle defined by $V_2$. By definition, the arc $\mO_T(e)$ is directed toward $v_2\in V_2$.   By Lemma \ref{lemma:cycle-h1-h1-h2-h2}, for all edge $e'=\{h_1',h_2'\}$ with $h_1'<h_2'$  in $D-e$, the arc   $\mO_T(e')=(h_2',h_1')$ is directed toward $V_2$ if and only if $e<e'$. Therefore, the fundamental cocycle $D$ is directed if and only if $e$ is minimal in $D$, that is, if $e$ is  $(\mG,T)$-active. \\ 
\ite  Suppose that $e$ is external.  We want to prove that the fundamental cycle $C$ of $e$ is directed if and only if  $e$ is  $(\mG,T)$-active. Recall that $C-e$ is the path in $T$ between $v_1$ and $v_2$. Since $\mO_T(e)$ is directed toward $v_1$, the cycle  $C$ is directed if and only if the path $C-e$ is directed from $v_1$ to $v_2$. Since every edge in $C-e\subseteq T$ is directed from father to son (Lemma \ref{lemma:between-h1-and-h2}), the cycle $C$ is directed if and only if $v_1$ is an ancestor of $v_2$. This is precisely the characterization of external $(\mG,T)$-active edges given by Lemma \ref{lem:external-active}.
\findem

Up to this point we have considered mappings defined on the set of spanning trees. In order to extend these mappings to general subgraphs we will associate a spanning tree to every subgraph. This is the task of the next section. \\

\section{A PARTITION OF THE SET OF SUBGRAPHS}\label{section:partition} 
In this section we define a partition of the set of subgraphs for any embedded graph. Each part of this partition is associated with a spanning tree.\\

\begin{figure}[ht!]
\begin{center}
\input{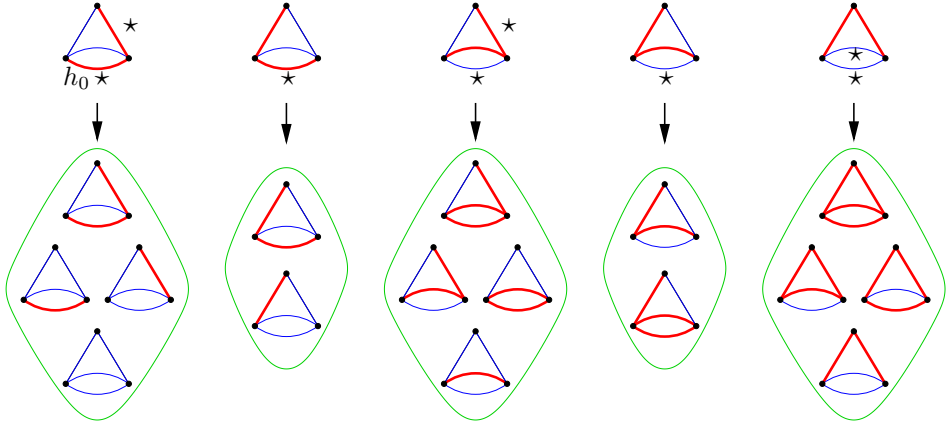}
\caption{The tree-intervals corresponding to each spanning tree. The active edges are indicated by a $\star$.} \label{fig:interval-exp}
\end{center}
\end{figure}

Let $\mG$ be an embedded graph. Given a spanning tree $T$, we consider the set of subgraphs that can be obtained from $T$ by removing some internal $(\mG,T)$-active edges and adding some external $(\mG,T)$-active edges. Observe that this set is an interval in the boolean lattice of the subgraphs of $\mG$ (i.e. subsets of edges). We call \emph{tree-interval} and denote by $[T^-,T^+]$ the set of subgraphs obtained from a spanning tree $T$. We represented the tree-intervals corresponding to each of the 5 spanning trees of the embedded graph in Figure \ref{fig:interval-exp}. \\

We prove some properties of the subgraphs in the tree-interval $[T^-,T^+]$.

\begin{lemma}\label{thm:fdmtl-cycle-is-in}
Let $\mG$ be an embedded graph and let $T$ be a spanning tree. Let $e$ be an internal (resp. external) $(\mG,T)$-active edge.  The fundamental cocycle (resp. cycle) of $e$ is contained in $\B{S}+e$ (resp. $S+e$) for any subgraph $S$ in $[T^-,T^+]$.  
\end{lemma}

\dem
If $e$ is internal and $(\mG,T)$-active, no edge in its fundamental cocycle $D$ is $(\mG,T)$-active (since their fundamental cycle contains $e$). Since no edge of $D-e$ is in $T$ nor is $(\mG,T)$-active, none is in $S$. Hence, $D\subseteq \B{S}+e$. Similarly, if $e$ is external $(\mG,T)$-active, its fundamental cycle is contained in $S+e$. 
\findem

\begin{lemma} \label{lem:keylink}
Let $\mG$ be an embedded graph. Let $T$ be a spanning tree and let $S$ be a subgraph in $[T^-,T^+]$ having $c(S)$ connected components. Then $c(S)-1$, (resp. $e(S)+c(S)-|V|$) is the number of edges in $\B{S}\cap T$ (resp. $S\cap\B{T}$).
\end{lemma}

\dem
Consider any subgraph $S$ in $[T^-,T^+]$. By Lemma \ref{thm:fdmtl-cycle-is-in}, removing an internal $(\mG,T)$-active edge from $S$ increases  $c(S)$ by one and leaves $e(S)+c(S)$ unchanged. Similarly, adding an external $(\mG,T)$-active edge to $S$ leaves $c(S)$ unchanged and increases $e(S)+c(S)$ by one. Moreover, $c(T)-1=0$ and $e(T)+c(T)-|V|=0$. Therefore, Lemma  \ref{lem:keylink} holds for every subgraph $S$ in $[T^-,T^+]$ by induction on the number of edges in $S\vartriangle T$.
\findem

\begin{figure}[ht!]
\begin{center}
\input{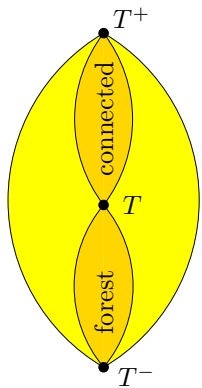}
\caption{The tree-interval $[T^-,T^+]$, the sub-interval $[T,T^+]$ of connected subgraphs and the sub-interval $[T ^-,T]$ of forests.} \label{fig:generic-tree-interval}
\end{center}
\end{figure}

By Lemma \ref{lem:keylink}, the connected subgraphs in $[T^-,T^+]$ are the subgraphs in the interval $[T,T^+]$ (the subgraphs obtained from $T$ by adding some external $(\mG,T)$-active edges). Similarly, the forests in $[T^-,T^+]$ are the subgraphs in the interval $[T^-,T]$ (the subgraphs obtained from $T$ by removing some internal $(\mG,T)$-active edges). These properties are illustrated in Figure \ref{fig:generic-tree-interval}. \\

We are now ready to state and comment on the main result of this section.\\

\begin{thm}\label{thm:partition}
Let  $G=(V,E)$ be a graph and let $\mG$ be an embedding of $G$. The tree-intervals form a partition of the set of subgraphs of $G$:
$$2^E=\biguplus_{T \textrm{ spanning tree}} [T^-,T^+],$$
where the disjoint union is over all spanning trees of $G$.
\end{thm}

The counterpart of Theorem \ref{thm:partition} is known for the notion of (internal and external)\emph{activities} defined by Tutte in \cite{Tutte:dichromate}. 
This property has been used to extract informations about the Tutte polynomial in \cite{Bari:Tutte-poly-partition,Crapo-Tutte-poly,gordon:generalized-activities-Tutte-poly}. \\

Theorem  \ref{thm:partition} constitutes the key link between the subgraph expansion \Ref{eq:intro-Tutte-subgraph} and spanning tree expansion \Ref{eq:intro-Tutte-OB} of the Tutte polynomial. Indeed, given Lemma \ref{lem:keylink}, we get  
$$\sum_{S \in [T^-,T^+]}(x-1)^{c(S)-1}(y-1)^{e(S)+c(S)-|V|}=(x-1+1)^{\mI(T)}(y-1+1)^{\mE(T)}=x^{\mI(T)}y^{\mE(T)},$$
where $\mI(T)$ (resp. $\mE(T)$) is the number of internal (resp. external) $(\mG,T)$-active edges. Summing over all spanning trees gives the identity:
$$\sum_{S \textrm{ subgraph}}(x-1)^{c(S)-1}(y-1)^{e(S)+c(S)-|V|}=\sum_{T \textrm{ spanning tree}}x^{\mI(T)}y^{\mE(T)}.$$

\noindent \textbf{Remark.} As observed in \cite{gordon:generalized-activities-Tutte-poly}, the partition of the set of subgraphs gives several other expansions of the Tutte polynomial. For instance, the tree-intervals can be partitioned into \emph{forest-intervals}. The \emph{forest-interval} of a forest $F$ in $[T^-,T^+]$ is the set $[F,F^+]$ of subgraphs obtained from $F$ by adding some external $(\mG,T)$-active edges. Since 
$$[T^-,T^+]=\biguplus_{F \textrm{ forest in } [T^-,T^+]} [F,F^+],$$ 
the partition into tree-intervals given by Theorem \ref{thm:partition} leads to a partition into  forest-intervals:
$$2^E=\biguplus_{F \textrm{ forest}} [F,F^+].$$
Given Lemma \ref{lem:keylink}, we get  
$$\sum_{S \in [F,F^+]}(x-1)^{c(S)-1}(y-1)^{e(S)+c(S)-|V|}=(x-1)^{c(F)-1}(y-1+1)^{\mE(T)}=(x-1)^{c(F)-1}y^{\mE(T)},$$
for any forest in $[T^-,T^+]$. Summing up over forests,  gives the \emph{forest expansion}   
$$T_G(x,y)=\sum_{F \textrm{ forest}}(x-1)^{c(F)-1}y^{\mE(F)},$$
where $\mE(F)$ is the number of $(\mG,T)$-active edges for the spanning tree $T$ such that $F\in[T^-,T^+]$. Let us mention that several alternative notions of \emph{external activities} have been defined, each of which gives a forest expansion \cite{Gessel:Tutte-poly+DFS,Kostic:multiparking-function} which can be used to obtain enumerative results about the Tutte polynomial \cite{Gessel:enumerative-csq-DFS,Gessel:Tutte-poly+DFS}.  \\


In order to  prove Theorem \ref{thm:partition} we define a mapping $\Delta$ from subgraphs to spanning trees.

\begin{Def}\label{def:Delta}
Let $\mG$ be an embedded graph rooted on $h_0$ and  let $S$ be a subgraph. The spanning tree $T=\Delta(S)$ is defined by the following procedure:

\noindent \textbf{Initialization:} Initialize the \emph{current half-edge}  $h$ to be the root $h_0$. Initialize the tree $T$ and the set of visited edges $F$ to be empty.\vspace{0cm}\\
\textbf{Core:} Do: \vspace{.0cm}\\
\textbf{C1:} If the edge $e$ containing $h$ is not in $F$, then decide whether to add $e$ to $T$  according to the following rule:\\
\indent \indent \hspace{.8cm} If ($e$ is in $S$ and is in no cycle $C\subseteq S\cap \B{F}$) or \\
\indent \indent \hspace{2.4cm} ($e$ is not in $S$ and is in a cocycle $D\subseteq \B{S}\cap\B{F}$), \\
\indent \indent \hspace{1.6cm} Then add $e$ to $T$.\\
\indent \indent \hspace{.8cm} Endif.\\
\indent \indent Endif.\\
\indent \indent Add $e$ to $F$.\\
\textbf{C2:}  Move to the next half-edge around $T$: \\
\indent \indent If $e$ is in $T$, then set the current half-edge $h$ to be $\sigma\alpha(h)$, else set it to be $\sigma(h)$.\vspace{.0cm}\\
Repeat until the current half-edge $h$ is $h_0$.\vspace{.0cm}\\
\textbf{End:} Return the tree $T$.
\end{Def}

\begin{figure}[ht!]
\begin{center}
\input{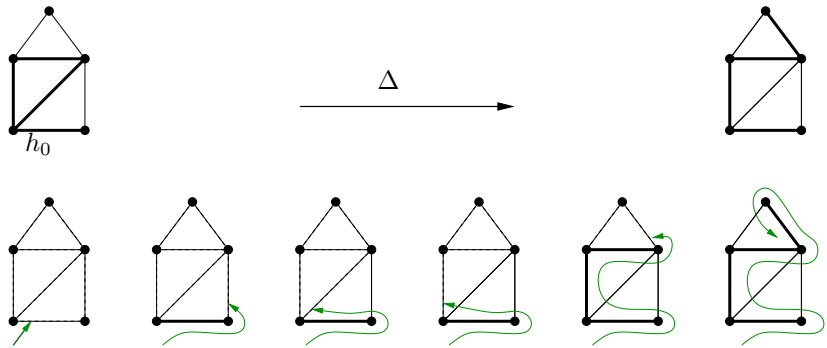}
\caption{The mapping $\Delta$ and some intermediate steps. The dashed lines correspond to the set $\B{F}$ of unvisited edges.} \label{fig:exp-Delta}
\end{center}
\end{figure}

An execution of the procedure $\Delta$ is illustrated in Figure \ref{fig:exp-Delta}.\\

There is a \emph{direct} proof that the mapping $\Delta$ is well defined on every subgraph (that is, the procedure terminates and returns a spanning tree). But we shall only prove an (a priory weaker) result: the mapping $\Delta$ is well defined on every tree-interval and  $\Delta(S)=T$ for any subgraph $S$ in $[T^-,T^+]$ (Proposition \ref{thm:Delta-on-intervals}). This will prove that the tree-intervals are disjoint. Moreover, the cardinality of the tree-interval $[T^-,T^+]$ is $2^{\mI(T)+\mE(T)}$, where $\mI(T)$ and $\mE(T)$ are the number of internal and external $(\mG,T)$-active edges. Therefore, the number of subgraphs contained in some tree-intervals is 
$$\left|\bigcup_{T \textrm{ spanning tree}} [T^-,T^+]\right|=\sum_{T \textrm{ spanning tree}} \left|[T^-,T^+]\right|=\sum_{T \textrm{ spanning tree}} 2^{\mI(T)+\mE(T)}.$$
By Theorem \ref{thm:activity-OB}, this sum is the specialization $T_G(2,2)$ of the Tutte polynomial counting the subgraphs of $G$ (as is clear from \Ref{eq:intro-Tutte-subgraph}). This counting argument proves that every subgraph belongs to a tree-interval. Thus, we only need to prove the following proposition.

\begin{prop}\label{thm:Delta-on-intervals}
Let $\mG$ be an embedded graph. Let $T$ be a spanning tree and let $S$ be a subgraph in the tree-interval $[T^-,T^+]$. The procedure $\Delta$ is well defined on $S$ and returns the tree $T$. 
\end{prop}
 
Before proving Proposition \ref{thm:Delta-on-intervals}, we need to recall a classical result of graph theory.

\begin{lemma}\label{lem:two-cycles-unoriented} 
The symmetric difference of  two cycles (resp. cocycles) $C$ and $C'$ is a union of cycles (resp. cocycles).
\end{lemma}

\begin{figure}[ht!]
\begin{center}
\input{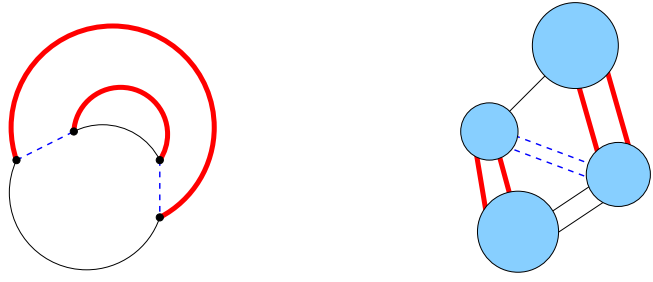}
\caption{Left: Two cycles (thin and thick lines) and their intersection (dashed lines). Right: two cocycles.} \label{fig:two-cycles-unoriented}
\end{center}
\end{figure}

Lemma \ref{lem:two-cycles-unoriented}  is illustrated by Figure \ref{fig:two-cycles-unoriented}.\\

We now characterize the edges in the symmetric difference $S\vartriangle T$.

\begin{lemma} \label{thm:in-sym-difference}
Let $\mG$ be an embedded graph. Let $T$ be a spanning tree and let $S$ be a subgraph in the tree-interval $[T^-,T^+]$. \\
(i) An edge $e$ is in $S\cap\B{T}$ if and only if $e$ is minimal (for the $(\mG,T)$-order) in a cycle $C\subseteq S$.\\
(ii) An edge $e$ is in $\B{S}\cap T$ if and only if $e$ is minimal (for the $(\mG,T)$-order) in a cocycle $D\subseteq \B{S}$.
\end{lemma}

\begin{figure}[ht!]
\begin{center}
\input{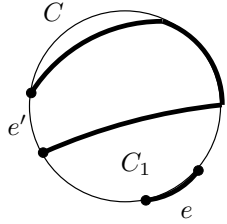}
\caption{The cycle $C$ (circle), some edges in the tree $T$ (indicated by thick lines) and the edges $e$ and $e'$.} \label{fig:symdifference}
\end{center}
\end{figure}

\dem 
We give the proof of $(i)$; the proof of $(ii)$ is similar.\\
\ite Suppose $e$ is in $S\cap\B{T}$. Then $e$ is $(\mG,T)$-active, that is, $e$ is minimal in its fundamental cycle $C$. Moreover, by Lemma \ref{thm:fdmtl-cycle-is-in}, $C$ is contained in $S$.\\
\ite Suppose $e$ is minimal in a cycle $C\subseteq S$. We want to prove that $e$ is in $\B{T}$. Suppose the contrary. Then, there is an edge $e'\neq e$ in $C\cap\B{T}$ (since $T$ has no cycle). Take the least edge $e'$ in  $C\cap\B{T}$  and consider its fundamental cycle $C'$. The edge $e'$ is $(\mG,T)$-active, that is, $e'$ is minimal in $C'$. In particular, $e$ is not in $C'$. This situation is represented in Figure \ref{fig:symdifference}. Since  $e$ is in $C\vartriangle C'$ and $e'$ is not, there is a cycle $C_1\subseteq C\vartriangle C'$ containing $e$ and not $e'$ (Lemma \ref{lem:two-cycles-unoriented}). By Lemma  \ref{thm:fdmtl-cycle-is-in},  the fundamental cycle $C'$ of $e'$ is contained is $S+e'$, thus $C_1\subseteq C\vartriangle C' \subseteq S$.  Note that $e$ is minimal in the cycle $C_1\subseteq S$ (since $e$ is minimal in $C$ and $e'>e$ is minimal in $C'$). Moreover, the least edge in $C_1\cap\B{T}$ (this edge exists since $T$ has no cycle) is in $C\cap \B{T}-e'$ (since $C'\subseteq T+e'$), hence is greater than $e'$. We can repeat this operation again in order to produce an infinite sequence $C_0=C,C_1,C_2,\ldots$ of cycles with $e$ minimal in $C_i$ and $C_i \subseteq S$ for all $i\geq 0$. But the minimal element of $C_i\cap\B{T}$ is strictly increasing with $i$. This is impossible.\vspace{.2cm}
\findem

\noindent \textbf{Proof of Proposition \ref{thm:Delta-on-intervals}.} We consider a subgraph $S$ in the tree-interval $[T_0^-,T_0^+]$. We denote by $H$ the set of half-edges. We denote by $t$ the motion function associated with spanning tree $T_0$ and we denote by $h_i=t^i(h_0)$ the $i^{th}$ half-edge for the $(\mG,T_0)$-order. For any half-edge $h$, we denote $F_h=\{e=\{h_1,h_2\}/ {\min}(h_1,h_2)<h \}$ and $T_h=T_0\cap F_h$.\\
We adopt the notations $h,~e,~F$ and $T$ of the procedure $\Delta$ (for instance, $h$ denotes the current half-edge) and we compare half-edges according to the $(\mG,T_0)$-order. We want to prove that, for all $i\leq |H|$,  \emph{at the beginning of the $i^{th}$ core step, $h=h_i$, $F=F_h$ and $T=T_h$}. We proceed by induction on $i$. The property holds for the first core step ($i=0$) since $h=h_0$ and $F_{h_0}=T_{h_0}=\emptyset$. Consider now the $i^{th}$ core step. Suppose first that the edge $e$ containing the current half-edge $h$ is not in $F$. By the induction hypothesis, $F=F_h$ thus $e$ is greater than any edge in $F$ and less than any edge in $\B{F}-e$.  By Lemma \ref{thm:in-sym-difference}, if $e$ is in $S$, then it is in $\B{T_0}$ if and only if it is in a cycle $C\subseteq S\cap \B{F}$. Also, if $e$ is in $\B{S}$, then it is in $T_0$ if and only if it is in a cocycle $D\subseteq\B{S}\cap \B{F}$. Therefore, the edge $e$ is added to $T$ at the step \textbf{C1} if and only if it is in $T_0$. Suppose now that the edge $e$ is already in $F$ at the beginning of the $i^{th}$ core step. Then, by the induction hypothesis, $e$ is in $T=T_h=T_0\cap F_h=T_0\cap F$ if and only if it is in $T_0$. Whether the edge $e$ is in $F$ or not at the beginning of the step \textbf{C1}, the edge $e$ is in $T$ at the beginning of the step \textbf{C2} if and only if it is in $T_0$. Therefore, the current half-edge at the  beginning of the $(i+1)^{th}$ core step, is $t(h)=h_{i+1}$.  Thus, the property holds for all $i\leq |H|$ by induction.  In particular, the procedure $\Delta$ stops after $|H|$ core steps and returns the spanning tree $T=T_{h_{|H|-1}}=T_0$. 
\findem

This concludes the proof of Theorem \ref{thm:partition}. \findem

Before we close this section we define some families of subgraphs counted by the evaluations $T_G(i,j),0\leq i,j\leq 2$ of the Tutte polynomial. Consider an embedded graph $\mG$ and a spanning tree $T$. Recall that the spanning tree $T$ is said to be \emph{internal} (resp. \emph{external}) if it has no external (resp. internal) $(\mG,T)$-active edge. For instance, among the spanning trees represented in Figure \ref{fig:interval-exp}, the two first (resp. last) are internal (resp. external). We say that a subgraph $S$ in $[T^-,T^+]$ is \emph{internal} or \emph{external} if the spanning tree $T$ is. The notion of \emph{internal subgraph} is close to Whitney's notion of subgraphs \emph{without broken circuit} \cite{Whitney:broken-circuit}. Observe that by Lemma \ref{lem:keylink} any internal subgraph is a forest and any external subgraph is connected (the converse is, of course, false).  In Figure \ref{fig:subgraphs-organized} we represented the subgraphs of figure \ref{fig:interval-exp} in each of the categories defined by the four criteria \emph{forest}, \emph{internal}, \emph{connected}, \emph{external}. \\

\begin{prop}\label{thm:counting-specilizations}
Let $\mG$ be an embedded graph. The number of subgraphs in each category defined by the criteria \emph{forest}, \emph{internal}, \emph{connected}, \emph{external} is given by the following evaluation of the Tutte polynomial: \\
\begin{center}
\begin{tabular}{|c|c|c|c|}
\hline
& General &Connected & External\\ \hline
General & $T_G(2,2)=2^{|E|}$ & $T_G(1,2)$ &  $T_G(0,2)$ \\ \hline
 Forest & $T_G(2,1)$ & $T_G(1,1)$ &  $T_G(0,1)$ \\ \hline
Internal & $T_G(2,0)$ & $T_G(1,0)$ &  $T_G(0,0)=0$ \\ \hline
\end{tabular}
\end{center}
\end{prop}

\dem Let $T$ be a spanning tree with $\mI(T)$ internal and $\mE(T)$ external  $(\mG,T)$-active edges. By Lemma \ref{lem:keylink}, the connected subgraphs in $[T^-,T^+]$ are obtained by adding some external $(\mG,T)$-active edges to $T$. Hence, there are $1^{\mI(T)}2^{\mE(T)}$ connected subgraphs in $[T^-,T^+]$. Thus, given the partition of the set of subgraphs into tree-intervals given by Theorem \ref{thm:partition}, the graph $\mG$ has 
$$\displaystyle \sum_{T\textrm{ spanning tree} }1^{\mI(T)}2^{\mE(T)}$$
connected subgraphs. This sum is equal to $T_G(1,2)$ by the characterization \Ref{eq:Tutte-embedded} of the Tutte polynomial. Observe that there are $0^{\mI(T)}2^{\mE(T)}$ external (connected) subgraphs in the interval $[T^-,T^+]$\footnote{Here, as everywhere in this paper, the convention is that $0^0=1$.}. Hence there are $T_G(0,2)$ external subgraphs of $G$. Every other category admits a similar treatment.
\findem

In the next section we will define a bijection $\Phi$ between subgraphs and orientations. In the following one we will study how $\Phi$ specializes to each of the  families of subgraphs defined by the criteria \emph{forest}, \emph{internal}, \emph{connected}, \emph{external} and deduce from it an interpretation for each of the evaluations $T_G(i,j),0\leq i,j\leq 2$ of the Tutte polynomial in terms of orientations.\\

\section{A BIJECTION BETWEEN SUBGRAPHS AND ORIENTATIONS}\label{section:subgraphs-orientations}

\begin{figure}[ht!]
\begin{center}
\input{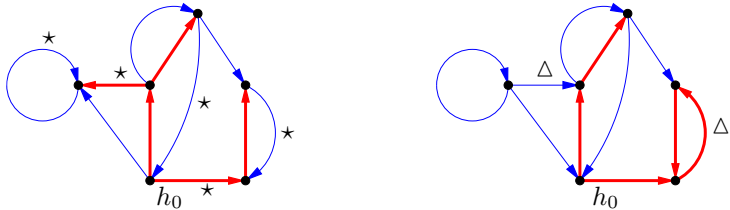}
\caption{Right. The orientation $\mO_T$ associated with a spanning tree $T$. The active edges are indicated by a $\star$. Left. The orientation $\mO_S$ associated with a subgraph $S$ in $[T^-,T^+]$. The edges in the symmetric difference $S\vartriangle T$ are indicated by a $\vartriangle$. }\label{fig:canonical-orientation-subset}
\end{center}
\end{figure}

In this section we define a bijection $\Phi$ between subgraphs and orientations. 
The bijection $\Phi$ is an extension of the correspondence $T\mapsto \mO_T$ between spanning trees and orientations defined in Section~\ref{section:glimpse}. For instance, the image by $\Phi$ of the spanning tree $T$ and the image of a subgraph $S$ in $[T^-,T^+]$ are shown in Figure \ref{fig:canonical-orientation-subset}.

\begin{Def}\label{def:Phi}
Let $\mG$ be an embedded graph. Let $T$ be a spanning tree and let $S$ be a subgraph in the tree-interval $[T^-,T^+]$.
The orientation $\mO_S=\Phi(S)$ is defined as follows. For any edge $e=\{h_1,h_2\}$ with $h_1<h_2$ (for the $(\mG,T)$-order), the arc $\mO_S(e)$ is $(h_1,h_2)$ if and only if - either $e$ is in $T$ and its fundamental cocycle contains no edge in the symmetric difference $S\vartriangle T$ - or if $e$ is not in $T$ and its fundamental cycle contains some edges in  $S\vartriangle T$; the arc $\mO_S(e)$ is $(h_2,h_1)$ otherwise.  
\end{Def}

Recall that a subgraph $S$ is in the tree-interval $[T^-,T^+]$ if and only if every edge in  the symmetric difference  $S \vartriangle T$ is $(\mG,T)$-active. Let $S$ be a subgraph in $[T^-,T^+]$ and let $e$ be any edge of $\mG$. We say that the arc $\mO_S(e)$ is \emph{reverse} if $\mO_S(e)\neq \mO_T(e)$. Observe that the arc $\mO_S(e)$  is reverse  if and only if the  fundamental cycle or cocycle of $e$ (with respect to the spanning tree $T$) contains an edge of $S \vartriangle T$ (compare for instance the orientations $\mO_S$ and $\mO_T$ in Figure \ref{fig:canonical-orientation-subset}). In particular, Definition \ref{def:Phi} of the mapping $\Phi$ extends the definition \Ref{def:O_T} given for spanning trees in Section~\ref{section:glimpse}.

\begin{figure}[ht!]
\begin{center}
\input{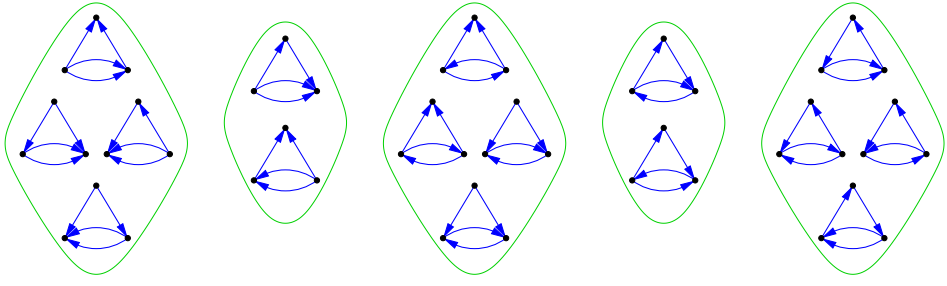}
\caption{The image by $\Phi$ of the subgraphs in Figure \ref{fig:interval-exp}.} \label{fig:Phi-bijection-exp}
\end{center}
\end{figure}

The main result of this section is that the mapping $\Phi$ is a bijection between subgraphs and orientations. For instance, we have represented in Figure \ref{fig:Phi-bijection-exp} the image by $\Phi$ of the subgraphs represented in Figure \ref{fig:interval-exp}.

\begin{thm}\label{thm:Phi-bijection}
Let $\mG$ be an embedded graph. The mapping $\Phi$ establishes a bijection between the subgraphs and the orientations of $G$.
\end{thm}

In order to prove Theorem \ref{thm:Phi-bijection}, we define a mapping $\Psi$ from orientations to subgraphs. We shall prove that $\Psi$ is the inverse of $\Phi$.

\begin{Def}\label{def:Psi-subgraphs-orientations}
Let  $\mG$ be an embedded graph and let $\mO$ be an orientation. We define the subgraph $S=\Psi(\mO)$ by the procedure described below. The procedure $\Psi$ visits the half-edges in sequential order. The set of visited edges is denoted by~$F$. If $C$ is a set of edges that intersects the set $F$ of visited edges, we denote by $e_\fir(C)$ and $h_\fir(C)$ the first visited edge and half-edge of $C$ respectively ($e_\fir(C)$ contains $h_\fir(C)$). In this case, $C$ is said to be \emph{tail-first} if $h_\fir(C)$  is a tail  and \emph{head-first} otherwise.

\noindent \textbf{Initialization:} Initialize the \emph{current half-edge} $h$ to be the root $h_0$. Initialize the subgraph $S$, the tree $T$ and the set of visited edges $F$ to be empty.\\
\textbf{Core:} Do: \\
\textbf{C1:} If the edge $e$ containing $h$ is not in $F$, then decide whether to add $e$ to $S$ and~$T$:
\begin{itemize}
\item If $h$ is a tail, then\\
$(a)$ If $e$ is in a directed cycle $C\subseteq \B{F}$, then add $e$ to $S$ but not to $T$.\\
$(b)$ If $e$ is in a head-first directed cocycle  $D\nsubseteq \B{F}$ such that for all directed cocycle $D'$ with  $e_{\fir}(D')=e_{\fir}(D)$ either $e\in D'$ or ($D\vartriangle D'\nsubseteq \B{F}$ and $e_{\fir}(D\vartriangle D')\in D'$), then do not add $e$ to $S$ nor to $T$.\\
$(c)$ Else, add $e$ to $S$ and to $T$.\\
\item If $h$ is a head, then\\
$(a')$ If $e$ is in a directed cocycle $D\subseteq \B{F}$, then add $e$ to $T$ but not to $S$.\\
$(b')$ If $e$ is in a tail-first directed cycle  $C\nsubseteq \B{F}$ such that for all directed cycle $C'$ with  $e_{\fir}(C')=e_{\fir}(C)$ either $e\in C'$ or ($C\vartriangle C'\nsubseteq \B{F}$ and $e_{\fir}(C\vartriangle C')\in C'$), then add $e$ to $S$ and to $T$.\\
$(c')$ Else, do not add $e$ to $S$ nor to $T$.
\end{itemize}
\indent \indent  Add $e$ to $F$.\vspace{.0cm}\\
\textbf{C2:}  Move to the next half-edge around $T$: \\
\indent \indent If $e$ is in $T$, then set the current half-edge $h$ to be $\sigma\alpha(h)$, else set it to be $\sigma(h)$.\vspace{.0cm}\\
Repeat until the current half-edge $h$ is $h_0$.\vspace{.0cm}\\
\textbf{End:} Return the subgraph $S$.
\end{Def}

In the procedure $\Psi$ the conditions $(a)$ and $(b)$ (resp. $(a')$ and $(b')$) are incompatible. Indeed the following lemma is a classical result of graph theory \cite{Minty:cycle-cocycle}.

\begin{lemma}[\cite{Minty:cycle-cocycle}]\label{lem:minty} 
Every arc (of an oriented graph) is either in a directed cycle or a directed cocycle but not both.
\end{lemma}

\dem (Hint) is the origin of the arc reachable from its end? \findem

We are now going to prove that $\Phi$ and $\Psi$ are inverse mappings.
 
\begin{prop}\label{thm:Psi-Phi-Id}
Let $\mG$ be an embedded graph and let $S$ be a subgraph. The mapping $\Psi$ is well defined on the orientation  $\Phi(S)$ (the procedure terminates) and  ${\Psi\circ\Phi(S)=S}$.
\end{prop}

Proposition \ref{thm:Psi-Phi-Id} implies that the mapping $\Phi$ is injective. Since there are as many subgraphs and orientations ($2^{|E|}$), it implies that $\Phi$ is bijective and that $\Psi$ and $\Phi$ are reverse mappings. The rest of this section is devoted to the proof of proposition~\ref{thm:Psi-Phi-Id}. Observe that $\Psi$ is a variation on the procedure \textbf{Construct-tree} presented in Section \ref{section:glimpse}. The difference lies in the extra Conditions $(a)$, $(b)$, $(a')$, $(b')$ which are now needed in order to cope with reverse edges. In Lemmas \ref{lem:fundamental-is-directed-subgraph} to \ref{lem:fundamental-is-tight} we express some properties characterizing reverse edges. \\

\begin{figure}[ht!]
\begin{center}
\input{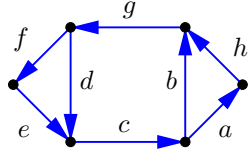}
\caption{The directed cycles $(a,h,g,f,e,c)$ and $(b,g,f,e,c)$ are tight whereas $(a,h,g,d,c)$ is not.} \label{fig:tight-or-not}
\end{center}
\end{figure}

We first need some definitions. Let $\mG$ be an embedded graph and $\mO$ be an orientation. Suppose that the edges and half-edges of $\mG$ are linearly ordered. For any set of edges $C$, we denote by $e_{\min}(C)$ and $h_{\min}(C)$ the minimal edge and half-edge of $C$ respectively. We say that $C$ is \emph{tail-min} if $h_{\min}(C)$ is a tail and \emph{head-min} otherwise. A directed cycle (resp. cocycle) is \emph{tight} if any directed cycle (resp. cocycle) $C'\neq C$ with  $e_{\min}(C')=e_{\min}(C)$ satisfies  $e_{\min}(C\vartriangle C')\in C'$. For instance, if the edges of the graph in Figure \ref{fig:tight-or-not} are ordered by $a<b<c<d<e<f<g$, the directed cycles $(a,h,g,f,e,c)$ and $(b,g,f,e,c)$ are tight whereas $(a,h,g,d,c)$  is not. \\

In Lemmas   \ref{lem:fundamental-is-directed-subgraph} to \ref{lem:fundamental-is-tight} we consider an embedded graph $\mG$, a spanning tree $T$ and a subgraph $S$ in the tree-interval $[T^-,T^+]$. We consider the orientation $\mO_S=\Phi(S)$ and compare edges and half-edges according to the $(\mG,T)$-order. \\

\begin{lemma}\label{lem:fundamental-is-directed-subgraph}
The fundamental cycle (resp. cocycle) of any edge in $S\cap \B{T}$ (resp. $\B{S}\cap T$) is $\mO_S$-directed and tail-min (resp. head-min).
\end{lemma}

\dem If $e$ is in $S\cap \B{T}$ (resp. $\B{S}\cap T$), then every edge $e'$ in its fundamental cycle (resp. cocycle) $C$ is reverse ($\mO_S(e')\neq\mO_T(e')$). By Lemma \ref{lem:fundamental-is-directed}, the cycle (resp. cocycle) $C$ is $\mO_T$-directed, hence it is $\mO_S$-directed. Since $e$ is $(\mG,T)$-active, the minimal edge $e_{\min}(C)$ is $e$. Hence,  $h_{\min}(C)$ is the least half-edge of $e$. By definition of $\mO_S$, the least half-edge of $\mO_S(e)$ is a tail (resp. head). Hence, $C$ is tail-min (resp. head-min). 
\findem

\begin{lemma}\label{lem:reverse-S-cyclic}
Let $e$ be a reverse edge ($\mO_S(e)\neq \mO_T(e)$). Then, $e$ is in $S$ if an only if it is in a directed cycle (otherwise it is in a directed cocycle by Lemma  \ref{lem:minty}). 
\end{lemma}

\dem \\
\ite Suppose that $e$ is in $S$. We want to prove that $e$ is in a directed cycle. If $e$ is in $S\cap \B{T}$, its fundamental cycle is directed by Lemma \ref{lem:fundamental-is-directed-subgraph}. If $e$ is in $S\cap T$ there is an edge $e'\in S\cap \B{T}$ in its fundamental cocycle (since $e$ is reverse). Therefore, $e$ is in the fundamental cycle of $e'$ which is directed by Lemma \ref{lem:fundamental-is-directed-subgraph}.\\
\ite A similar argument proves that if $e$ is in $\B{S}$, then it is in a directed cocycle. In this case, $e$ is not in a directed cycle by Lemma  \ref{lem:minty}. 
\findem

We now need to recall a classical result of graph theory (which is closely related to the axioms of oriented matroid theory \cite{Bjorner:oriented-matroids}).

\begin{lemma}[Orthogonality]\label{lem:one-cycle-one-cocycle}
Let $D$ be a cocycle and let $V_1$ and $V_2$ be the connected components after deletion of $D$. If a directed cycle $C$ contains an arc oriented from $V_1$ to $V_2$ then it also contains an arc oriented from $V_2$ to $V_1$. 
\end{lemma}

\begin{figure}[ht!]
\begin{center}
\input{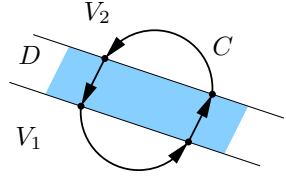}
\caption{A directed cycle crossing a cocycle.} \label{fig:one-cycle-one-cocycle}
\end{center}
\end{figure}

Lemma \ref{lem:one-cycle-one-cocycle} is illustrated by Figure \ref{fig:one-cycle-one-cocycle}.

\begin{lemma}\label{lem:symetric-difference}
An edge $e$ is in $S\cap\B{T}$ (resp. $\B{S}\cap T$)  if and only if it is minimal in a tail-min (resp. head-min) directed cycle (resp. cocycle).
\end{lemma}

\begin{figure}[ht!]
\begin{center}
\input{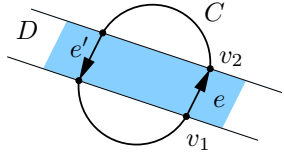}
\caption{The directed cycle $C$, the fundamental cocycle $D$ and the edges $e$ and $e'$.} \label{fig:symetric-difference}
\end{center}
\end{figure}

\dem  We only prove that if an edge is minimal in a tail-min  directed cycle then it is in $\in S\cap\B{T}$. The reverse implication is given by Lemma \ref{lem:fundamental-is-directed-subgraph}. The proof of the dual equivalence ($e$ is minimal in a tail-min  directed cycle if and only if $e$ is in $\B{S}\cap T$) is similar.\\
Let $e=\{h_1,h_2\}$ with $h_1<h_2$ be a minimal edge in a tail-min directed cycle $C$. We want to prove that $e$ is in $S\cap\B{T}$. Observe first that $\mO_S(e)=(h_1,h_2)$ (since $h_{\min}(C)=h_1$  and $C$ is tail-min).   We now prove successively the following points. \\
\iten \emph{The edge $e$ is not in $\B{S}\cap T$}. Otherwise, the edge $e$ would be both in a directed cycle $C$ and in a directed cocycle by Lemma \ref{lem:fundamental-is-directed-subgraph}. \\
\iten \emph{The edge $e$ is not in $S\cap T$}. Suppose the contrary. Since $e$ is in $T$, the arc $\mO_S(e)=(h_1,h_2)=\mO_T(e)$ is not reverse. Let $D$ be the fundamental cocycle of $e$. Let $v_1$ and $v_2$ be the endpoints of $h_1$ and $h_2$ respectively and let $V_2$ be set of descendants of $v_2$. Recall that $v_1$ is the father of $v_2$ in $T$ (Lemma \ref{lemma:between-h1-and-h2}) and that $D$ is the cocycle defined by $V_2$. Since the cycle $C$ is directed and the arc $\mO_S(e)$ in $C\cap D$ is directed toward $V_2$, there is an edge $e'$ in  $C\cap D$ with $\mO_S(e')$ directed away from $V_2$ by Lemma \ref{lem:one-cycle-one-cocycle}. This situation is represented in Figure \ref{fig:symetric-difference}. Since $e$ is minimal in the cycle $C$, we have $e<e'$. Therefore, the arc $\mO_T(e')$ is directed toward $V_2$ by Lemma \ref{lemma:cycle-h1-h1-h2-h2}. Thus, $e'$ is reverse. The edge $e'$ is reverse and contained in a directed cycle, therefore it is in $S$ by Lemma \ref{lem:reverse-S-cyclic}. We have shown that $e'$ is in $S\cap \B{T}$. But this is impossible since $e<e'$ is in the fundamental cycle of $e'$.\\
\iten \emph{The edge $e$ is in $S\cap\B{T}$}. We know from the preceding points that $e$ is in  $\B{T}$. Hence, $\mO_T(e)=(h_2,h_1)\neq \mO_S(e)$. Thus,  $e$ is reverse in  a directed cycle. Therefore, $e$ is in  $S$ by Lemma \ref{lem:reverse-S-cyclic}.
\findem

\begin{lemma}\label{lem:fundamental-is-tight}
The fundamental cycle (resp. cocycle) of any edge in $S\cap \B{T}$ (resp. $\B{S}\cap T$) is tight.
\end{lemma}

\begin{figure}[ht!]
\begin{center}
\input{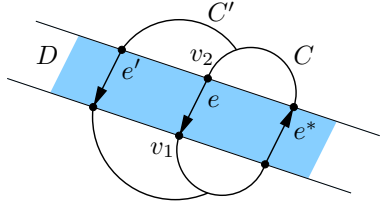}
\caption{The directed cycles $C$ and $C'$ and the cocycle $D$.} \label{fig:fundamental-is-tight}
\end{center}
\end{figure}

\dem We prove that the fundamental cycle of an edge in $S\cap \B{T}$ is tight. The proof of the dual property (concerning edges in $\B{S}\cap T$) is similar. Let $e^*$ be in $S\cap \B{T}$. Recall that $e^*=e_{\min}(C)$. By Lemma \ref{lem:fundamental-is-directed-subgraph}, the fundamental cycle $C$ of $e^*$  is directed.  We want to prove that $C$ is tight. Suppose not and consider a directed cycle $C'$ with $e_{\min}(C')=e_{\min}(C)=e^*$ and $e=e_{\min}(C\vartriangle C')\in C$. The edge $e$ is in the fundamental cycle $C$ of $e^*$, hence $e^*$ is in  fundamental cocycle $D$ of $e$. This situation is represented in Figure \ref{fig:fundamental-is-tight}. Let $v_1$ and $v_2$ be the endpoints of $e$ with $v_1$ father of $v_2$ in $T$. Let $V_2$ be the set of descendants of $v_2$. Recall that $D$ is the cocycle defined by $V_2$. The edge $e$ is in the fundamental cycle of $e^*$ which is $(\mG,T)$ active, hence $e^*<e$. Therefore, the arc $\mO_T(e^*)$ is directed away from $V_2$ by Lemma  \ref{lemma:cycle-h1-h1-h2-h2}. Since $e^*$ is in $S\cap \B{T}$, the arc $\mO_S(e^*)$ is reverse, hence is directed toward $V_2$. Since the cycle $C'$ is directed and the arc $\mO(e^*)$ in $C'\cap D$ is directed toward $V_2$, there is an arc $\mO_S(e')$ in $C'\cap D$ oriented away from $V_2$  by Lemma \ref{lem:one-cycle-one-cocycle}.  Observe that $e'$ is not in the fundamental cycle $C$ since $C\subseteq T+e^*$ and $D\subseteq \B{T}+e$.  Thus, $e'$ is in $C\vartriangle C'$ and $e'>e$. Hence, by Lemma \ref{lemma:cycle-h1-h1-h2-h2}, the arc $\mO_T(e')$ in the fundamental cocycle $D$ of $e$ is directed toward $V_2$. Thus, the arc  $\mO_S(e')\neq \mO_T(e')$ is reverse. Since $e'$ is reverse and contained in a directed cycle, it is in $S$ by Lemma \ref{lem:reverse-S-cyclic}. We have shown that $e'$ is in $S \cap\B{T}$. But this is impossible. Indeed $e'$ is not $(\mG,T)$-active since its fundamental cycle contains $e$ which is less than $e'$. 
\findem

\noindent \textbf{Proof of Proposition \ref{thm:Psi-Phi-Id}.}
We consider a subgraph $S_0$ in the tree-interval $[T_0^-,T_0^+]$ and the orientation $\mO_{S_0}=\Phi(S_0)$.  We want to prove that the procedure $\Psi$ returns the subgraph $S_0$.  We compare edges and half-edges according to the  $(\mG,T_0)$-order denoted by $<$: we say that an edge or half-edge is \emph{greater} or \emph{less} than another.  We also compare edges and half-edges according to their \emph{order of visit} during the algorithm: we say that an edge or half-edge is \emph{before} or \emph{after} another.  
We denote by  $t$ the motion function associated with $T_0$. We denote by $h_i=t^i(h_0)$ the $i^{th}$ half-edge for the $(\mG,T_0)$-order. Also, for every half-edge $h$, we denote $F_h=\{e=\{h_1,h_2\} \textrm{ such that } {\min}(h_1,h_2)< h \}$, $T_h=T_0\cap F_h$ and $S_h=S_0\cap F_h$.\\
We want to prove that \emph{at the beginning of the $i^{th}$ core step, $h=h_i$, $F=F_h$, $T=T_h$, $S=S_h$}, where $h$ is the current half-edge. We proceed by induction on the number of core steps. The property holds for the first $(i=0)$ core step since $h=h_0$ and $F_{h_0}=T_{h_0}=S_{h_0}=\emptyset$. Suppose the property holds for all $i\leq k$. By the induction hypothesis, the $(\mG,T_0)$-order and the order of visit coincide on the edges and half-edges of $F$. In particular, if $C$ is any set not contained in $\B{F}$, then $h_{\min}(C)=h_{\fir}(C)$ and $e_{\min}(C)=e_{\fir}(C)$. Suppose the edge $e$ containing the current half-edge $h$ is not in $F=F_h$. In this case, the current half-edge $h$ (resp. edge $e$) is less than any other half-edge (resp. edge) in $\B{F}$. We consider the different cases $(a),(b),(c),(a'),(b'),(c')$. We will prove successively the following properties.
\begin{itemize}
\item \emph{Condition $(a)$ is equivalent to $e\in S_0 \cap \B{T_0}$}.\\
\iten Suppose Condition $(a)$ holds: $h$ is a tail and $e$ is in a directed cycle $C\subseteq \B{F}$. Since, $C\subseteq \B{F}$, the current half-edge $h$ is minimal in $C$. Since $h$ is a tail, the directed cycle $C$ is tail-min. Thus, $e$ is in $S_0 \cap \B{T_0}$ by Lemma \ref{lem:symetric-difference}.\\
\iten Conversely, if $e$ is in $S_0 \cap \B{T_0}$, then $e$ is minimal in a tail-min directed cycle $C$ by Lemma \ref{lem:symetric-difference}. Therefore, $h$ is a tail and $C\subseteq \B{F}$. 
\item \emph{Condition $(a')$ is equivalent to $e\in \B{S_0} \cap T_0$}. \\
The proof is the similar to the proof of the preceding point.
\item \emph{Condition $(b)$ is equivalent to $e\in \B{S_0} \cap \B{T_0}$ and $\mO_{S_0}(e)$ is reverse}.\\
\iten Suppose Condition $(b)$ holds: $h$ is a tail and $e$ is in a head-first directed cocycle $D\nsubseteq \B{F}$ such that for all directed cocycle $D'$ with  $e_{\fir}(D')=e_{\fir}(D)$ either $e\in D'$ or $D\vartriangle D'\nsubseteq \B{F}$ and $e_{\fir}(D\vartriangle D')\in D'$. Since the $(\mG,T_0)$-order and the order of visit coincide on $F$ we have $h_{\min}(D)=h_{\fir}(D)$. Since the cocycle $D$ is head-first, it is tail-min. The edge $e^*:=e_{\min}(D)$ is minimal in a head-min directed cocycle, hence $e^*$ is in $\B{S_0}\cap T_0$ by Lemma~\ref{lem:symetric-difference}. Let $D^*$ be the fundamental cocycle of $e^*$. Recall that $e_{\min}(D^*)=e^*=e_{\min}(D)$ We want to prove that $e$ is in $D^*$. Suppose $e$ is not in $D^*$. By Condition $(b)$, we have $D\vartriangle D^*\nsubseteq \B{F}$ and $e_\fir(D\vartriangle D^*)\in D^*$. But this is impossible since $e_{\min}(D\vartriangle D^*)=e_\fir(D\vartriangle D^*)$ and $D^*$ is tight by Lemma~\ref{lem:fundamental-is-tight}. Thus, $e$ is indeed in the fundamental cocycle $D^*$ of $e^*$.  Since $e^*$ is in $\B{S_0}\cap T_0$, the edge $e$ is in $\B{T_0}$ and also in $\B{S_0}$ by Lemma \ref{thm:fdmtl-cycle-is-in}. Moreover the arc $\mO_{S_0}(e)$ is reverse.\\
\iten Conversely, suppose that $e$ is in $\B{S_0}\cap \B{T_0}$ and that the arc $\mO_{S_0}(e)$ is reverse. The current half-edge $h$ is the least half-edge of $e$. Since $e$ is external,  $h$ is the head of the arc $\mO_{T_0}(e)$ and the tail of the reverse arc $\mO_{S_0}(e)$. Since $\mO_{S_0}(e)$ is reverse, the external edge $e$ is in the fundamental cocycle $D$ of an edge $e^*\in \B{S_0}\cap T_0$.  The cocycle $D$ is head-min, directed and tight by Lemmas \ref{lem:fundamental-is-directed-subgraph} and \ref{lem:fundamental-is-tight}. Since $e^*=e_{\min}(D)$, the edge  $e^*$ is less than $e$. Therefore $e^*$ is before $e$ and $D\nsubseteq \B{F}$. The cocycle $D$ is head-first since $h_\fir(D)=h_{\min}(D)$. Consider any directed cocycle $D'$ such that $e_\fir(D')=e_{\fir}(D)=e^*$ and $e\notin D'$. We want to prove that $D\vartriangle D'\nsubseteq \B{F}$ and $e_{\fir}(D\vartriangle D')\in D'$. Since $D$ is tight, the edge $e'=e_{\min}(D\vartriangle D')$ is in $D'$. Since $e$ is in $D\vartriangle D'$, the edge $e'$ is less than $e$, hence it is in $F$. Therefore, $D\vartriangle D'\nsubseteq \B{F}$ and $e_{\fir}(D\vartriangle D')=e_{\min}(D\vartriangle D')=e'$ is in $D'$. 
\item \emph{Condition $(b')$ is equivalent to $e\in S_0 \cap T_0$ and $\mO_{S_0}(e)$ is  reverse}.\\
The proof is the similar to the proof of the preceding point.
\item \emph{Condition $(c)$ is equivalent to $e\in S_0 \cap T_0$ and is not reverse}.\\
\iten Suppose Condition $(c)$  holds. In this case, Conditions $(a)$, $(a')$, $(b)$, $(b')$ do not hold. Hence (by the preceding points), the edge $e$ is not in $S_0\vartriangle T_0$ and the arc $\mO_{S_0}(e)$ is not reverse. Since $\mO_{S_0}(e)$ is not reverse and the half-edge $h$ (which is the least half-edge of $e$) is a tail, the edge $e$ is in $T_0$. Since $e$ is not in $S_0\vartriangle T_0$, it is in $S_0$.\\
\iten  Conversely, suppose that $e$ is in $S_0 \cap T_0$  and that $\mO_{S_0}(e)$  is not reverse. By the preceding points, none of the conditions  $(a)$, $(a')$, $(b)$, $(b')$ holds. Moreover, the half-edge $h$ (which is the least half-edge of $e$) is a tail.
\item \emph{Condition $(c')$ is equivalent to $e\in \B{S_0} \cap \B{T_0}$ and is not reverse}.\\
The proof is the similar to the proof of the preceding point.
\end{itemize}
By the preceding points, $e$ is added to $S$ (resp. $T$) in the procedure $\Psi$ if and only if $e$ is in $S_0$ (resp. $T_0$). Hence,  the next half-edge will be $t(h)=\sigma\alpha(h)$ if $h$ is in $T_0$ and $\sigma(h)$ otherwise.  Thus, all the properties are satisfied at the beginning of the $(k+1)^{th}$ core step.
\findem

This concludes the proof of Theorem \ref{thm:Phi-bijection}. We have also proved the following property that will be useful in the next section. 
\begin{lemma}\label{lemPsi-visit-in-order}
During the execution of the procedure $\Psi$ on an orientation $\mO$, the half-edges are visited in $(\mG,T)$-order, where $T$ is the spanning tree $\Delta\circ\Psi(\mO)$.\\
\end{lemma}

\section{SPECIALIZATIONS OF THE BIJECTION BETWEEN SUBGRAPHS AND ORIENTATIONS}\label{section:specializations}
In this section we study several restrictions of the bijection $\Phi$ between subgraphs and orientations. More precisely we shall look at the restriction of  $\Phi$ to each family of subgraphs defined by combining the four criteria \emph{forest}, \emph{internal}, \emph{connected}, \emph{external}. In Figure \ref{fig:subgraphs-organized} we organized the subgraphs according to these criteria. We also represented the orientations associated to each subgraph by the mapping~$\Phi$. As Figure \ref{fig:subgraphs-organized} suggests, there are nice correspondence, which it is the goal of this section to explicit, between the properties of the subgraph and the properties of the associated orientations. Recall from  Proposition \ref{thm:counting-specilizations} that the families of subgraphs defined by combining the criteria \emph{forest}, \emph{internal}, \emph{connected}, \emph{external} are counted by the evaluations $T_G(i,j),0\leq i,j\leq 2$ of the Tutte polynomial. By studying the  restriction of  $\Phi$ to each of these families  we shall obtain a combinatorial interpretation for each of the evaluations $T_G(i,j),0\leq i,j\leq 2$ in terms of orientations or outdegree sequences (see Theorem \ref{thm:all-specializations}).\\


\begin{figure}[ht!]
\begin{center}
\input{subgraphs-organized2.pstex_t}
\caption{Subgraphs in each category  defined by the four criteria \emph{forest}, \emph{internal}, \emph{connected}, \emph{external} and the corresponding orientations. The categories goes from the most general to the most constrained from left to right and from up to down. 
The non-connected subgraphs (resp. non-external connected subgraphs, external subgraphs) are in column 1 (resp. 2, 3).  The subgraphs that are not forests (resp. the forests that are not internal, the internal forests) are in line 1 (resp. 2, 3).} \label{fig:subgraphs-organized}
\end{center}
\end{figure}


\subsection{Connected subgraphs and external subgraphs}\label{section:specialization-connected}
In this subsection we study the restriction of $\Phi$ to connected and to external subgraphs. 

\begin{prop}\label{thm:connected=root-connected}
Let $\mG$ be an embedded graph and let $v_0$ be the root-vertex. The orientation $\mO_S$ is $v_0$-connected if and only if the subgraph $S$ is connected.
\end{prop}

\begin{lemma}\label{lem:hmin-incident-mG0}
Let $\mG$ be an embedded graph and let $T$ be a spanning tree. Let $D$ be a cut and let $\mG_0$ be the connected component of $\mG$ containing the root-vertex $v_0$ after $D$ is removed. Then, the half-edge $h_{\min}(D)$ is incident to $\mG_0$. Moreover, every half-edge not in $\mG_0$ is greater than or equal to $h_{\min}(D)$.
\end{lemma}

\dem 
Let  $t$ be the motion function of $T$. If a half-edge $h$ is incident to $\mG_0$ and is not in $D$ then $t(h)$  is incident to $\mG_0$.  
 Since the root $h_0$ is incident to $\mG_0$, the half-edge $h_{\min}(D)$ is also incident to $\mG_0$ and is less than any half-edge not in $\mG_0$. $\square$ 

\begin{lemma}\label{lem:no-head-min=root-connected}
An orientation is $v_0$-connected if and only if it has no head-min directed cocycle.
\end{lemma}

\dem \\
\ite If there is a head-min directed cocycle, this cocycle is directed toward the component containing $v_0$ by Lemma \ref{lem:hmin-incident-mG0}. Therefore, the vertices in the other components are not reachable from $v_0$ and the orientation is not $v_0$-connected.  \\
\ite If the orientation is not $v_0$-connected we consider the cut $D$ defined by the set $V_0$ of vertices reachable from $v_0$. The cut $D$ is directed toward $V_0$, hence is head-min by Lemma \ref{lem:hmin-incident-mG0}. Let $v_1$ be the endpoint of the edge $e=e_{\min}(D)$ that is not in $V_0$. Let $V_1$ be the set of vertices in the connected component containing $v_1$ after the cut $D$ is deleted. The set of edges $D_1$ with one endpoint in $V_0$ and one endpoint in $V_1$ is a cocycle contained in $D$.  Since every edge in $D_1$ is directed away from $V_0$ the cocycle $D_1$ directed. Since $h_{\min}(D_1)=h_{\min}(D)$ is a head, the cocycle $D_1$ is head-min.
\findem

\noindent \textbf{Proof of Proposition \ref{thm:connected=root-connected}.}
Let $S$ be a subgraph in $[T^-,T^+]$. The orientation $\mO_S$ is $v_0$-connected if and only if there is no head-min directed cocycle by Lemma \ref{lem:no-head-min=root-connected}. An edge is in $\B{S}\cap T$ if and only if it is minimal in a head-min directed cocycle by Lemma ~\ref{lem:symetric-difference}. Thus,  $\mO_S$ is $v_0$-connected if and only if $\B{S}\cap T=\emptyset$. And  $\B{S}\cap T=\emptyset$ if and only if $S$ is connected by Lemma \ref{lem:keylink}. 
\findem

We now study the restriction of the bijection $\Phi$ to external subgraphs.
\begin{prop}\label{thm:external=strongly-connected}
Let $\mG$ be an embedded graph and let $S$ be a subgraph. The orientation $\mO_S$ is strongly connected if and only if $S$ is external.
\end{prop}

\begin{lemma}\label{lem:son-reachable}
Let $T$ be a spanning tree and let $e$ be an edge of $T$. Let $u$ and $v$ be the endpoints of $e$ with the convention that $u$ is the father of $v$. For any \emph{connected} subgraph $S$ in $[T^-,T^+]$, the vertex $v$ is $\mO_s$-reachable from its father $u$.
\end{lemma}

\dem
For any connected subgraph $S$ in  $[T^-,T^+]$, the set $\B{S}\cap T$ is empty by Lemma \ref{lem:keylink}. If the fundamental cocycle of the edge $e$ contains no edge of $S\cap \B{T}$, then the arc $\mO_S(e)$ is not reverse. In this case, the arc $\mO_S(e)=\mO_T(e)$ is directed from $u$ to $v$ by Lemma \ref{lemma:between-h1-and-h2}. Suppose now that the fundamental cocycle of $e$ contains an edge $e^*$ of $S\cap \B{T}$. In this case, $e$ is in the fundamental cycle $C^*$ of $e^*$ which is $\mO_S$-directed  by Lemma \ref{lem:fundamental-is-directed-subgraph}. Therefore, the vertex $v$ is $\mO_s$-reachable from~$u$ (and vice-versa).
\findem

\begin{lemma}\label{lem:mincocycle=active}
Let $\mG$ be an embedded graph. Let $T$ be a spanning tree and let $S$ be a \emph{connected} subgraph in $[T^-,T^+]$. An edge $e$ is minimal in an $\mO_S$-directed cocycle if and only if $e$ is an internal $(\mG,T)$-active edge. 
\end{lemma}

\dem Since the subgraph $S$ is connected, the subset $\B{S}\cap T$ is empty by Lemma \ref{lem:keylink} and the orientation $\mO_S$ is $v_0$-connected by Lemma \ref{thm:connected=root-connected}.\\
\ite Suppose that the edge $e$ is an internal $(\mG,T)$-active edge. The edge $e$ is minimal in its fundamental cocycle $D$. We want to prove that $D$ is $\mO_S$-directed. Note first that $e$ is not in  $S\vartriangle T$ (since $e$ is in $T$ and  $\B{S}\cap T=\emptyset$). No other edge of $D$ is in $S\vartriangle T$ since none is $(\mG,T)$-active. Hence, $\mO_S(e)=\mO_T(e)$.  Let $e'\neq e$ be an edge in the fundamental cocycle $D$ of $e$. The fundamental cycle of $e'$ does not contain any edge of $\B{S}\cap T$ since this edge is empty. Hence, $\mO_S(e')=\mO_T(e')$. Thus, the orientations $\mO_S$ and $\mO_T$ coincide on the cocycle $D$. By Lemma \ref{lem:fundamental-is-directed}, the cocycle $D$ is $\mO_T$-directed, hence it is $\mO_S$-directed.\\
\ite Suppose that $e=\{h_1,h_2\}$ with $h_1<h_2$ is minimal in an $\mO_S$-directed cocycle $D$.  We want to prove that  $e$ is an internal $(\mG,T)$-active edge. We prove successively the following properties:\\
\iten \emph{The half-edge $h_1$ is a tail}. Otherwise, the cocycle $D$  is head-min. (This is impossible by Lemma \ref{lem:no-head-min=root-connected} since $\mO_S$ is is $v_0$-connected.)
\iten  \emph{The edge $e$ is in $T$}. If $e$ is not in $T$, then the arc $\mO_S(e)=(h_1,h_2)$ is reverse. Thus, the fundamental cycle $C$ of $e$ contains an edge of $S\vartriangle T$. Since $C\subseteq T+e$ and  $\B{S}\cap T=\emptyset$, the edge $e$ is in $S\cap \B{T}$. Thus, the cycle $C$ is $\mO_S$-directed by Lemma \ref{lem:fundamental-is-directed-subgraph}. This is impossible since $e$ cannot be both is a directed cycle and a directed cocycle. \\
\iten  \emph{The edge $e$ is $(\mG,T)$-active}. Since the edge $e$ is in $T$, the arc $\mO_S(e)=(h_1,h_2)=\mO_T(e)$ is not reverse. Let $v_1$ and $v_2$ be the endpoints of $h_1$ and $h_2$ respectively. Let $\mG_2$ be the connected component of $\mG$ containing $v_2$ once the cocycle $D$ is removed. The arc $\mO_S(e)$ is directed toward $v_2$, thus the cocycle $D$ is directed toward $\mG_2$. By Lemma \ref{lem:son-reachable}, all the descendants of $v_2$ are reachable from $v_2$, hence they are all in $\mG_2$.  Let $e'$ be an edge in the fundamental cocycle $D'$ of $e$. Since one of the endpoints of $e'$ is a descendant of $v_2$, the edge $e'$ is either in $D$ or in $\mG_2$. Since the minimal half-edge $h_1$ of $D$ is not incident to $\mG_2$, every edge in  $D\cup \mG_2$ is greater than or equal to $e$ by Lemma \ref{lem:hmin-incident-mG0}. Thus, $e'$ is greater than $e$. The edge $e$ is minimal in its fundamental cocycle $D$, that is, $e$ is $(\mG,T)$-active.
\findem

\noindent \textbf{Proof of Proposition \ref{thm:external=strongly-connected}.} 
Let $S$ be a subgraph in $[T^-,T^+]$. \\
\ite Suppose that the subgraph $S$ is external. The subgraph $S$ is connected and there is no $(\mG,T)$-active edge, hence there is no $\mO_S$-directed cocycle by Lemma \ref{lem:mincocycle=active}. Thus, the orientation $\mO_S$ is strongly connected.\\
\ite  Suppose that the orientation $\mO_S$ is strongly connected. The subgraph $S$ is connected (since $\mO_S$ is $v_0$-connected) and there is no $\mO_S$-directed cocycle, hence there is no $(\mG,T)$-active edge  by Lemma \ref{lem:mincocycle=active}. Thus, the subgraph $S$ is external.
\findem

\subsection{Forests and internal forests}\label{section:specialization-forest}
In this subsection we study the restriction of the bijection $\Phi$ to forests and to internal subgraphs. \\

Let $\mG$ be an embedded graph and let $\mO$ be an orientation. We compare half-edges according to the $(\mG,T)$-order, where $T=\Delta\circ\Psi(\mO)$. We say that the orientation $\mO$ is \emph{minimal} if there is no tail-min $\mO$-directed cycle. We shall see (Lemma \ref{thm:unique-minimal}) that for any out degree sequence $\delta$ there is a unique minimal $\delta$-orientation.

\begin{prop}\label{thm:forest=minimal}
The orientation $\mO_S$ is minimal if and only if the subgraph $S$ is a forest. 
\end{prop}

\dem Let $T=\Delta(S)$. By Lemma  \ref{lem:symetric-difference}, an edge is in $S\cap \B{T}$ if and only if it is minimal in a tail-min directed cycle. Thus, the orientation $\mO_S$ is minimal if and only if $S\cap \B{T}=\emptyset$. And  $S\cap \B{T}=\emptyset$ if and only if $S$ is a forest by Lemma \ref{lem:keylink}. 
\findem

\begin{prop}\label{thm:internal=acyclic}
The orientation $\mO_S$ is acyclic if and only if the subgraph $S$ is internal.
\end{prop}

\begin{figure}[ht!]
\begin{center}
\input{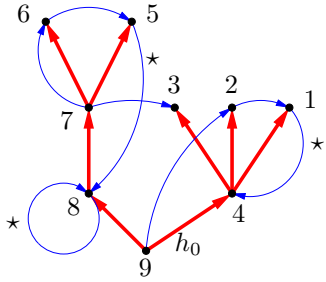}
\caption{A spanning tree $T$, the postfix order, the orientation $\mO_T$ and the external active edges (indicated by a $\star$).} \label{fig:exp-postfix}
\end{center}
\end{figure}

In order to prove Proposition \ref{thm:internal=acyclic} we need to define a linear order,  the \emph{postfix order}, on the vertex set.  For any vertex $v\neq v_0$ we denote by $h_v$ the half-edge incident to $v$ and contained in the edge linking $v$ to its father in $T$. The \emph{postfix order}, denoted by $<_\pos$, is defined by $v<_\pos v_0$ for $v\neq v_0$ and  $v<_\pos v'$ if $h_v<h_{v'}$ for $v,v'\neq v_0$. The postfix order is illustrated in Figure \ref{fig:exp-postfix}.

\begin{lemma}\label{lem:toward-greatest=active}
Let $T$ be a spanning tree and let $e$ be an edge. The arc $\mO_T(e)$ is directed toward its greatest endpoint (for the postfix order) if and only if the edge $e$ is external $(\mG,T)$-active.  
\end{lemma}

Lemma \ref{lem:toward-greatest=active} is illustrated by Figure \ref{fig:exp-postfix}.\\

\dem
Recall from Lemma \ref{lemma:cycle-h1-h1-h2-h2} that a half-edge $h$ is incident to a descendant of $v$ if and only if $h_v'<h\leq h_v$, where $h_v'=\alpha(h_v)$ is the other half of the edge containing $h_v$.  \\ 
\ite Consider an internal edge $e$. Let $u$ and $v$ be the endpoints of $e$ with $u$ father of $v$. By Lemma \ref{lemma:between-h1-and-h2}, the arc $\mO_T(e)$ is directed toward $v$. We want to prove that  $v<_\pos u$.  If  $u=v_0$, the inequality holds. Else, the half-edges  $h_u$ and $h_v$ exist. Moreover, the half-edge $h_v$ is incident to a descendant of $u$, hence $h_v<h_u$ and $v<_\pos u$.\\
\ite Consider an external edge $e$. We write $e=\{h_1,h_2\}$ with $h_1<h_2$ and denote by $u$ and $v$ the endpoints of $h_1$ and $h_2$ respectively. By definition, the arc $\mO_T(e)$ is directed toward $u$. We want to prove that $v\leq_\pos u$ if and only if $e$ is $(\mG,T)$-active.\\
\iten Suppose the edge $e$ is $(\mG,T)$-active. Then, the vertex $v$ is a descendant of $u$ by Lemma \ref{lem:external-active}. The half-edge $h_v$ is incident to a descendant of $u$, hence $h_v\leq h_u$ and $v\leq_\pos u$.\\
\iten Suppose that $v\leq_\pos u$. If $u=v_0$, the vertex $v$ is a descendant of $u$ and the edge $e$ is $(\mG,T)$-active by Lemma \ref{lem:external-active}.  Else, the half-edges $h_u$ and $h_v$ exist and $h_v\leq h_u$. In this case, $\alpha(h_u)<h_1<h_2<h_v\leq h_u$ (indeed, $h_2<h_v$ since $h_2$ is incident to $v$ and $\alpha(h_u)<h_1$ since $h_1$ is incident to $u$), hence $v$ is a descendant of $u$ by Lemma \ref{lemma:cycle-h1-h1-h2-h2}. Thus, the edge $e$ is $(\mG,T)$-active by Lemma \ref{lem:external-active}.
\findem

\noindent \textbf{Proof of Proposition \ref{thm:internal=acyclic}.} Let $S$ be a subgraph in the tree-interval $[T^-,T^+]$. We compare half-edges according to the $(\mG,T)$-order.\\
\ite Suppose that the subgraph $S$ is internal (i.e. the tree $T$ is internal). Recall that  $S\cap \B{T}=\emptyset$. We want to prove that the orientation $\mO_S$ is acyclic. Observe first that the orientation $\mO_T$ is acyclic since the vertices are strictly decreasing (for the postfix order) along any $\mO_T$-directed path by Lemma \ref{lem:toward-greatest=active}. Suppose now that there is an $\mO_S$-directed cycle $C$. The $\mO_S$-directed cycle $C$ contains a reverse arc $\mO(e)$ or $C$ would be $\mO_T$-directed. Since $S\cap \B{T}=\emptyset$, the reverse edges are in the fundamental cocycle of an edge of $\B{S}\cap T$. Thus, the edge $e$ is in the fundamental cocycle $D$ of an edge of $\B{S}\cap T$. The cocycle $D$ is directed by Lemma \ref{lem:fundamental-is-directed-subgraph}. This is impossible since $e$ cannot be both in a directed cycle and in a directed cocycle.\\
\ite  Suppose that the orientation $\mO_S$ is acyclic. We want to prove that the subgraph $S$ is internal (i.e. the tree $T$ is internal). Suppose there is an external $(\mG,T)$-active edge $e$. Let $C$ be the fundamental cycle of $e$. Since  $\mO_S$ is minimal, we know (by Proposition \ref{thm:forest=minimal}) that  $S\cap \B{T}$ is empty. Therefore, the reverse edges are in the fundamental cocycle of an edge of $\B{S}\cap T$. Since $e$ is active, it is not in the fundamental cocycle of an edge of $\B{S}\cap T$. 
Since the other edges of $C$ are not active (they are less than $e$) they are not in $\B{S}\cap T$. Moreover, since they are in $T$, they are not in the fundamental cocycle of an edge of $\B{S}\cap T$.  Thus, the orientations $\mO_S$ and $\mO_T$ coincide on the cycle $C$. By Lemma \ref{lem:fundamental-is-directed}, the cycle $C$ is $\mO_T$-directed, hence it is $\mO_S$-directed. This is impossible since  $\mO_S$ is acyclic.
\findem

\subsection{Minimal orientations and outdegree sequences}\label{section:specialization-outdegree}
In the previous subsection we proved that the bijection $\Phi$ induces a bijection between forests and minimal orientations (Proposition \ref{thm:forest=minimal}). We are now going to link minimal orientations and outdegree sequences.

\begin{prop}\label{thm:unique-minimal}
Let $\mG$ be an embedded graph. For any outdegree sequence $\delta$ there exists a unique minimal $\delta$-orientation. 
\end{prop}

The rest of this subsection is devoted to the proof of  Proposition \ref{thm:unique-minimal}.  We first recall the link between outdegree sequences and the \emph{cycle-flips}. \\ 

Consider an orientation $\mO$ and an  $\mO$-directed cycle (resp. cocycle) $C$.  \emph{Flipping} the $\mO$-directed cycle (resp. cocycle) $C$ means reversing every arc in $C$. We shall talk about \emph{cycle-flips} and \emph{cocycle-flips}.  Observe that flipping a directed cycle does not change the outdegree sequence. Therefore, any orientation $\mO'$ obtained from $\mO$ by a sequence of cycle-flips has the same outdegree sequence as $\mO$. It was proved in \cite{Felsner:lattice} that the converse is also true.
\begin{lemma}\cite{Felsner:lattice}\label{lem:outdegree=sequence-flip}
Two orientations $\mO$ and $\mO'$ have the same outdegree sequence if and only if they can be obtained from one another by a sequence of cycle-flips. Moreover, the flipped cycles can be chosen to be contained in the set $\{e/ \mO(e)\neq \mO'(e)\}$.
\end{lemma}

Lemma \ref{lem:outdegree=sequence-flip} is a direct consequence of the following result proved in \cite{Felsner:lattice}.
\begin{lemma}\cite{Felsner:lattice}\label{lem:same-outdegree-differ-on-a-cycle}
Let $G$ be a graph and let $\mO$ and $\mO'$ be two orientations having the same outdegree sequence. For any edge $e$ in the set $K=\{e'/ \mO(e)\neq \mO'(e)\}$, there is an $\mO$-directed cycle $C\subseteq K$ containing $e$. 
\end{lemma}

\dem (Hint) Start from the end $v$ of $\mO(e)$ and look for an edge $e_1$ in $K$ directed away from $v$. This edge exists except if $v$ is also the origin  of $e$ (since the number of edges directed away from $v$ is the same in $\mO$ and $\mO'$). Repeat the process until arriving to the origin of $e$.
\findem

Recall  that any very arc of an oriented graph is either in a directed cycle or a directed cocycle but not both (Lemma \ref{lem:minty}).
We say that an arc $a$ is \emph{cyclic} or \emph{acyclic} depending on $a$ being in a directed cycle or in a directed cocycle. 
We call \emph{cyclic part} (resp. \emph{acyclic part}) of an orientation the set of cyclic (resp. acyclic) edges. \\

It is well known that the cyclic and acyclic parts are unchanged by a cycle-flip or a cocycle flip \cite{Felsner:lattice,Gioan:enumerating-degree-sequences,Propp:lattice}. Indeed, it is easily seen that the cyclic part of an orientation can only grow when a directed cocycle $D$ is flipped (since no directed cycle intersects with $D$). Since we return to the original orientation by flipping $D$ twice, we conclude that the cyclic and acyclic parts are unchanged by a cocycle-flip. Similarly, the cyclic and acyclic parts are unchanged by a cycle-flip. \\

We will also need the following classical result (closely related to an axioms of oriented matroids theory \cite{Bjorner:oriented-matroids}). 

\begin{lemma}[Elimination]\label{lem:two-cycles-oriented} 
Let  $\mO$ be an orientation and let  $C$ and $C'$ be two $\mO$-directed cycles (resp. cocycles). Let $\mO'$ be the orientation obtained from $\mO$ by flipping $C'$. Then, the symmetric difference of $C$ and $C'$ is a union of $\mO'$-directed cycles (resp. cocycles). In particular, any edge in the $\mO$-directed cycle (resp. cocycle) $C$ is in an $\mO'$-directed cycle (resp. cocycle) $C''\subseteq C\cup C'$. 
\end{lemma}

\begin{figure}[ht!]
\begin{center}
\input{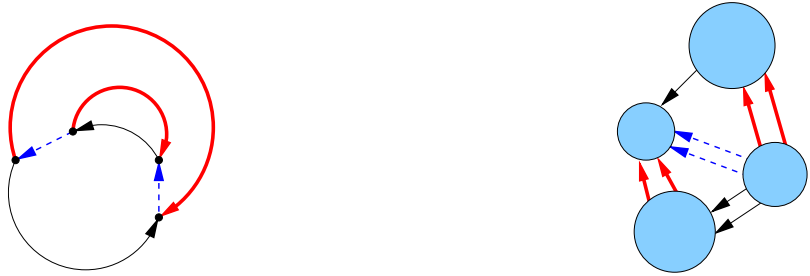}
\caption{The $\mO$-directed cycles (resp. cocycles) $C$ and $C'$ (thin and thick lines) and their intersection (dashed lines).} \label{fig:two-cycles-oriented}
\end{center}
\end{figure}

Lemma \ref{lem:two-cycles-oriented}  is illustrated by Figure \ref{fig:two-cycles-oriented}.

We are now ready to prove Proposition \ref{thm:unique-minimal}. A \emph{false} proof of the uniqueness of the minimal $\delta$-orientation in this proposition is as follows. If there are two different $\delta$-orientations $\mO$ and $\mO'$, then these orientations differ on a directed cycle $C$. Hence, the cycle $C$ is tail-min in either $\mO$ or $\mO'$.  A \emph{false} proof of the existence (of a minimal $\delta$-orientation) is as follows. Take any $\delta$-orientation and starts flipping cycles until no more tail-min directed cycle remains. Of course, both the uniqueness and existence proofs are false in this version since flipping a cycle changes the associated subgraph, hence the spanning tree and the order on the half-edges. However being a bit careful, one can make both proofs correct.\\

We consider the procedure $\Psi$ on orientations (see Definition  \ref{def:Psi-subgraphs-orientations}). For an orientation $\mO$ we denote by $\Psi[\mO]$ the execution of $\Psi$ on $\mO$.  Recall (from Lemma \ref{lemPsi-visit-in-order}) that the half-edges are visited in $(\mG,T)$-order during $\Psi[\mO]$, where $T$ is the spanning tree $\Delta\circ\Psi(\mO)$. Therefore, \emph{the orientation $\mO$ is minimal if and only if Condition $(a)$ never holds during the execution $\Psi[\mO]$}.  \\

\begin{lemma}\label{lem:minimal=outdegree1}
Let $\mO$ be an orientation. 
Consider the current half-edge $h$, the edge $e$ and the sets $F$, $S$ and $T$  at the beginning of a given core step of the execution $\Psi[\mO]$. Let $C_f\subseteq \B{F+e}$ be an $\mO$-directed cycle and let $\mO'$ be the orientation obtained from $\mO$ by flipping $C_f$. We want to prove that \emph{Condition $(a)$ (resp. $(b)$, $(c)$, $(a')$, $(b')$, $(c')$) holds for the orientation $\mO$ if and only if it holds for the orientation $\mO'$}. (Let us insist that when evaluating the Conditions   $(a), \cdots, (c')$ for the orientation $\mO'$, the symbols $F$, $S$, $T$, $h_\fir$ and $e_\fir$ continue to refer to the execution of $\Psi[\mO]$.) 
\end{lemma}

\dem
 Note first that the orientations $\mO$ and $\mO'$ coincide on the current half-edge $h$ since $e\notin C_f$. We now study separately the different conditions.\\
\ite Recall that $\mO$ and $\mO'$ coincide on their acyclic part: the directed cocycles of $\mO$ and $\mO'$ are the same. Therefore, Condition $(b)$ (resp. $(a')$) holds for $\mO$ if and only if it holds for $\mO'$. \\
\ite Suppose now that Condition  $(a)$ holds for $\mO$: the current half-edge $h$ is a tail and the edge $e$ is in an $\mO$-directed cycle $C\subseteq \B{F}$. By Lemma \ref{lem:two-cycles-oriented}, the edge $e$ is also in an $\mO'$-directed cycle $C'\subseteq C\cup C_f\subseteq \B{F}$. Thus,  Condition  $(a)$ holds for $\mO'$. The same argument proves that if Condition $(a)$ holds for $\mO'$, then it holds for $\mO$ ($\mO$ is obtained from $\mO'$ by flipping the $\mO'$-directed cycle $C_f$).\\
\ite  Suppose now that Condition  $(b')$ holds for $\mO$: the current half-edge $h$ is a head and the edge $e$ is in a tail-first $\mO$-directed cycle $C\nsubseteq \B{F}$ such that for all $\mO$-directed cycle $C'$ with  $e_{\fir}(C')=e_{\fir}(C)$ either $e\in C'$ or ($C\vartriangle C'\nsubseteq \B{F}$ and $e_{\fir}(C\vartriangle C')\in C'$). By Lemma  \ref{lem:two-cycles-oriented},  the edge $e^*=e_{\fir}(C)$ is in an $\mO'$-directed cycle $C_1\subseteq C\cup C_f$. Note that $e_{\fir}(C_1)=e^*$. We want to prove that Condition $(b')$  holds for $\mO'$ by considering the  $\mO'$-directed cycle $C_1$. We prove successively the following properties.
\begin{itemize}
\item \emph{The edge $e$ is in $C_1$}. \\
The edge $e^*$ is in the $\mO'$-directed cycle $C_1$ and not in $C_f$. By Lemma  \ref{lem:two-cycles-oriented}, there is an $\mO$-directed cycle $C_2\subseteq C_1\cup C_f$ containing $e^*$ (since $\mO$ is obtained from $\mO'$ by flipping $C_f$). Note that  $e_{\fir}(C_2)=e^*$. Suppose that $e$ is not in $C_2$. By Condition $(b')$ on $C$,  we have $C\vartriangle C_2 \nsubseteq \B{F}$ and $e_{\fir}(C\vartriangle C_2)\in C_2$.  This is impossible since $\B{C}\cap C_2\subseteq C_f$ (since $C_2\subseteq C_1 \cup C_f\subseteq C \cup C_f$) and the edge $e$ in $C\cap\B{C_2}$ is visited before any edge in $C_f$. Thus $e\in C_2$. Since  $e\in C_2\subseteq C_1 \cup C_f$ and $e$ is not in $C_f$, it is in $C_1$.  
\item \emph{For all $\mO'$-directed cycle $C_1'$ with  $e_{\fir}(C_1')=e_{\fir}(C_1)$ either $e\in C_1'$ or  ($C_1\vartriangle C_1'\nsubseteq \B{F}$ and $e_{\fir}(C_1 \vartriangle C_1')\in C_1'$}). (This proves that Condition $(b')$ is satisfied for $\mO'$). \\
Let $C_1'$ be an $\mO'$-directed cycle not containing $e$ and such that  $e_{\fir}(C_1')=e_{\fir}(C_1)=e^*$. We want to prove that $C_1\vartriangle C_1'\subseteq \B{F}$ and $e_{\fir}(C_1 \vartriangle C_1')\in C_1'$. The edge $e^*$ is in the $\mO'$-directed cycle $C_1'$ but not in $C_f$. By Lemma  \ref{lem:two-cycles-oriented}, there exists an $\mO$-directed cycle $C'\subseteq C_1'\cup C_f$ containing $e^*$. Note that $e_{\fir}(C')=e^*$ and that $e\notin C'$ (since $e$ is not in $C_f$ nor in $C_1'$ by hypothesis). By Condition $(b')$ on $C$, we have $C\vartriangle C'\nsubseteq \B{F}$  and $e^\vartriangle=e_{\fir}(C \vartriangle C')\in C'$. We now prove the following properties.\\
\iten \emph{The edge $e^\vartriangle$ is in $\B{C_1}\cap C_1'$. }\\
 The edge $e^\vartriangle$ is in $C_1'$ since $e^\vartriangle\notin C_f$ and  $e^\vartriangle\in C'\subseteq C_1'\cup C_f$. Moreover, $e^\vartriangle$ is not in $C_1$ since $e^\vartriangle\notin C$,  $e^\vartriangle\notin C_f$ and $C_1\subseteq C\cup C_f$. Thus, $e^\vartriangle$ is in $\B{C_1}\cap C_1'$. \\
\iten \emph{Any edge in $C_1\cap \B{C_1'}$ is visited after $e^\vartriangle$ during the execution $\Psi[\mO]$.}\\
Let $e'$ be an edge in $C_1\cap \B{C_1'}$. If $e'$ is in $C_f$, it  is visited  after $e^\vartriangle$. Else,  $e'$ is in $C$ since $e' \in C_1$, $e' \notin C_f$ and  $C_1\subseteq C\cup C_f$. Moreover, $e'$ is not in $C'$ since $e'\notin C_1'$,  $e'\notin C_f$ and  $C'\subseteq C_1\cup C_f$.   Since $e' \in  C\vartriangle C'$, the edge $e'$ is visited after $e^\vartriangle=e_\fir(C\vartriangle C')$ during the execution $\Psi[\mO]$.\\
Since  $e^\vartriangle$ is in $\B{C_1}\cap C_1'$ and any edge  in $C_1\cap \B{C_1'}$ is visited after $e^\vartriangle$, the edge $e_{\fir}(C_1 \vartriangle C_1')$ is in $C_1'$. Thus, Condition $(b')$ holds for $\mO'$.
\end{itemize}
We have proved that if Condition $(b')$ holds for $\mO$, then it holds for $\mO'$. The same argument proves that if Condition $(b')$ holds for $\mO'$, then it holds for $\mO$.\\
\ite Condition  $(c)$ holds for $\mO$ if $h$ is a tail and Conditions $(a)$ and $(b)$ do not hold for $\mO$  By the preceding points this is true if and only if $h$ is a tail and Conditions $(a)$ and $(b)$ do not hold for $\mO'$. Therefore, Condition $(c)$ holds for $\mO$ if and only if it holds for $\mO'$. Similarly, Condition  $(c')$ holds for $\mO$ if and only if it holds for $\mO'$. 
\findem


\begin{lemma}\label{lem:minimal=outdegree2}
Consider two orientations $\mO$ and $\mO'$ having the same outdegree sequence.  We consider the executions  $\Psi[\mO]$ and $\Psi[\mO']$. 
For all $0\leq i<|H|$, we denote by $h_i$, $F_i$, $T_i$ and $S_i$  the current half-edge and the sets $F$, $T$ and $S$ at the beginning of the $i^{th}$ core step of the execution $\Psi[\mO]$ (see Definition  \ref{def:Psi-subgraphs-orientations}). We define $h_i'$, $F_i'$, $T_i'$ and $S_i'$ similarly for the orientation $\mO'$. We want to prove that \emph{if the orientations $\mO$ and $\mO'$  coincide on $h_i$ for all $i<k$ (that is, $\mO(e_i)=\mO'(e_i)$ where $e_i$ is the edge containing $h_i$), then the $k$ first core steps of the executions  $\Psi[\mO]$ and $\Psi[\mO']$ are the same. In particular,  $h_i=h_i'$, $F_i=F_i'$, $S_i=S_i'$, and $T_i=T_i'$ for all  $i\leq k$}.
\end{lemma}

\dem
We proceed by induction on $k$. Recall from Lemma \ref{lem:outdegree=sequence-flip} that the orientation $\mO'$ can be obtained from $\mO$ by a sequence of cycle-flips such that the flipped cycles are contained in the set $K=\{e/ \mO(e)\neq \mO'(e)\}$.  For $k=0$ the property obviously holds. Now suppose that the property holds for $k$ and suppose that $\mO$ and $\mO'$  coincide on $h_i,i<k+1$. By the induction hypothesis the current half-edge $h_k=h_k'$ and the sets $F=F_k=F_k'$, $S=S_k=S_k'$, and $T=T_k=T_k'$ are the same at the  beginning of the $(k+1)^{th}$ core step of the procedures $\Psi[\mO]$ and $\Psi[\mO']$. Moreover, the set $K=\{e'/ \mO(e')\neq \mO'(e')\}$ of reverse edges is contained in $\B{F+e}$. Since $\mO'$ is obtained from $\mO$ by a sequence of flips of cycles contained in $\B{F+e}$, we know by induction on Lemma \ref{lem:minimal=outdegree1} that Condition $(a)$ (resp. $(b)$, $(c)$, $(a')$, $(b')$, $(c')$) holds for the orientation $\mO$ if and only if it holds for the orientation $\mO'$. Therefore, the $(k+1)^{th}$ core step is the same for the two executions $\Psi[\mO]$ and $\Psi[\mO']$. In particular, the sets $F$, $S$, and $T$ are modified in the same way in both executions and  $h_{k+1}=h_{k+1}'$. Thus, the property holds by induction. 
\findem

\noindent \textbf{Proof of Proposition \ref{thm:unique-minimal}.} 
Recall that an orientation $\mO$ is minimal if and only if Condition $(a)$ never holds during the execution $\Psi[\mO]$. Thus, we need to prove that for any outdegree sequence $\delta$ there exists a unique $\delta$-orientation $\mO$ such that Condition $(a)$ never holds during the execution $\Psi[\mO]$.\\
\ite \textbf{Uniqueness:}  Let $\mO$ and $\mO'$ be two (distinct) orientations having the same outdegree sequence. We take the same notations $h_i$, $F_i$, $T_i$, $S_i$, $h_i'$, $F_i'$, $T_i'$, $S_i'$ as in Lemma \ref{lem:minimal=outdegree2}. 
Let $k$ be the first index such that $\mO$ and $\mO'$ differ on $h_k$. By Lemma \ref{lem:minimal=outdegree2}, we have $h_k=h_k'$ and $F_k=F_k'$, $T_k=T_k'$, $S_k=S_k'$. We can suppose without loss of generality that $h_k$ is a tail in  $\mO$ and a head in $\mO'$. We now prove that Condition $(a)$ holds for $\mO$. By hypothesis, the edge $e$ containing $h$ is such that $\mO(e)\neq \mO'(e)$. Hence, by Lemma \ref{lem:same-outdegree-differ-on-a-cycle}, the edge $e$ is contained in an $\mO$-directed cycle $C\subseteq K=\{e/ \mO(e)\neq \mO'(e)\}$. Since  $\mO$ and $\mO'$ coincide on $h_i$ for $i<k$, the set $K$ is contained in $\B{F_i}$. Since  $C\subseteq\B{F_i}$ is $\mO$-directed, Condition $(a)$ holds for $\mO$.\\
\ite \textbf{Existence:} Let $\delta$ be an outdegree sequence. We want to find a $\delta$-orientation $\mO$ such that Condition $(a)$ never holds during the execution $\Psi[\mO]$. Let $\mO_0$ be any $\delta$-orientation. We are going to define a set of $\delta$-orientations $\mO_0,\mO_1,\ldots,\mO_{|H|}$ such that Condition $(a)$ is not satisfied during the $i$ first core steps of the execution $\Psi[\mO_i]$. We prove that $\mO_k$ exists by induction on $k$. Suppose the $\delta$-orientation $\mO_{k-1}$ exists. We consider the current half-edge $h$, the edge $e$ and the sets $F$, $S$ and $T$  at the beginning of the $k^{th}$  core step of $\Psi[\mO_{k-1}]$. If either $e\in F$ or Condition $(a)$ does not hold, we define $\mO_{k}=\mO_{k-1}$. Else, the current half-edge $h_k$  is a tail (for the orientation $\mO_{k-1}$) and there is an $\mO_k$-directed cycle $C\subseteq \B{F}$ containing $e$. In this case, we define $\mO_k$ to be the orientation obtained from $\mO_{k-1}$ by flipping the cycle $C$. Observe that $\mO_{k}$ is a $\delta$-orientation in which $h_k$ is a head. Moreover, since $C\subseteq \B{F}$  the two orientations $\mO_{k-1}$ and $\mO_{k}$  coincide on the half-edges $h_i$ for $i<k$, where  $h_i$ is the current half-edge at the beginning of the $i^{th}$ core step of the execution $\Psi[\mO_{k-1}]$. Thus, by Lemma \ref{lem:minimal=outdegree2}, the $k$ first core steps of the executions $\Psi[\mO_{k-1}]$ and $\Psi[\mO_{k}]$ are the same. Moreover, the current half-edge  $h=h_k$ at the beginning of the $k^{th}$ core step of $\Psi[\mO_{k}]$ is a head (for the orientation $\mO_{k}$). Hence, Condition $(a)$ does not hold at this core step. Thus, $\mO_k$ is a $\delta$-orientation such that Condition $(a)$ does not hold during the $k^{th}$ first core steps of the execution $\Psi[\mO_{k}]$. The orientations  $\mO_0,\mO_1,\ldots,\mO_{|H|}$ exist by induction. In particular, the $\delta$-orientation $\mO_{|H|}$ is such that Condition $(a)$ never holds during the execution $\Psi[\mO_{|H|}]$. 
\findem

From Proposition \ref{thm:forest=minimal} and  \ref{thm:unique-minimal} one obtains the following bijection between outdegree sequences and forests. 
\begin{prop}\label{thm:bij-outdegree}
Let $\mG$ be an embedded graph. The mapping $\Gamma$ which associates with any subgraph $S$ the outdegree sequence of  the orientation $\mO_S$ establishes a bijection between the forests and the outdegree sequences of $\mG$.
\end{prop}
Another bijection between outdegree sequences and forests was established in \cite{Kleitman:forests-score-vectors} after Stanley asked for such a bijection \cite{Stanley:rationnal-polytopes}. \\

\subsection{Summary of the specializations and further refinements}\label{section:all-specializations}
From Propositions  \ref{thm:connected=root-connected}, \ref{thm:external=strongly-connected}, \ref{thm:forest=minimal} and \ref{thm:internal=acyclic} we can characterize the orientations associated with each class of subgraphs defined by the criteria \emph{forest}, \emph{internal}, \emph{connected}, \emph{external}. Each class of subgraphs is counted by a specialization of the Tutte polynomial given in Proposition \ref{thm:counting-specilizations}. Our results are summarized in the following theorem.

\begin{thm}\label{thm:all-specializations}
Let $\mG$ be an embedded graph and let $v_0$ be the root-vertex. 
\begin{enumerate}
\item The $v_0$-connected orientations are in bijection with the connected subgraphs counted by $T_G(1,2)$. 
\item The strongly connected orientations are in bijection with the external subgraphs counted by $T_G(0,2)$.
\item The outdegree sequences are in bijection with minimal orientations, which are in bijection with forests, counted by $T_G(2,1)$.
\item The acyclic orientations are in bijection with internal forests counted by $T_G(2,0)$.
\item The $v_0$-connected outdegree sequences are in bijection with $v_0$-connected minimal orientations which are in bijection with spanning trees counted by $T_G(1,1)$.
\item The strongly connected outdegree sequences  are in bijection with strongly connected minimal orientations which are in bijection with external spanning trees counted by $T_G(0,1)$.
\item The $v_0$-connected acyclic orientations are in bijection with internal spanning trees counted by $T_G(1,0)$.
\end{enumerate}
\end{thm}

Theorem \ref{thm:all-specializations} is illustrated by Figure \ref{fig:subgraphs-organized}. The enumeration of acyclic orientations by $T_G(2,0)$ was first established by Winder in 1966 \cite{Winder:hyperplane-orient-acycliques} and rediscovered by Stanley 1973 \cite{Stanley:acyclic-orientations}. The result of Winder was stated as an enumeration formula for the number of faces of hyperplanes arrangements and was independently extended to reel arrangements by Zaslavsky \cite{Zaslavsky:face-count-hyperplane} and to orientable matroids by Las Vergnas \cite{Vergnas:matroides-orientables}. The enumeration of $v_0$-connected acyclic orientations by $T_G(1,0)$ was found by Greene and Zaslavsky \cite{Greene:interpretation-Tutte-poly}. In \cite{Gessel:Tutte-poly+DFS}, Gessel and Sagan gave a bijective proof of both results. In  \cite{Gebhard-sagan:acyclic-orientations}, Gebhard and Sagan gave three other proofs of Greene and Zaslavsky's result. The enumeration of strongly connected orientations by $T_G(0,2)$ is a direct consequence of Las Vergnas' characterization of the Tutte polynomial \cite{Vergnas:Morphism-matroids-2}. The enumeration of outdegree sequences by $T_G(2,1)$ was discovered by Stanley \cite{Brylawski:Tutte-poly,Stanley:rationnal-polytopes} and a bijective proof was established in \cite{Kleitman:forests-score-vectors}. The enumeration of $v_0$-connected orientations by $T_G(1,2)$, the enumeration of $v_0$-connected outdegree sequences by  $T_G(1,1)$ and the enumeration of strongly connected outdegree sequences by $T_G(0,1)$ were proved by Gioan \cite{Gioan:enumerating-degree-sequences}. \\

\noindent \textbf{Refinements.} It is possible to refine the results of Theorem \ref{thm:all-specializations}. For instance, we have proved that the acyclic orientations of a graph $G$ are counted by $T_G(2,0)$. This is the sum of the coefficients of the polynomial $T_G(1+x,0)$ (which is closely related to the chromatic polynomial of $G$). We denote by $[x^i]P(x)$ the coefficient of $x^i$ in a polynomial $P(x)$. The identities
$$\sum_{i\in \mathbb{N}} [x^i]T_G(1+x,0)=T_G(2,0)=|\{\textrm{acyclic orientations}\}|,$$
and 
$$\sum_{i\in \mathbb{N}} [x^i]T_G(x,0)=T_G(1,0)=|\{\textrm{$v_0$-connected acyclic orientations}\}|,$$
make it appealing to look for a partition of the acyclic orientations (resp. root-connected acyclic orientations) in parts of size $[x^i]T_G(1+x,0)$ (resp.  $[x^i]T_G(x,0)$). Such partitions were defined by Lass in \cite{Lass:interpretation-Tutte-poly} using set functions algebra. The partition defined by Lass is linked to former constructions by Cartier, Foata, Gessel, Stanley and Viennot (see references in \cite{Lass:interpretation-Tutte-poly}). 
More generally, one can try to interpret the coefficients of $T_G(x,1)$, $T_G(1+x,1)$, $T_G(x,2)$, $T_G(1+x,2)$ \emph{etc.} in terms of orientations in order to interpolate between the different specializations $T_G(i,j),0\leq i,j\leq 2$. Observe that the coefficients of each of these polynomials can be given an interpretation in terms of subgraphs. For instance, $[x^i]T_G(1+x,0)$ counts internal forests with $i+1$ trees (by Theorem \ref{thm:partition} and Lemma \ref{lem:keylink}) and $[x^i]T_G(x,0)$ counts internal spanning trees with $i$ internal embedding-active edges (by Theorem \ref{thm:activity-OB}).\\

\begin{figure}[htb!]
\begin{center}
\input{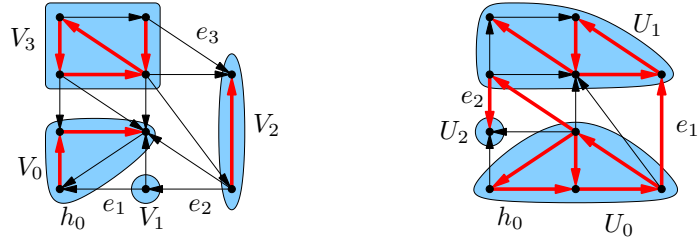}
\caption{Left: root-components of an orientation. Right: root-strong-components of a $v_0$-connected orientation. The thick edges correspond to the subgraph associated with the orientation by the bijection $\Psi$.}\label{fig:root-components}
\end{center}
\end{figure}

We will give an interpretation of the coefficients  $[x^i]T_G(1+x,j)$ for $i\geq 0$ and $j=0,1,2$ in terms of orientations. Let $\mO$ be an orientation. We define the partition of the vertex set $V$ into \emph{root-components} $V=\biguplus_{0\leq i\leq k} V_i$ as follows. The first root-component $V_0$ is the set of vertices reachable from the root-vertex $v_0$. If $W_k=\cup_{0\leq i\leq k} V_i\subsetneq V$, we consider the minimal edge $e_k$ with one vertex in $W_k$ and one vertex $v_k$ in $\B{W_k}$ (the edges are compared according to the $(\mG,T)$-order, where $T=\Delta(\Psi(\mO))$). Then, the $(k+1)^{th}$ root-component is the set of vertices in $\B{W_k}$ that are reachable from $v_k$. For instance, the root-components have been indicated for the orientation in Figure \ref{fig:root-components} (left). It is clear that $v_0$-connected orientations have only one root-component. Given a $v_0$-connected orientation $\mO$, we define the partition of the vertex set $V$ into \emph{root-strong-components} $V=\biguplus_{0\leq i\leq k} U_i$ as follows. The first root-strong-component $U_0$ is the set of vertices that can reach the root-vertex $v_0$. If $W_k=\cup_{0\leq i\leq k} U_i\subsetneq V$, we consider the minimal edge $e_k$ with one vertex in $W_k$ and one vertex $v_k$ in $\B{W_k}$. Then, the $(k+1)^{th}$ root-strong-component is the set of vertices in $\B{W_k}$ that can reach $v_k$.  For instance, the root-strong-components have been indicated for the $v_0$-connected orientation in Figure \ref{fig:root-components} (right).\\

\begin{thm}\label{thm:refinements}
Let $\mG$ be an embedded graph and let $v_0$ be the root-vertex. 
The coefficient $[x^i]T_G(1+x,2)$ (resp. $[x^i]T_G(1+x,1)$, $[x^i]T_G(1+x,0)$) counts orientations (resp. minimal orientations, acyclic orientations) with $i+1$ (non-empty) root-components.  The coefficient $[x^i]T_G(x,2)$ (resp. $[x^i]T_G(x,1)$, $[x^i]T_G(x,0)$) counts $v_0$-connected orientations (resp. minimal $v_0$-connected orientations, acyclic $v_0$-connected orientations) with $i+1$ (non-empty) root-strong-components. \\
\end{thm}

As mentioned above, the coefficients $[x^i]T_G(1+x,0)$ and $[x^i]T_G(x,0)$ had already been interpreted by Lass in \cite{Lass:interpretation-Tutte-poly}. We now prove Theorem \ref{thm:refinements}.\\

\begin{lemma}\label{lem:ei=min-directed-cocycle}
Let $\mG$ be an embedded graph and let $\mO$ be an orientation. We consider the spanning tree  $T=\Delta(\Psi(\mO))$ and compare the half-edges and edges according to the $(\mG,T)$-order. Let $V_0,\dots,V_k$ be the root-components and let $W_i=\cup_{0\leq j\leq i}V_j$. Let $D_i$ for $i=1\ldots k$ be the cut defined by $W_{i-1}$ and let $e_i$ be the minimal edge in $D_i$. Then, an edge is minimal in a head-min directed cocycle if and only if it is in the set $\{e_1,\ldots,e_k\}$. 
\end{lemma}

\dem \\
\ite We first prove that for all $1\leq i \leq k$ \emph{the edge $e_k$ is minimal in a head-min directed cocycle.} Clearly, every edge in the set $D_i$ is directed toward the vertices in $W_{i-1}$.  Let $v_i$ be the endpoint of $e_i=e_{\min}(D)$ which is not in $W_{i-1}$. Let $X_i$ be the set of vertices contained in the connected component containing $v_i$ once the cut $D$ is removed. The set $D$ of edges with one endpoint in $W_{i-1}$ and one endpoint in $X_i$ is a directed cocycle contained in $D_i$. Thus, the edge $e_i$ is minimal in the directed cocycle $D$ directed toward $W_{i-1}$. Since the cocycle $D$ is directed toward the component containing the root-vertex, it is head-min by Lemma \ref{lem:hmin-incident-mG0}.\\
\ite Consider an edge $e$ minimal in a head-min directed cocycle $D$. We want to prove that $e$ is in   $\{e_1,\ldots,e_k\}$. Let $\mG_0$ and $\mG_1$ be the connected components after the cocycle $D$ is removed with the convention that $\mG_0$ contains the root-vertex $v_0$. The directed cocycle $D$ is head-min, hence it is directed toward $\mG_0$ by Lemma \ref{lem:hmin-incident-mG0}. Let $i$ be the first index such that the root-component $V_i$ contains a vertex $v$ of $\mG_1$. The cocycle $D$ is directed toward $\mG_0$, hence no edge of $\mG_1$ is reachable from $v_0$ and the index $i$ is positive. Let $u_i$ and $v_i$ be the endpoints of $e_i$ in $W_{i-1}$ and $\B{W_{i-1}}$ respectively. By definition, the endpoint $u_i$ is in $\mG_0$. Moreover, the vertex $v\in \mG_1$ is reachable from $v_i$, hence the endpoint $v_i$ is in $\mG_1$. Thus, the edge $e_i$ is in $D$ and $e_i \geq e=e_{\min}(D)$. We will now prove that $e_i\leq e$. The subset of vertices $W_{i-1}$ contains the root-vertex and the subset of edges $D_i$ separate $W_{i-1}$ and $\B{W_{i-1}}$, hence every edge with one endpoint in $\B{W_{i-1}}$ is greater than $e_i=e_{\min}(D_i)$ by Lemma \ref{lem:hmin-incident-mG0}. The edge $e$ has one endpoint in $\mG_1\subseteq \B{W_{i-1}}$, hence $e_i\leq e$. Thus, $e=e_i$.
\findem

Here is a counterpart of Lemma \ref{lem:ei=min-directed-cocycle} for root-strong-components.

\begin{lemma}\label{lem:ei'=min-directed-cocycle'}
Let $\mG$ be an embedded graph and let $\mO$ be a $v_0$-connected orientation. We consider the spanning tree  $T=\Delta(\Psi(\mO))$ and compare the half-edges and edges according to the $(\mG,T)$-order. Let $U_0,\dots,U_k$ be the root-strong-components and let $W_i=\cup_{0\leq j\leq i}U_j$. Let $D_i$ for $i=1\ldots k$ be the cut defined by $W_{i-1}$ and let $e_i$ be the minimal edge in  $D_i$. Then, an edge is minimal in a directed cocycle if and only if it is in the set $\{e_1,\ldots,e_k\}$. 
\end{lemma}

\dem The proof of Lemma \ref{lem:ei'=min-directed-cocycle'} very similar to the proof of Lemma \ref{lem:ei=min-directed-cocycle} and is left to the reader.
\findem

\noindent \textbf{Proof of Theorem \ref{thm:refinements}.}\\
\ite We first prove that \emph{the coefficient $[x^i]T_G(1+x,2)$ (resp. $[x^i]T_G(1+x,1)$, $[x^i]T_G(1+x,0)$) counts orientations (resp. minimal orientations, acyclic orientations) with $i+1$ root-components.} 
Let $T$ be a spanning tree with $\mI(T)$ internal and $\mE(T)$ external $(\mG,T)$-active edges. By Lemma \ref{lem:keylink}, the coefficient  $[x^i](1+x)^{\mI(T)}2^{\mE(T)}$ counts the  subgraphs $S$ in the tree-interval $[T^-,T^+]$ having $i$ edges in $\B{S}\cap T$. Given that the tree-intervals form a partition of the set of subgraphs, the coefficient $[x^i]\sum_{T \textrm{ spanning tree }}(1+x)^{\mI(T)}2^{\mE(T)}$ counts the subgraphs $S$ having $i$ edges in $\B{S}\cap \Delta(S)$. Moreover, by the characterization \Ref{eq:Tutte-embedded} of the Tutte polynomial, the sum $\sum_{T}(1+x)^{\mI(T)}2^{\mE(T)}$ is equal to $T_G(1+x,2)$. Similarly,  the coefficient $[x^i]T_G(1+x,1)$ (resp. $[x^i]T_G(1+x,0)$) counts the forests (resp. internal forests) $S$ having $i$ edges in $\B{S}\cap \Delta(S)$. By Theorem \ref{thm:all-specializations} and Lemma \ref{lem:symetric-difference}, the coefficient  $[x^i]T_G(1+x,2)$  (resp.  $[x^i]T_G(1+x,1)$,  $[x^i]T_G(1+x,0)$) counts the orientations (resp. minimal orientations, acyclic orientations) having exactly $i$ edges which are minimal in some head-min directed cocycle. Moreover, by Lemma \ref{lem:ei=min-directed-cocycle}, an orientation has $i$ edges which are minimal in some head-min directed cocycle if and only if it has $i+1$ root-components.\\
\ite We now prove that \emph{the coefficient $[x^i]T_G(x,2)$ (resp. $[x^i]T_G(x,1)$, $[x^i]T_G(x,0)$) counts $v_0$-connected orientations (resp. minimal  $v_0$-connected orientations, acyclic  $v_0$-connected orientations) with $i+1$ root-strong-components.} Let $T$ be a spanning tree with $\mI(T)$ internal $(\mG,T)$-active edges and $\mE(T)$ external $(\mG,T)$-active edges. By Lemma \ref{lem:keylink}, the coefficient  $[x^i]x^{\mI(T)}2^{\mE(T)}$  is the number of connected subgraphs in the tree-interval $[T^-,T^+]$ if $\mI(T)= i$ and $0$ otherwise. Given that the tree-intervals form a partition of the set of subgraphs, the coefficient $[x^i]\sum_{T\textrm{ spanning tree }} x^{\mI(T)}2^{\mE(T)}$ counts the connected subgraphs $S$ such that the tree $T=\Delta(S)$ has $i$ internal $(\mG,T)$-active edges.  Moreover, by the characterization \Ref{eq:Tutte-embedded} of the Tutte polynomial, the sum $\sum_{T} x^{\mI(T)}2^{\mE(T)}$ is equal to $T_G(x,2)$.  Similarly,  the coefficient $[x^i]T_G(x,1)$ (resp. $[x^i]T_G(x,0)$) counts the spanning trees (resp. internal spanning trees) $T$ having $i$  internal $(\mG,T)$-active edges. By Theorem \ref{thm:all-specializations} and Lemma  \ref{lem:mincocycle=active}, the coefficient  $[x^i]T_G(x,2)$  (resp.  $[x^i]T_G(x,1)$,  $[x^i]T_G(x,0)$) counts the  $v_0$-connected orientations (resp. minimal  $v_0$-connected orientations, acyclic  $v_0$-connected orientations) having  exactly $i$ edges which are minimal in some directed cocycle. Moreover, by Lemma \ref{lem:ei=min-directed-cocycle}, an orientation has  $i$ edges which are minimal in some directed cocycle if and only if it has $i+1$ root-strong-components.
\findem

One specialization of this result is of special interest: the coefficient $[x^1]T_G(x,0)$ counts \emph{bipolar orientations}. Given two vertices $u$ and $v$, a \emph{$(u,v)$-bipolar} orientation is an acyclic orientation such that $u$ is the unique source and $v$ is the unique sink. The bipolar orientations are important for many graph algorithms \cite{Mendez:these}. In addition, a bijection between spanning trees having activities $(1,0)$ with respect to Tutte's definition \cite{Tutte:dichromate} and bipolar orientations is the building block used in \cite{Gioan-bij-tree-orientation} in order to define a general correspondence between spanning trees and orientations. This correspondence explains the link between the activities of spanning trees defined by Tutte in \cite{Tutte:dichromate} and the activities of orientations defined by Las Vergnas in \cite{Vergnas:Morphism-matroids-2}.\\

\begin{prop}\label{thm:bipolar}
Let  $\mG$ be an embedded graph, let $v_0$ be the root-vertex and let $v_1$ be the  other endpoint of the root-edge. The mapping $\Phi$ establishes a bijection between the spanning trees having embedding-activities $(\mI(T),\mE(T))=(1,0)$ (counted by $[x^1]T_G(x,0)$) and the $(v_0,v_1)$-bipolar orientations.
\end{prop}

\begin{figure}[htb!]
\begin{center}
\input{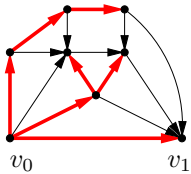}
\caption{A bipolar orientation and the corresponding spanning tree (indicated by thick lines).}\label{fig:bipolar-orientation}
\end{center}
\end{figure}

Proposition \ref{thm:bipolar} is illustrated by  Figure \ref{fig:bipolar-orientation}.\\

\dem Observe first that  an acyclic orientation $\mO$ is $(v_0,v_1)$-bipolar if and only if any vertex is reachable from $v_0$ and can reach $v_1$. By Theorem \ref{thm:refinements} the coefficient  $[x^1]T_G(x,0)$ counts acyclic $v_0$-connected orientation having 2 root-strong-components. No vertex $v\neq v_0$ can reach $v_0$ in an acyclic $v_0$-connected orientation (there would be a directed path from $v_0$ to $v$ and back). Hence the first root-component $U_0$ of an acyclic $v_0$-connected orientation is reduced to $\{v_0\}$. The minimal edge with one endpoint in $U_0=\{v_0\}$ and one endpoint outside $U_0$ is the root-edge. Hence an acyclic $v_0$-connected orientation has 2 root-strong-components if and only if every vertex can reach $v_1$. Thus, the coefficient  $[x^1]T_G(x,0)$ counts $(v_0,v_1)$-bipolar orientations.
\findem

\section{A BIJECTION BETWEEN BETWEEN SPANNING TREES AND RECURRENT SANDPILE CONFIGURATIONS}\label{section:bij-sandpile}

In Section \ref{section:glimpse}, we defined a mapping $\Lambda:T\mapsto \mS_T$ from spanning trees to sandpile configurations. Recall from Definition \ref{def:Lambda} that the number of grains $\mS_T(v)$ on the vertex $v$ in the configuration $\mS_T=\Lambda(T)$ is the number of tails plus the number of external $(\mG,T)$-active heads incident to $v$ in the orientation~$\mO_T=\Phi(T)$. In this section, we prove that the mapping $\Lambda$ is a bijection between spanning trees and recurrent sandpile configurations. 


\begin{thm}\label{thm:bijection-sandpile}
Let $\mG$ be an embedded graph. The mapping $\Lambda:T\mapsto \mS_T$ is a bijection between the spanning trees and the recurrent sandpile configurations of $\mG$.
\end{thm}

Let $G=(V,E)$ be the graph underlying the embedding $\mG$. Observe that the \emph{level} of the configuration $\mS_T$ , that is, $\sum_{v\in V}\mS_T(v)-|E|$, is the number of external $(\mG,T)$-active edges. Indeed, every edge of $G$ has contribution 1 to the sum $\sum_{v}\mS_T(v)$ except the external $(\mG,T)$-active edges which have contribution 2. 

\begin{cor}\label{cor:level-preserving}
Let $\mG$ be an embedded graph.  The number of recurrent  sandpile configurations at level $i$ is the number $[y^i]T_G(1,y)$ of spanning trees having $i$ external $(\mG,T)$-active edges.
\end{cor}

\begin{figure}[htb!]
\begin{center}
\input{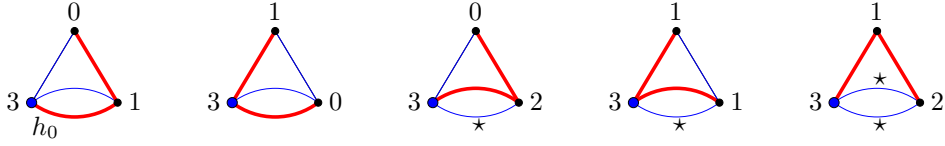}
\caption{The spanning trees (thick lines) and the corresponding  sandpile configurations. The external active edges are indicated by a $\star$.}\label{fig:sandpile-exp}
\end{center}
\end{figure}

As mentioned above, Corollary \ref{cor:level-preserving} is not new. It was first proved recursively in \cite{Merino:external-activity=sandpile-level} and then bijectively in \cite{Borgne-cori-activite-externe-sable} (using Tutte's notion of \emph{activity} \cite{Tutte:dichromate}). The Theorem \ref{thm:bijection-sandpile} and Corollary \ref{cor:level-preserving} are illustrated by Figure \ref{fig:sandpile-exp}.\\

We first prove that the image of any spanning tree is a recurrent sandpile configuration.
\begin{prop}
Let $\mG$ be an embedded graph. For any spanning tree $T$, the  sandpile configuration $\mS_T=\Lambda(T)$ is recurrent. 
\end{prop}


\dem Let $v_0$ be the root-vertex. We consider the orientation $\mO_T$ and prove successively the following properties.\\
\ite \emph{The  sandpile configuration $\mS_T$ is stable.} Let $v$ be any vertex distinct from $v_0$. We want to prove that  $\mS_T(v)< \deg(v)$. Observe that any half-edge incident to $v$ has contribution at most one to $\mS_T(v)$. Moreover, the half-edge $h_v$ incident to $v$ and contained in the edge of $T$ linking $v$ to its father is a head by Lemma \ref{lemma:between-h1-and-h2}. Thus, $h_v$ has no contribution to $\mS_T(v)$, and  $\mS_T(v)\leq \deg(v)-1$.\\
\ite \emph{$\mS_T(v_0)=\deg(v_0)$.} We must prove that every half-edge incident to $v_0$ has contribution 1 to  $\mS_T(v_0)$.  By Lemma \ref{lemma:between-h1-and-h2}, the internal edges are oriented from father to son in $\mO_T$. Therefore any internal half-edge incident to $v_0$ is a tail, hence has contribution 1 to $\mS_T(v_0)$. Let $h$ be an external half-edge incident to $v_0$. By definition, if the half-edge $h$ is greater than the half-edge $h'=\alpha(h)$, then $h$ is a tail. Else, the edge $e=\{h,h'\}$ is  $(\mG,T)$-active by Lemma \ref{lem:external-active} (since the endpoint $v_0$ of $h$ is an ancestor of the endpoint of $h'$). Thus, any external half-edge incident to $v_0$ has contribution 1 to $\mS_T(v_0)$. \\
\ite  \emph{The  sandpile configuration $\mS_T$ is recurrent.} We want to prove that there is a labeling of the vertices $v_0,v_1,\ldots,v_{|V|-1}$ such that the sequence of topplings $\mS_T\topple{v_0}\mS_T^1\topple{v_1}\cdots\topple{v_{|V|-1}}\mS_T^{|V|}$ is valid. Observe that in this case  the configuration $\mS_T$ is recurrent. Indeed, the final configuration $\mS_T^{|V|}$ is equal to $\mS_T$ since every vertex $v$ has been toppled once, hence has sent and received exactly $\deg(v,*)$ grains during the sequence of topplings (recall that $\deg(v,*)$ is the number of non-loop edges incident to $v$). In Section \ref{section:specializations}, we defined a linear order, the \emph{postfix order}, on the vertex set $V$ (see Lemma \ref{lem:toward-greatest=active}). The root-vertex $v_0$ is the maximal element for this order. We want to prove that taking the unique labeling such that  $v_0>v_1>\cdots>v_{|V|-1}$ for the postfix order, the sequence of topplings  $\mS_T\topple{v_0}\mS_T^1\topple{v_1}\cdots\topple{v_{|V|-1}}\mS_T^{|V|}$ is valid. From the preceding point, the toppling of $v_0$ is valid. Suppose that the sequence  $\mS_T\topple{v_0}\mS_T^1\topple{v_1}\cdots \topple{v_{i-1}}\mS_T^{i}$ is valid. After these topplings, the number of  grains on the vertex $v_i$ is $\mS_T^i(v_i)=\mS_T(v_i)+\sum_{j<i}\deg(v_i,v_j)$ (recall that $\deg(v_i,v_j)$ is the number of edges linking $v_i$ and $v_j$). We want to prove that $v_i$ can be toppled, that is, $\mS_T^i(v_i)\geq \deg(v_i)$. By Lemma \ref{lem:toward-greatest=active}, any arc $\mO_T(e)$ is directed toward its least endpoint (for the postfix order) unless $e$ is external $(\mG,T)$-active. Let $h$ be an half-edge in an edge  linking $v_i$  to a vertex $v_j,~j\geq i$. The vertex $v_j$ is less than or equal to $v_i$ for the postfix order, hence $h$ is either a tail or an external $(\mG,T)$-active half-edge. In both cases, the half-edge $h$ has contribution  1 to $\mS_T(v_i)$. Hence,  
$$\mS_T(v_i)\geq \sum_{j\geq i}\deg(v_i,v_j).$$
Thus, 
$$\mS_T^i(v_i)=\mS_T(v_i)+\sum_{j\geq i}\deg(v_i,v_j)\geq \sum_{j\geq 0}\deg(v_i,v_j)=\deg(v_i)$$
and $v_i$ can be toppled. By induction, the sequence of topplings $\mS_T\topple{v_0}\mS_T^1\topple{v_1}\cdots$\\$\cdots\topple{v_{|V|-1}}\mS_T^{|V|}$ is valid.
\findem

It remains to prove that $\Lambda:T\mapsto \mS_T$ is a bijection between the spanning trees and the recurrent sandpile configurations. For this purpose we define a mapping  $\Upsilon$ that we shall prove to be the inverse of $\Lambda$. The mapping $\Upsilon$ is a variant of the \emph{burning algorithm} introduced by Dhar in order to distinguish between recurrent and non-recurrent  sandpile configurations \cite{Dhar:sandpile-self-organized}. The spanning tree returned by the algorithm can be seen as the path through which the \emph{fire} (the sequence of topplings) propagates. The intuitive principle of the algorithm is to decompose each toppling and consider its effect grain after grain. When a grain makes another vertex topple, we add the edge by which the grain has traveled into the tree. Different variants of this algorithm have been proposed \cite{Borgne-cori-activite-externe-sable,Chebikin:family-bijection-sandpile}. These variants differ by the rule used for choosing the next grain to be sent, and also differ from the procedure $\Upsilon$ given below.  Let us insist that the variants considered in \cite{Borgne-cori-activite-externe-sable,Chebikin:family-bijection-sandpile} do not contain our bijection $\Lambda$ as a special case.\\

If $v$ is a vertex and $F\subseteq E$ be a subgraph, we denote by $\deg_F(v)$ the degree of $v$ in the subgraph $F$. 

\begin{Def}\label{def:Upsilon}
Let $\mG=(H,\sigma,\alpha,h_0)$ be an embedded graph. The mapping $\Upsilon$ associates with a recurrent sanpile configuration $\mS$ the spanning tree defined by the following procedure.

\noindent \textbf{Initialization:} Initialize the \emph{current half-edge} $h$ to be $h_0'=\sigma^{-1}(h_0)$. Initialize the tree $T$ and the set of visited edges $F$ to be empty. \\
\textbf{Core:} Do:\\
\textbf{C1:} 
Let $e$ be the edge containing $h$, let $u$ be the vertex incident to $h$ and let $v$ be the other endpoint of $e$.  \\
\indent  If $e$ is not in $F$, then \\
\indent  \indent \iten Add $e$ to $F$.\\
\indent  \indent \iten If $u$ is not connected to $v$ by $T$ and $\mS(v)+\deg_F(v)\geq \deg(v)$ then \\
\indent \indent \indent  \indent Add $e$ to $T$.\\
\textbf{C2:}  Move to the next half-edge clockwise around $T$: \\
\indent  If $e$ is in $T$, then set the current half-edge $h$ to be $\sigma^{-1}\alpha(h)$, else set it to be~$\sigma^{-1}(h)$.\vspace{.0cm}\\
Repeat until the current half-edge $h$ is $h_0'$.\vspace{.0cm}\\
\textbf{End:} Return the tree $T$.
\end{Def}

\begin{figure}[htb!]
\begin{center}
\input{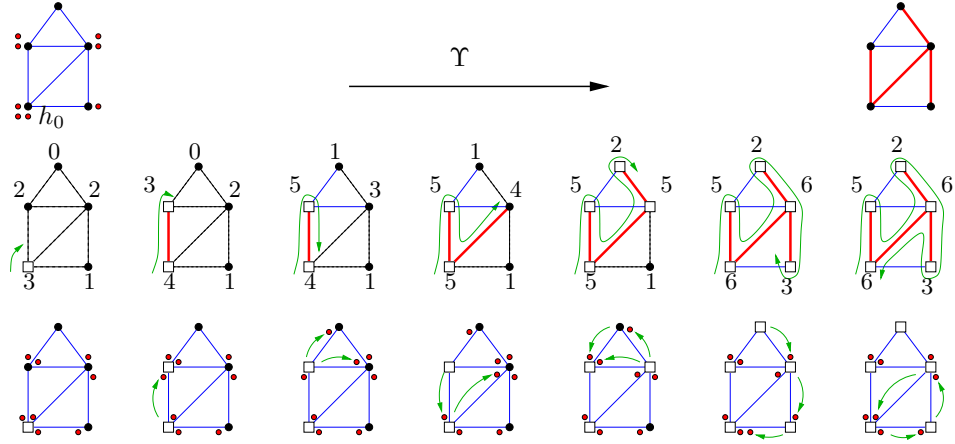}
\caption{The mapping $\Upsilon$. In the middle line, some intermediate steps are represented.  The  set $\B{F}$ of unvisited edges is indicated by dashed lines.  The number associated to each vertex $v$ is equal to $\mS(v)+deg_F(v)$. In the bottom line, the \emph{burning algorithm} representation of each of the intermediate steps is given.}\label{fig:exp-Upsilon}
\end{center}
\end{figure}

We represented the intermediate steps of the procedure $\Upsilon$ in Figure \ref{fig:exp-Upsilon}.\\

Observe that during the procedure $\Upsilon$ our motion (step \textbf{C2}) around the spanning tree is reverse (compared to our previous algorithms). This way of visiting the half-edges would be the usual tour of the spanning tree in the embedded graph $\mG'=(H,\sigma^{-1},\alpha,h_0')$. \\

We will now prove that $\Upsilon$  and $\Lambda$ are inverse bijections. We first prove that the mapping $\Upsilon$ is well defined on recurrent configurations and returns a spanning tree (Proposition \ref{thm:Upsilon-returns-tree}). Then we prove that  $\Upsilon$  and $\Lambda$ are inverse mappings (Propositions \ref{thm:LambdaUpsilon=Id} and \ref{thm:UpsilonLambda=Id}).


\begin{prop}\label{thm:Upsilon-returns-tree}
The procedure $\Upsilon$ is well defined on recurrent configurations and returns a spanning tree. 
\end{prop}

\begin{lemma}\label{lem:Tstay-connected}
Let $\mS$ be a recurrent configuration. Then, at any time of the execution of the procedure $\Upsilon$ on $\mS$, the endpoint $u$ of the current half-edge $h$ is connected to $v_0$ by $T$.
\end{lemma}

\dem
The property holds at the beginning of the execution. Clearly, it remains true each time a step \textbf{C2} is performed. 
\findem

\noindent \textbf{Proof of Proposition \ref{thm:Upsilon-returns-tree}.} Let $\mS$ be a recurrent configuration. We denote by $\Upsilon[\mS]$ the execution of the procedure $\Upsilon$ on $\mS$.
We prove successively the following properties on the execution $\Upsilon[\mS]$.\\
\ite \emph{At any time of the execution, the subgraph $T$ is a tree incident to $v_0$.}
The property holds at the beginning of the execution. Suppose that it holds at the beginning of a given core step and consider the edge $e$ with endpoints $u$ and $v$ containing the current half-edge. If the edge $e$ is added to $T$,  the subgraph $T$ remains acyclic since $u$ is not connected to $v$ by $T$. Moreover the subgraph $T$ remains connected and incident to $v_0$ since (by Lemma \ref{lem:Tstay-connected}) the vertex $u$ is connected to $v_0$ by $T$.\\
\ite \emph{No half-edge is visited twice, hence the execution terminates.}
Suppose that a half-edge $h$ is visited twice during the execution. We consider the first time this situation happens. First note that $h\neq h_0'$ or the execution would have stopped just before the second visit to $h$. Let $h_1$ and $h_2$ be respectively the current half-edge just before the first and second visit to $h$..  Let $T_1$ and $T_2$ be the trees constructed by the procedure $\Upsilon$ at the time of the first and second visit to $h$.  Let $e$ be the edge containing $\sigma^{-1}(h)$. For $i=1,2$ we have $h=\sigma^{-1}\alpha(h_i)$  if $e$ is in $T_i$ and $h=\sigma^{-1}(h_i)$ otherwise. Since $h_1\neq h_2$ and $T_1\subseteq T_2$, the edge  $e$ is in  $T_2$ but not in $T_1$. This is impossible since after the visit of $h_1$ the edge $e$ is in $F$ and cannot be added to the tree $T$ anymore.\\
We denote by $T_0$ the tree returned by the execution $\Upsilon[\mS]$ and by $F_0$ the set of visited edges at the end of this execution. \\
\ite \emph{If $e=\{h_1,h_2\}$ is an edge in $T_0=\Upsilon(\mS)$ and the endpoint of $h_1$ is the father of the endpoint of $h_2$, then $h_1$ is visited during the execution $\Upsilon[\mS]$.} 
Consider the core step at which the edge $e$ is added to the tree $T$. Let  $h$ be the current half-edge, let $u$ be the vertex incident to $h$ and let $v$ be the other endpoint of $e$. By Lemma \ref{lem:Tstay-connected}, the vertex  $u$ is connected to  $v_0$ by $T\subseteq T_0-e$, hence $u$ is the father of $v$. Hence $h_1=h$ is visited  during the execution $\Upsilon$. \\
\ite \emph{At the end of the execution, any edge adjacent to $T_0$ is in $F_0$.} 
We want to show that any half-edge incident to $T_0$ is visited during the execution $\Upsilon[\mS]$.  First observe that no edge can be added to $T$ after its first visit. Therefore, when a step \textbf{C2} is performed, the edge $e$ containing  the current half-edge is in $T$ if and only if it is in $T_0$. Let $h$ be a half-edge incident to $T_0$ which has not been visited during the execution $\Upsilon$. If the half-edge $\sigma^{-1}(h)$ is not in $T_0$ then it has not been visited (or $h$ would have been the next half-edge visited during the execution). Thus by applying $\sigma^{-1}$ repeatedly we find an unvisited half-edge $h$ such that  $\sigma^{-1}(h)$ is in $T_0$. Then, the half-edge $\alpha \sigma^{-1}(h)$ has not been visited during the execution $\Upsilon$ (or $h$ would have been the next half-edge visited during the execution). Thus (by the preceding point) the endpoint of $\alpha \sigma^{-1}(h)$ is the son of the endpoint of $\sigma^{-1}(h)$. We have proved that if there is an unvisited half-edge $h$ incident to $T_0$, then there is an unvisited half-edge incident to one of its sons in $T_0$. We reach an impossibility.\\
\ite \emph{The tree $T_0=\Upsilon(\mS)$ is spanning.}
Let $v_0,v_1,\ldots,v_{|V|-1}$ be a labeling of the vertices such that the sequence $\mS\topple{v_0}\mS^1\topple{v_1}\cdots\topple{v_{|V|-1}}\mS^{|V|}$ is valid. In the configuration $\mS_i$, the number of sand grains on the vertex $v_i$ is $\mS^i(v_i)=\mS(v_i)+\sum_{j<i}\deg(v_j,v_i)$ and is more than the degree of $v_i$. Suppose now that the tree $T_0$ is not spanning and consider the least index $i$ such that $v_i$ is not connected to $v_0$ by $T$.  Each vertex $v_j$ for $j<i$ is incident to $T$, hence (by the preceding point) every edge joining $v_j$ and $v_i$ is in $F_0$. Moreover $v_i$  is adjacent to at least one of the vertices $v_j,j<i$  since $\mS(v_i)$ is less than its degree and $\mS^i(v_i)$ is not. Consider the last edge $e$ (in order of visit) joining $v_i$ to a vertex $v_j,j<i$. When the edge $e$ is visited, we have $\deg_F(v_i)\geq \sum_{j<i} \deg(v_i,v_j)$. Therefore, the condition  $\mS(v_i)+\deg_F(v_i)\geq \deg(v_i)$ holds and the edge $e$ should have been added to the tree $T$. We reach a contradiction.
\findem

We proceed to prove that $\Lambda$ and $\Upsilon$ are inverse mappings.

\begin{lemma}\label{lem:inequality-is-equality}
Consider a given core step of the procedure $\Upsilon$. Let $e$ be the edge containing the current half-edge $h$ and let $v$ be the endpoint of $\alpha(h)$. If the edge $e$ is added to $T$, then the inequality $\mS(v)+\deg_F(v)\geq \deg(v)$ (tested in the procedure $\Upsilon$) is an equality.
\end{lemma}

\dem Observe first that the vertex $v$ is distinct from $v_0$, otherwise adding $e$ to the tree $T$ would create a cycle by Lemma \ref{lem:Tstay-connected}. While $v$ is not connected to $v_0$ by $T$, it is not the endpoint of the current half-edge $h$ (Lemma \ref{lem:Tstay-connected}). Thus, each time the quantity $\deg_F(v)$ increases, that is,  each time an edge incident to $v$ is added to $F$, the condition  $\mS(v)+\deg_F(v)\geq \deg(v)$ is tested and the edge is added to $T$ if the condition holds. 
\findem

\begin{lemma} \label{lem:lambda-in-terms-of-half-edges}
Let $\mG=(H,\sigma,\alpha,h_0)$ be an embedded graph and let $T$ be a spanning tree. We consider the $(\mG,T)$-order on half-edges. Let $v$ be a vertex distinct from $v_0$ and let $h_v$ be the half-edge incident to $v$ in the edge of $T$ linking $v$ to its father. 
Any half-edge $h$ incident to $v$ and such that  $\alpha(h)>h_v$ is external. Moreover, there are  $\deg(v)-\mS_T(v)-1$ such half-edges. 
\end{lemma}

\dem We consider the orientation $\mO_T$. Recall from Lemma \ref{lemma:between-h1-and-h2} that $\alpha(h_v)<h_v$ and that the half-edges $h$ incident to a descendant of $v$ are characterized by  $\alpha(h_v)<h\leq h_v$. In particular, the inequalities $\alpha(h_v)<h\leq h_v$ hold for the half-edges incident to $v$. We now prove successively the following properties. \\
\ite \emph{Any half-edge $h$ incident to $v$ and such that  $\alpha(h)>h_v$ is external.} Suppose that the half-edge $h$ is internal and consider the edge $e$ containing $h$. If $e$ links $v$ to its father, then $h=h_v$ and $\alpha(h)=\alpha(h_v)<h_v$. If $e$ links $v$ to one of its sons, then $\alpha(h)$ is incident to a descendant of $v$ and $\alpha(h)\leq h_v$. In either cases, the hypothesis  $\alpha(h)>h_v$ does not hold.\\
\ite \emph{An external half-edge $h$ incident to $v$ is a non-active head if and only if $\alpha(h)>h_v$.}
The three following properties are sufficient to prove the equivalence:\\
\iten \emph{If $h$ is a tail then $\alpha(h)<h_v$}. Indeed, we have $\alpha(h)<h$ since $h$ is a tail and $h\leq h_v$ since $h$ is incident to  $v$.\\
\iten \emph{If $h$ is a head and $\alpha(h)<h_v$ then  $h$ is $(\mG,T)$-active.} Since $h$ is a head, we have $h<\alpha(h)$ hence, $\alpha(h_v)<h<\alpha(h)<h_v$. Thus, $\alpha(h)$ is incident to a descendant of $v$ and the edge $e=\{h,\alpha(h)\}$ is $(\mG,T)$-active by Lemma \ref{lem:external-active}.\\
\iten \emph{If $h$ is a head and $\alpha(h)>h_v$ then $h$ is not $(\mG,T)$-active.}  Since $h$ is a head we have $h<\alpha(h)$. Since  $\alpha(h)>h_v$, the half-edge $\alpha(h)$ is not incident to a descendant of $v$ and the edge  $e=\{h,\alpha(h)\}$ is not $(\mG,T)$-active by Lemma \ref{lem:external-active}.\\
\ite \emph{There are  $\deg(v)-\mS_T(v)-1$ half-edges $h$ incident to $v$ and such that  $\alpha(h)>h_v$.} By definition, $\mS_T(v)$ is the number of tails plus the number of external $(\mG,T)$-active heads incident to $v$.  Hence, $\deg(v)-\mS_T(v)$  is the number of heads incident to $v$ which are not external  $(\mG,T)$-active. By Lemma \ref{lemma:between-h1-and-h2}, internal edges are oriented from father to son. Hence, the vertex $v$ is incident to exactly one internal head. Thus  $\deg(v)-\mS_T(v)-1$ is the number of external non-active heads. By the preceding point, these half-edges are  characterized by the condition  $\alpha(h)>h_v$.
\findem

We now define the \emph{clockwise-tour} of a tree. Let $\mG=(H,\sigma,\alpha,h_0)$ be an embedded graph. Given a spanning tree $T$, we define the \emph{clockwise-motion function} $\tau$ on half-edges by 
$$\tau(h)=\sigma^{-1}\alpha(h) \textrm{ if } h  \textrm{ is internal and } \tau(h)=\sigma^{-1}(h) \textrm{ otherwise.}$$
As observed above, the clockwise-motion function $\tau$  is the usual motion function for the embedded graph $\mG^{-1}=(H,\sigma^{-1},\alpha,\sigma^{-1}(h_0))$. This defines the \emph{$(\mG^{-1},T)$-order} on the half-edge set $H$ for which $h_0'=\sigma^{-1}(h_0)$ is the least element. The $(\mG,T)$-order denoted by $<$ and the $(\mG^{-1},T)$-order denoted by $<^{-1}$ are closely related.
 
\begin{lemma}\label{lem:equiv-cw-order}
Let $\mG$ be an embedded graph and let $T$ be a spanning tree. The $(\mG,T)$-order and  $(\mG^{-1},T)$-order are related by  $h< h'$ if and only if $\beta(h')<^{-1}\beta(h)$, where $\beta$ is the involution defined by $\beta(h)=h$ if  $h$ is external and  $\beta(h)=\alpha(h)$  otherwise.
\end{lemma}

\dem
Let $t$ be the usual motion function and let $\tau$ be the clockwise-motion function. Observe that  $t \beta= \sigma$ and $\tau \beta= \sigma^{-1}$. Thus, $\tau=\beta t^{-1} \beta$. Let us write $t=(h_0,h_1,\ldots,h_{|H|-1})$ in cyclic notation. Then $t^{-1}=(h_{|H|-1},\ldots,h_1,h_{0})$ and   $\tau=\beta t^{-1} \beta=(\beta(h_{|H|-1}),\ldots,\beta(h_1),\beta(h_0))$. Moreover, $\sigma\beta(h_{|H|-1})=t(h_{|H|-1})=h_0$, hence $\beta(h_{|H|-1})=h_0'=\sigma^{-1}(h_0)$. Therefore, $h_i< h_j$ if and only if $i<j$ if and only if $\beta(h_j)<^{-1}\beta(h_i)$.
\findem

\begin{lemma}\label{lem:visited-in-cw-order}
Let $\mS$ be a recurrent configuration and let $T_0=\Upsilon(\mS)$ be the spanning tree returned by the procedure $\Upsilon$. The  half-edges of $\mG$ are visited in $(\mG^{-1},T_0)$-order during the procedure $\Upsilon$.
\end{lemma}

\dem During the procedure $\Upsilon$, no edge can be added to the tree $T$ after its first visit. Therefore, when a step \textbf{C2} is applied, the edge $e$ containing  the current half-edge is in $T$ if and only if it is in $T_0$. Hence, a step \textbf{C2} corresponds to an application of the clockwise-motion function $\tau$ of the spanning tree $T_0$. Since the first visited half-edge is $h_0'=\sigma^{-1}(h_0)$, the half-edges are visited in $(\mG^{-1},T_0)$-order.
\findem

\begin{lemma} \label{lem:number-visited}
Let $\mG$ be an embedded graph and let $T$ be a spanning tree. Let $v$ be a vertex distinct from $v_0$ and let $e_v$ be the edge of $T$ linking $v$ to its father. There are  $\deg(v)-\mS_T(v)-1$ edges incident to $v$ and less than $e_v$ for the $(\mG^{-1},T)$-order.
\end{lemma}

\dem
Let $h_v$ be  the half-edge of $e_v$ incident to $v$. Let $h\neq h_v$ be a half-edge incident to $v$ and let $e$ be the edge containing $h$. We prove successively the following properties. \\
\ite \emph{The edge $e$ is less than $e_v$ if and only if $\alpha(h)<^{-1} \alpha(h_v)$. Moreover, in this case $e$ is not a loop.} 
By Lemma \ref{lemma:between-h1-and-h2} applied to the embedded graph $\mG^{-1}$, the half-edges $h$ incident to $v$ are such that $\alpha(h_v)<^{-1}h\leq^{-1}h_v$. Hence, the edge containing $h$ is less than $e_v$ for the $(\mG^{-1},T)$-order if and only if $\alpha(h)<^{-1} \alpha(h_v)$. In this case, $\alpha(h)$ is not incident to $v$ by Lemma \ref{lemma:between-h1-and-h2}, that is, $e$ is not a loop. \\
\ite \emph{The conditions  $\alpha(h)<^{-1} \alpha(h_v)$ and $\alpha(h)> h_v$ are equivalent. Moreover,  there are $\deg(v)-\mS_T(v)-1$ half-edges satisfying this condition.} 
Suppose $\alpha(h)<^{-1} \alpha(h_v)$. In this case, $h$ external. Indeed, $h$ is not in $e_v$ and is not incident to a son of  $v$ by Lemma \ref{lemma:between-h1-and-h2} applied to the embedded graph $\mG^{-1}$. Hence, by Lemma \ref{lem:equiv-cw-order}, we get $\alpha(h)> h_v$.  Conversely, if  $\alpha(h)> h_v$, the edge $e$ is external by Lemma \ref{lem:lambda-in-terms-of-half-edges}, hence $\alpha(h)<^{-1} \alpha(h_v)$ by Lemma \ref{lem:equiv-cw-order}. Moreover, there are $\deg(v)-\mS_T(v)-1$ half-edges satisfying this condition by Lemma   \ref{lem:lambda-in-terms-of-half-edges}.
\findem

\begin{prop}\label{thm:LambdaUpsilon=Id}
The mapping $\Lambda\circ\Upsilon$ is the identity on recurrent configurations.
\end{prop}

\dem
Let $\mS$ be a recurrent configuration and let $T=\Upsilon(\mS)$. We want to prove that the recurrent configuration $\mS_T=\Lambda(T)$ is equal to $\mS$. We already know that $\mS_T(v_0)=\deg(v_0)=\mS(v_0)$ since  $\mS_T$ and $\mS$ are recurrent configurations. Let $v$ be a vertex distinct from $v_0$ and let $e_v$ be the edge of $T$ linking $v$ to its father. Let $F$ be the set of visited edges when $e_v$ is added to $T$ during the execution $\Upsilon[\mS]$. We know that $\mS(v)=\deg(v)-\deg_F(v)$ by Lemma \ref{lem:inequality-is-equality}. It remains to prove that $\mS_T(v)=\deg(v)-\deg_F(v)$. By Lemma \ref{lem:visited-in-cw-order}, the half-edges are visited in $(\mG^{-1},T)$-order during the execution $\Upsilon[\mS]$. Therefore, the value $\deg_F(v)$ is the number of edges incident to $v$ which are less or equal to $e_v$ for the $(\mG^{-1},T)$-order. There are $\deg(v)-\mS_T(v)$ such edges by Lemma \ref{lem:number-visited}. We obtain $\deg_F(v)=\deg(v)-\mS_T(v)$, or equivalently, $\mS_T(v)=\deg(v)-\deg_F(v)$. Thus, $\mS_T(v)=\mS(v)$.
\findem

\begin{prop}\label{thm:UpsilonLambda=Id}
The mapping $\Upsilon\circ \Lambda$ is the identity on spanning trees.
\end{prop}

\dem
Let $T_0$ be a spanning tree. We denote by $T_1=\Upsilon(\mS_{T_0})$ the image of $T_0$ by $\Upsilon\circ \Lambda$ and want to prove that $T_1=T_0$. Recall that every edge of $\mG$ is visited during the execution $\Upsilon[\mS_{T_0}]$. Hence, it is sufficient to prove that at the beginning of any core step of the execution $\Upsilon[\mS_{T_0}]$, \emph{the tree $T$ constructed by the procedure $\Upsilon$ is $T_0\cap F$, where $F$ denotes the set of visited edges.} We proceed by induction on the number of core steps. The property holds at the beginning of the first core step. Suppose that it holds at the beginning of the $k^{th}$ core step. If the edge $e$ containing the current half-edge is already in the set $F$ of visited edges, then the set $F$ and the tree $T$ are unchanged during this core step and the property holds at the beginning of the $k+1^{th}$ core step. Suppose now that the edge $e$ is not in $F$ at the beginning of the $k^{th}$ core step. By the induction hypothesis, the tree $T$ constructed by the procedure $\Upsilon$ is $T_0\cap F$. Moreover, no edge is added to the tree $T$ after its first visit, hence  $T=T_1\cap F$. In other words, the spanning trees $T_0$ and $T_1$ coincide on $F$. By Lemma \ref{lem:visited-in-cw-order}, the half-edges are visited in $(\mG^{-1},T_1)$-order during the execution $\Upsilon[\mS_{T_0}]$, hence the edges visited before $e$ during the execution $\Upsilon[\mS_{T_0}]$ have been visited in $(\mG^{-1},T_0)$-order. Thus, the edges visited before $e$ during the execution $\Upsilon[\mS_{T_0}]$ are the edges which are less than $e$ for the  $(\mG^{-1},T_0)$-order. Suppose now that the edge $e$ is in the tree $T_0$. In this case the endpoints $u$ and  $v$ of $e$ are not connected by $T\subseteq T_0-e$.  Moreover, the value  $\deg_{F+e}(v)$ which corresponds to the number of edges incident to $v$ and visited before $e$ during the execution $\Upsilon[\mS_{T_0}]$, that is, the edge which are less or equal to $e$ for the  $(\mG^{-1},T_0)$-order, is $\deg(v)-\mS_{T_0}(v)$ by Lemma \ref{lem:number-visited}. Thus, the condition $\mS_{T_0}(v)+\deg_{F+e}(v)\geq \deg(v)$ (tested by the procedure $\Upsilon$) holds and the edge  $e$ is added to the tree $T$. Suppose now that $e$ is not in $T_0$. In this case, the edge $e_v$ linking $v$ to its father in $T_0$ is greater than $e$ for the $(\mG^{-1},T_0)$-order. Hence, the value $\deg_{F+e}(v)$ is less or equal to the number of edges incident to $v$ which are less than $e_v$ for the  $(\mG^{-1},T_0)$-order. Thus, $\deg_{F+e}(v)< \deg(v)-\mS_{T_0}(v)-1$ by Lemma \ref{lem:number-visited}. The condition $\mS_{T_0}(v)+\deg_{F+e}(v)\geq \deg(v)$ (tested by the procedure $\Upsilon$) does not hold, hence the edge $e$ is not added to the tree $T$. In any case, the property holds at the beginning of the $k+1^{th}$ core step. 
\findem

This concludes our proof of Theorem \ref{thm:bijection-sandpile}. 
\findem

\section{CONCLUDING REMARKS}\label{section:conclusion}
\subsection{The cycle and cocycle reversing systems} \label{section:cycle-reversing-system}
We consider the \emph{cycle reversing system} and the \emph{cocycle reversing system}. A transition in the cycle (resp. cocycle) reversing system consists in flipping a directed cycle (resp. cocycle).  The cycle and cocycle reversing systems appear implicitly in many works (e.g. \cite{Felsner:lattice,Fraysseix:Topological-aspect-orientations,Propp:lattice,Bonichon:realizers}).\\

It is known from  \cite{Propp:lattice} that there is a unique $v_0$-connected orientation (equivalently, orientation without head-min directed cocycle by Lemma \ref{lem:no-head-min=root-connected}) in each equivalence class of the cocycle reversing system.  The counterpart of this property for the cycle reversing system is given by Proposition \ref{thm:unique-minimal}. Indeed, the equivalence classes of the cycle reversing system are in one-to-one correspondence with outdegree sequences \cite{Felsner:lattice}. Thus, Proposition \ref{thm:unique-minimal} proves that there is a unique minimal orientation (that is, orientation without tail-min directed cycle) in each equivalence class of the cycle reversing system.\\

The \emph{cycle-cocycle reversing system} in which a transition consists in flipping  either a directed cycle or a directed cocycle was introduced in \cite{Gioan:enumerating-degree-sequences}. It was observed in this paper that the cycle and cocycle flips are really independent since they act on the cyclic part and acyclic part  respectively and do not modify the other part. As a consequence it was shown that the   equivalence classes of the cycle-cocycle reversing system are in one-to-one correspondence with root-connected outdegree sequences. Since the  cycle and cocycle flips are independent, the  unicity of the $v_0$-connected orientation in the classes of the cocycle reversing system (\cite{Propp:lattice}) and the unicity of minimal orientation in the classes of the cycle reversing system (Proposition \ref{thm:unique-minimal}) proves that there is a unique $v_0$-connected minimal orientation in each equivalence class of the cycle-cocycle reversing system.\\

As observed in \cite{Gioan:enumerating-degree-sequences}, the enumerative results of Theorem \ref{thm:all-specializations} can be expressed in terms of cycle/cocycle reversing systems. For instance, the  equivalence classes of the cocycle reversing system  (in bijection with minimal orientations) are counted by $T_G(1,2)$, the equivalence classes of the cocycle reversing system reduced to one element (equivalently, the strongly connected orientations) are counted by $T_G(0,2)$ \emph{etc}.

\subsection{The planar case and duality} 
In this subsection we restrict our attention to planar graphs. Our goal is to highlight some nice properties of our bijections with respect to duality. Therefore we will handle simultaneously a planar embedding and its \emph{dual}. In order to avoid confusion we shall indicate the implicit embedding $\mG$ for the tree-intervals and the mapping $\Phi$ by writing $[T^-,T^+]_\mG$ and $\Phi_\mG$. \\

Let $G=(V,E)$ be a planar graph. The graph $G$ can be embedded in the sphere, that is, drawn in such a way the edges only intersect at their endpoints. An embedding of $G$ in the oriented sphere  defines a combinatorial  embedding $\mG=(H,\sigma,\alpha)$ where the permutation $\sigma$ corresponds to the counterclockwise order around each vertex. There is a one-to-one  correspondence between the  embedding of graphs in the oriented sphere and  combinatorial embeddings having \emph{Euler characteristic} 0, where the Euler characteristic is the number of vertices (cycles of  $\sigma$) plus the number of \emph{faces} (cycles of  $\sigma\alpha$) minus the number of \emph{edges} (cycles of  $\alpha$) minus 2. We call these embeddings planar. If $\mG=(H,\sigma,\alpha,h_0)$ is a (combinatorial) planar embedding, then $\mG^*=(H,\sigma\alpha ,\alpha,h_0)$ correspond to the graphical dual of $\mG$ in the reverse-oriented sphere (the graphical dual of a graph embedded in the sphere is obtained by putting a vertex in each face and an edge across each edge). Observe, by the way that $\mG^{**}=\mG$. \\

Consider a planar embedding $\mG$. Observe that the edges, subgraphs and orientations of $\mG$ can also be considered as edges, subgraphs and orientations of $\mG^*$. Given a subgraph $S$ of $\mG$ we denote by $\B{S}^*$ the \emph{co-subgraph}, that is, the complement of $S$ considered as a subgraph of $\mG^*$. Given an orientation $\mO$ of $\mG$ we denote by $\B{\mO}^*$ the \emph{co-orientation}, that is, the orientation obtained from $\mO$ by reversing all arcs considered as an orientation of $\mG^*$. Observe that for any subgraph $S$ and any orientation $\mO$, we have $\B{\B{S}^*}^*=S$ and $\B{\B{\mO}^*}^*=\mO$. From the Jordan Lemma, a subgraph $S$ is connected if and only if the co-subgraph $\B{S}^*$ is acyclic. This implies the well known property  (see \cite{Mullin:tree-rooted-maps}) that a subgraph $T$ is a spanning tree of $\mG$ if and only if the co-subgraph $\B{T}^*$ is a spanning tree of $\mG^*$.  From this property, it follows that the fundamental cycle (resp. cocycle) of an internal (resp. external) edge $e$ with respect to $\mG$ and $T$ is the fundamental cocycle (resp. cycle) of  $e$ with respect to $\mG^*$ and $\B{T}^*$. Moreover, it follows directly from the definitions that the motion function of the spanning tree $T$ of $\mG$ and the motion function of the spanning tree $\B{T}^*$ of $\mG^*$ are equal. In particular, the $(\mG,T)$-order and the $(\mG^*,\B{T}^*)$-order are the same. Hence, an edge is $(\mG,T)$-active if and only if it is $(\mG^*,\B{T}^*)$-active. Thus, the mapping $S\mapsto \B{S}^*$ induces a bijection between the tree-intervals  $[T^-,T^+]_\mG$ and  $[{\B{T}^*}^-,{\B{T}^*}^+]_{\mG^*}$. It follows directly from this property and the definitions that the mappings $\Phi_\mG$ and $\Phi_{\mG^*}$ are related by :
$$\textrm{for any subgraph } S \textrm{ of }\mG,~~\B{\Phi_\mG(S)}^*=\Phi_{\mG^*}(\B{S}^*).$$

\vspace{.3cm}

\noindent \textbf{Acknowledgments:} This work has benefited from discussions with \'Eric Fusy, Emeric Gioan,  Yvan Le Borgne, Igor Pak, Gilles Schaeffer and Michel Las Vergnas. I am deeply indebted to Mireille Bousquet-Mélou for her very patient reading of early versions of this paper and for her constant support and guidance.

\bibliography{../../../biblio/allref}

\bibliographystyle{plain}

\end{document}